\documentclass[11pt,a4paper]{article}
\usepackage{amsmath}
\usepackage{amssymb}
\usepackage{amsthm}
\usepackage{verbatim}
\usepackage[all]{xy}

\newtheorem{prop}{Proposition}[section]

\newtheorem{theorem}{Theorem}
{\theoremstyle{definition}
\newtheorem{rem}[prop]{Remark}}
\newtheorem{thm}[prop]{Theorem}
\newtheorem{lem}[prop]{Lemma}
\newtheorem{de}[prop]{Definition}
\newtheorem{cor}[prop]{Corollary}
\newtheorem{corollary}[theorem]{Corollary}
{\theoremstyle{definition}
\newtheorem{se}[prop]{}}

\numberwithin{equation}{prop}

\begin{document}

\def \A{\mathbb{A}}
\def \calA{\mathcal{A}}
\def \B{\mathcal{B}}
\def \C{\mathcal{C}}
\def \CC{\mathbb{C}}
\def \D{\mathcal{D}}
\def \F{\mathbb{F}}
\def \bF{\overline{\F}}
\def \calF{\mathcal{F}}
\def \G{\mathcal{G}}
\def \bG{\mathbb{G}}
\def \H{\mathcal{H}}
\def \hH{\operatorname{H}}
\def \I{\mathcal{I}}
\def \J{\mathcal{J}}
\def \L{\mathcal{L}}
\def \bM{\mathbb{M}}  
\def \calM{\mathcal{M}}
\def \M{{\boldsymbol M}}
\def \tM{\tilde{\boldsymbol M}}
\def \hM{\hat{\boldsymbol M}}
\def \N{\mathcal{N}}
\def \O{\mathcal{O}}
\def \P{\mathbb{P}}
\def \calP{\mathcal{P}}
\def \Q{\mathbb{Q}}
\def \R{\mathbb{R}}
\def \calS{\mathcal{S}}
\def \T{\mathcal{T}}
\def \oT{\overline{\mathcal{T}}}
\def \U{\mathbb{U}}
\def \calU{\mathcal{U}}
\def \bU{\boldsymbol U}
\def \V{\mathcal{V}}
\def \X{\mathbb{X}}
\def \Z{\mathbb{Z}}
\def \e{\tilde{e}}
\def \oe{\overline{e}}
\def \f{\tilde{f}}
\def \lmi{{\lambda_i,\mu_i}}

\setlength{\unitlength}{0.7mm}
\def \a{ \begin{picture}(5,5)
\put(1,0){\circle*{0.5}}
\put(3,2){\circle*{0.5}}
\put(5,4){\circle*{0.5}}
\end{picture}}
\def \charpol{\operatorname{charpol}}
\def \Fr{\operatorname{Fr}}
\def \Res{\operatorname{Res}}
\def \id{\operatorname{id}}
\def \End{\operatorname{End}}
\def \Hom{\operatorname{Hom}}
\def \GL{\operatorname{GL}}
\def \SL{\operatorname{SL}}
\def \Lie{\operatorname{Lie}}
\def \Grass{\operatorname{Grass}}
\def \Spec{\operatorname{Spec}}
\def \Proj{\operatorname{Proj}}
\def \inv{\operatorname{inv}}
\def \diag{\operatorname{diag}}
\def \det{\operatorname{det}}
\def \ord{\operatorname{ord}}
\def \vol{\operatorname{vol}}
\def \mod{\operatorname{mod}}
\def \type{\operatorname{type}}
\def \length{\operatorname{length}}
\def \dim{\operatorname{dim}}
\def \rk{\operatorname{rk}}
\def \inv{\operatorname{inv}}
\def \GSp{\operatorname{GSp}}
\def \GU{\operatorname{GU}}
\def \d{^\vee}
\def \gr{\operatorname{gr}}
\def \height{\operatorname{ht}}
\def \Stab{\operatorname{Stab}}
\def \Qisg{\operatorname{Qisg}}

\title{The supersingular locus of the Shimura variety for $\GU(1,s)$}
\author{Inken Vollaard}
\maketitle

\begin{abstract}
In this paper we study the supersingular locus of the reduction modulo
$p$ of the Shimura variety for $\GU(1,s)$ in the case of an inert prime
$p$. Using Dieudonn\'e theory we define  a stratification of the
corresponding moduli 
space of $p$-divisible groups. We describe the incidence relation of
this stratification in terms of the Bruhat-Tits building of a unitary
group. 

In the case of $\GU(1,2)$, we show that the supersingular locus is
equi-dimensional of dimension 1 and is of complete intersection. We give
an explicit 
description of the irreducible components and their intersection behaviour.\\

\end{abstract}

\section*{Introduction}
In this paper we study the 
supersingular locus of the reduction modulo $p$ of the Shimura variety for
$\GU(1,s)$ in the case of an inert prime $p$. For
$\GU(1,2)$  
this is a purely 1-dimensional variety, and
we describe explicitly the irreducible components and their
intersection behaviour.

The results of this paper are thus part of the general program of
giving an explicit description of the supersingular (or, more
generally, the basic) locus of the
reduction modulo $p$ of a Shimura variety. Let us review previous work
on this problem.

We fix a prime $p$ and
denote by $\calA_g$ the moduli space over $\bF_p$ of principally
polarized abelian varieties 
of dimension 
$g>0$ with a small enough level structure prime to $p$. Let
$\calA_g^{ss}$ be its 
supersingular locus.  
If $g=1$, this is a finite set of points. 
If $g=2$, Koblitz (\cite{Kb}) shows
that the irreducible components of
$\calA_2^{ss}$ are smooth curves which intersect pairwise transversally
at the superspecial points, i.e., the points of $\calA_2$ where the
underlying abelian variety is superspecial (i.e., is isomorphic to a power 
of a supersingular elliptic curve). Each superspecial point is
the intersection of $p+1$ irreducible components. 

Katsura and Oort
prove in \cite{KO1} that  each irreducible
component of $\calA_2^{ss}$ is isomorphic to $\P^1$. In
\cite{KO2} they calculate the 
dimension of the irreducible
components of $\calA_3^{ss}$ and the
number of irreducible components of $\calA_2^{ss}$ and $\calA_3^{ss}$. 
In the case $g=3$, Li and Oort show in \cite{LO} that the irreducible
components are birationally equivalent to  a
$\P^1$-bundle over a Fermat curve. Furthermore, they
compute for general $g$ the dimension of the supersingular locus and
the number of irreducible components.  

The $\bF_p$-rational points of $\calA_2^{ss}$  are
described independently by
Kaiser (\cite{Ka}) and by Kudla and Rapoport (\cite{KR2}).
They fix a supersingular principally polarized
abelian variety $A$ over $\bF_p$ of dimension 2  and
study the $\bF_p$-rational points 
of the moduli space of quasi-isogenies of $A$. They
cover the $\bF_p$-rational points of these moduli spaces with subsets
which are  
in bijection with the $\bF_p$-rational points of $\P^1$.
The
incidence relation of these subsets is described by the
Bruhat-Tits building of an algebraic group  over $\Q_p$. 
The number of superspecial points of $\calA_2^{ss}$ is calculated in
\cite{KR2}. Each
irreducible component of $\calA_2^{ss}$ contains $p^2+1$ superspecial points.

The $\bF_p$-rational points of the moduli space of quasi-isogenies of
a supersingular principally polarized abelian variety of dimension 3
are described by 
Richartz (\cite{Ri}). In analogy to the case $g=2$, she defines subsets
of the $\bF_p$-rational points of this moduli space and proves that
the incidence relation of these sets is given by the Bruhat-Tits
building of an algebraic group over $\Q_p$. She identifies some of
these sets with the $\bF_p$-rational points of Fermat curves over $\bF_p$. 

Now consider the supersingular locus of a Hilbert-Blumenthal variety
associated to a totally real field extension of degree $g$ of $\Q$ 
in the case of an inert prime $p$. 
For $g=2$ the supersingular locus 
is studied by Stamm (\cite{St}, comp.~\cite{KR1}). He shows that the
irreducible 
components of the supersingular
locus are isomorphic to $\P^1$ and  contain $p^2+1$ superspecial points. 
Two
components intersect transversally in at most one superspecial point
and each superspecial point is the intersection of two irreducible components.
The number of
irreducible components and the number of superspecial points of the
supersingular 
locus are calculated in \cite{KR1}.

Bachmat and Goren analyse the case $g=2$ for arbitrary prime $p$ (\cite{BG}). 
They give another proof of Stamm's results and prove that in case of a 
ramified prime $p$ the components of the supersingular locus are smooth 
rational curves. Furthermore, they calculate the number of components in 
case of an inert or ramified prime $p$. 

The case of $g=3$ and of an inert prime $p$ is studied by Goren in \cite{Gn}. 
He proves that the irreducible components are smooth rational curves. Each 
superspecial point is the intersection of three irreducible components and 
each irreducible component contains only one superspecial point.

Yu analyses in \cite{Yu} the case of $g=4$. He computes the number of
irreducible components of the supersingular locus and the completion
of the local ring at every superspecial point. Furthermore, he shows
that every irreducible component is isomorphic to
a ruled surface over $\P^1$. 

Goren and Oort analyze in \cite{GO} the Ekedahl-Oort stratification
for general $g$. They
prove that the supersingular locus is equi-dimensional of
dimension $[g/2]$.

Finally, consider the
supersingular locus of the reduction
modulo $p$ of the Shimura 
variety for $\GU(1,s)$.
In case of an inert prime $p$, 
B\"ultel and Wedhorn prove that the dimension of the supersingular locus is
equal to 
$[s/2]$  (\cite{BW}).\\

We will now recall the definition of the moduli space of abelian
varieties for $\GU(r,s)$. We first review the corresponding PEL-data,
comp.~\cite{Ko}. 
Let $E$ be an imaginary quadratic extension of $\Q$  such that $p$ is inert in
$E$ and let $\O_E$ 
be its ring of integers. 
Denote by $^*$ the nontrivial
Galois automorphism of $E$. 
Let $V$ be an $E$-vector space of dimension
$n>0$ with perfect alternating $\Q$-bilinear form 
$$(\cdot,\cdot):V\times V\rightarrow \Q$$
such that 
$$(xv,w)=(v,x^*w)$$
for all $x\in E$ and $v,w\in V$. 
Denote by $G$ the algebraic group
over $\Q$ such that
\begin{align*}
G(R)=\{g\in\GL_{E\otimes_\Q R}(V\otimes_\Q R)\mid (gv,gw)=c(g)(v,w);\
c(g)\in R^\times\}
\end{align*}
for every $\Q$-algebra $R$. Let $E_p$ be the completion of $E$ with
respect to the $p$-adic topology.
We assume that there exists an $\O_{E_p}$-lattice $\Lambda$ in
$V\otimes_\Q \Q_p$ such that the form $(\cdot,\cdot)$ induces a perfect
$\Z_p$-form on $\Lambda$. 
We fix an embedding $\varphi_0$ of $E$ into $\CC$ whereby identifying $E\otimes_\Q \R$ with $\CC$. 
We assume that there exists an isomorphism of $\CC$-vector spaces of $V\otimes_\Q\R$ with $\CC^n$
such that the form $(\cdot,\cdot)$ induces the hermitian form  given by the 
diagonal matrix $\diag(1^r,(-1)^s)$ on $\CC^n$. We fix such an
isomorphism. Then 
the group $G_\R$ is equal to the group $\GU(r,s)$
of unitary similitudes of $\diag(1^r,(-1)^s)$. The nonnegative
integers $r,s$ 
satisfy $r+s=n$. 
Let $h$ be the homomorphism of
real algebraic groups
$$h:\Res_{\CC/\R}(\bG_{m,\CC})\rightarrow G_\R$$
which maps an element $z\in\CC^\times$ to the matrix
$\diag(z^r,\overline{z}^s)$.
Then $(E,\O_E,^*,V,(\cdot,\cdot),\Lambda,G,h)$ is a PEL-datum.  
Let $K$ be the reflex field of this PEL-datum. Then $K$ is isomorphic 
to $E$ if $r\neq s$ and is equal to $\Q$ if
$r=s$. Denote by 
$K_p$ the completion of $K$ with respect to the $p$-adic topology and
by $\F$ the residue field of $K$.

Let $A$ be an abelian scheme over an $\O_{K_p}$-scheme $S$ of dimension
$n$
with
$\O_E$-action, i.e., with a morphism
$$\iota:\O_E\rightarrow\End A.$$
Let $\varphi_0$ and $\varphi_1$ be the
two $\Q$-embeddings of $E$ into $\overline{K}_p$. Then the polynomial 
$(T-\varphi_0(a))^r(T-\varphi_1(a))^s$ is an element of $\O_{K_p}[T]$.
We say that $(A,\iota)$ satisfies the Kottwitz determinant condition
of signature $(r,s)$
(\cite{Ko} Chap.~5) if 
\begin{align*}
\charpol(a,\Lie A)=(T-\varphi_0(a))^r(T-\varphi_1(a))^s\in
\O_S[T]  
\end{align*}
for all $a\in\O_E$. 

We recall the definition of the moduli problem defined by Kottwitz in
\cite{Ko} for these PEL-data. 
Denote by $\A_f^p$ the ring of finite
adeles of $\Q$ away from $p$.
We fix a compact open
subgroup $C^p$ of $G(\A_f^p)$. We denote by $A\d$ the dual abelian
scheme of $A$.
Let $\calM$ be the moduli problem which associates to any 
$\O_{K_p}$-scheme $S$ the isomorphism classes of the following data. 
\begin{itemize}
\item An abelian scheme $A$ over $S$ of dimension $n$.
\item A $\Q$-subspace  $\overline{\lambda}$ of
  $\Hom(A,A\d)\otimes_\Z\Q$ such that $\overline{\lambda}$ contains a
  $p$-principal polarization (i.e. a polarization of order prime to $p$). 
\item A homomorphism $\iota\colon\O_E\rightarrow\End(A)\otimes_\Z\Z_{(p)}$ 
  such that the
  Rosati involution given by 
  $\overline{\lambda}$ on $\End(A)\otimes_\Z\Z_{(p)}$ induces
  the involution $^*$ on $E$.
\item A  $C^p$-level structure
  $\overline{\eta}:\hH_1(A,\A_f^p)\stackrel{\sim}{\longrightarrow}
  V\otimes_\Q\A_f^p\mod C^p$. 
\end{itemize}
We assume that 
  $(A,\iota)$ satisfies the
  determinant condition of signature $(r,s)$.

Two such data $(A,\overline{\lambda},\iota, \overline{\eta})$ and 
$(A',\overline{\lambda}',\iota', \overline{\eta}')$ are isomorphic 
if there exists an isogeny prime to $p$ from $A$ to $A'$, commuting with the 
action of $\O_E$, carrying $\overline{\eta}$ into $\overline{\eta}'$ 
and carrying $\overline{\lambda}$ into $\overline{\lambda}'$. For simplicity we write $(A\otimes_\Z\Z_{(p)},\overline{\lambda},\iota\otimes_\Z \Z_{(p)}, \overline{\eta})$ for such an isomorphism class. \\

The moduli problem $\calM$ is represented by a smooth, quasi-projective scheme 
over $\Spec\O_{K_p}$ if $C^p$ is small enough (\cite{Ko} Chap.~5). The
relative dimension of $\calM$ is 
equal to $rs$. 
Denote by $\calM^{ss}$ the supersingular locus of the
special fibre $\calM_\F$ of $\calM$. It is a closed subscheme of
$\calM_\F$ which is proper over $\Spec\F$. 
Our goal is to describe the
irreducible components of $\calM^{ss}$ and their 
intersection behaviour. 

The supersingular locus
$\calM^{ss}$ contains an $\bF_p$-rational point (\cite{BW}
Lem.~5.2). We will view $\calM^{ss}$ as a scheme over $\bF_p$.
Let
$x=(A\otimes_\Z\Z_{(p)},\iota\otimes_\Z\Z_{(p)},\overline{\lambda},
\overline{\eta})\in\calM^{ss}(\bF_p)$ and 
denote by $(\X,\iota)$ the supersingular $p$-divisible group
of height $2n$ 
corresponding
to $x$ with $\O_{E_p}$-action $\iota$. We choose a $p$-principal
polarization $\lambda\in\overline{\lambda}$ and denote again by
$\lambda$ the induced $p$-principal polarization of $\X$.
By construction $(\X,\iota)$ 
satisfies the determinant condition of type $(r,s)$.

We recall the definition of the moduli space $\N$ of quasi-isogenies
of $p$-divisible groups in characteristic $p$ (\cite{RZ} Def.~3.21) in 
the case of the group $\GU(r,s)$.
The moduli space $\N$ over $\Spec \overline{\F}_p$ is given by
the following data up to isomorphism for an $\overline{\F}_p$-scheme
$S$.
\begin{itemize}
\item A $p$-divisible group $X$ over $S$ of height $2n$ with
   $p$-principal polarization
  $\lambda_X$ and $\O_{E_p}$-action $\iota_X$ such that the Rosati
   involution induced by 
   $\lambda$ induces the involution $^*$ on $\O_E$. We assume that 
   $(X,\iota)$ satisfies the determinant condition of type $(r,s)$.
\item An $\O_{E_p}$-linear quasi-isogeny
  $\rho:X\rightarrow\X\times_{\Spec\overline{\F}_p}S$ such that
  $\rho\d\circ\lambda\circ\rho$ is a $\Q_p$-multiple of $\lambda_X$
  in $\Hom_{\O_{E_p}}(X,X\d)\otimes_\Z\Q$. 
\end{itemize}

The moduli space $\N$ is represented by a separated formal scheme which is
locally formally of finite type over $\overline{\F}_p$ (\cite{RZ}
Thm.~3.25).

We recall the uniformization theorem of Rapoport and Zink (\cite{RZ}
Thm.~6.30). We will formulate this theorem only for the underlying
schemes, not for the formal schemes.
Let $I(\Q)$ be the group of quasi-isogenies in
$\End_{\O_E}(A)\otimes\Q$ which respect the homogeneous polarization
$\overline{\lambda}$. 
As $GU(r,s)$ satisfies the Hasse principle (cf. \ref{a1}), 
there exists an isomorphism of schemes over $\bF_p$
\begin{align*}
I(\Q)\backslash(\N^{red}\times G(\A_f^p)/C^p)\stackrel{\sim}{\longrightarrow}
\calM^{ss}.
\end{align*}
In Section~\ref{a} we will show, in the case of 
$\GU(1,2)$, that $\calM^{ss}$ is locally
isomorphic to $\N^{red}$  if $C^p$ is small enough.\\

We now state our results.
We assume that $p\neq 2$.
Let $k$ be an algebraically closed field extension of $\bF_p$  and
let $(X,\rho)$ be an element of $\N(k)$. As the height of the
quasi-isogeny $\rho$ is divisible by $n$ (Lem.~\ref{d14}), we may
define the morphism 
$\kappa:\N\rightarrow\Z$ by sending an element of $\N$ to the
height of the quasi-isogeny divided by $n$. 
The fibres $\N_i$ of $\kappa$  define a disjoint decomposition of
$\N$ into open and closed formal subschemes. In fact, $\N_i$ 
is empty 
if $ni$ is odd and $\N_i$ is isomorphic to
$\N_0$ if $ni$ is even (Lem.~\ref{d17}, Prop.~\ref{d21}). For the rest of this
introduction, we fix an integer $i$ with $ni$ even. 

From now on we assume that $r$ is equal to 1. Let $C$ be a $\Q_{p^2}$-vector
space of dimension $n$.
We choose a perfect skew-hermitian form $\{\cdot,\cdot\}$ on $C$ such that
there exists a self-dual $\Z_{p^2}$-lattice in 
$C$ if $n$ is odd and such that there exists no self-dual
$\Z_{p^2}$-lattice if $n$ is even.
Denote
by $H$ the special unitary group of $(C,\{\cdot,\cdot\})$ over $\Q_p$ and
denote by $\B(H,\Q_p)_{simp}$ the simplicial complex of the
Bruhat-Tits building of $H$.  
We associate to each
vertex $\Lambda$ of $\B(H,\Q_p)_{simp}$ a subset $\V(\Lambda)(k)$ of
$\N_i(k)$.
In Section \ref{i} we attach to each vertex $\Lambda$ an
odd integer $l$ with $1\leq l\leq n$, the type of
$\Lambda$. The type classifies the 
different orbits of the action of $H(\Q_p)$ on the set of vertices of
$\B(H,\Q_p)_{simp}$. 
We call a point of $\N_i(k)$ superspecial if
the underlying $p$-divisible
group is superspecial, i.e., if the
corresponding Dieudonn\'e module $M$ satisfies $FM=VM$.
Vertices of type 1 correspond to superspecial points of
$\N_i(k)$.We prove the following theorem (Prop.~\ref{s3},
Prop.~\ref{i4}, Thm.~\ref{i5})

\begin{theorem}\label{1}
The sets $\V(\Lambda)(k)$
cover $\N_i(k)$. 

Let $\Lambda$ and $\Lambda'$ be two different vertices of
$\B(H,\Q_p)_{simp}$ of 
type $l$ and $l'$  respectively.
Then the intersection of $\V(\Lambda)(k)$ and $\V(\Lambda')(k)$ is nonempty
if and only if one vertex is a neighbour of the other or if the
corresponding vertices have a common neighbour  
of type 
$l''\leq\min\{l,l'\}$.  
\end{theorem}

We associate to each vertex $\Lambda\in\B(H,\Q_p)_{simp}$ a variety
$Y_\Lambda$ over 
$\bF_p$ such
that for each algebraically closed field extension $k$ of $\bF_p$, we
have a bijection of 
$Y_\Lambda(k)$ with $\V(\Lambda)(k)$. Let $l$ be
the type of $\Lambda$ and let $U$ be the unitary group of an
$l$-dimensional hermitian space over $\F_{p^2}$. 
We obtain the
following theorem (Prop.~\ref{s1}, Thm.~\ref{s8}).

\begin{theorem}\label{2}
The variety $Y_\Lambda$ is projective, smooth and irreducible and its
dimension is equal to
$d=(l-1)/2$. 
\begin{enumerate}
\item[a)]
There exists a decomposition of $Y_\Lambda$ into a disjoint union of 
locally closed subvarieties
\begin{align*}
Y_\Lambda=\biguplus_{j=0}^{d}X_{P_j}(w_j),
\end{align*}
where each $X_{P_j}(w_j)$ is isomorphic to a Deligne-Lusztig variety 
with respect to the group $U$  and a parabolic subgroup
$P_j$ of $U$.
\item[b)]
For every $c$ with $0\leq c\leq d$, the locally closed subvariety 
$X_{P_c}(w_c)$ is equi-dimensional of dimension $c$ and its closure in
$Y_\Lambda$ is equal to
$\biguplus_{j=0}^{c}X_{P_j}(w_j)$.
The variety $X_{P_d}(w_d)$ is open, dense and irreducible of
dimension $d$ 
in $Y_\Lambda$. 
\item[c)] For every $c$ with $0\leq c< d$, the subset
  $\biguplus_{j=0}^{c}X_{P_j}(w_j)(k)$ of $Y_\Lambda(k)$ corresponds
  to the subset $\bigcup_{\Lambda'}\V(\Lambda')(k)$ in $\N_i(k)$ where
  the union is 
  taken over all neighbours $\Lambda'$ of $\Lambda$ of type $(2c+1)$.
\end{enumerate}
\end{theorem}

In the case of $GU(1,0)$ and $GU(1,1)$, the scheme $\N_0^{red}$ is a
disjoint union of infinitely many superspecial points.

In the case of $\GU(1,2)$, we define for each vertex
$\Lambda$ of type 3 a closed embedding 
of $Y_\Lambda$ into $\N_0$ such that for every algebraically closed
field extension 
$k$ of $\bF_p$ the image of $Y_\Lambda(k)$ in $\N_0(k)$ is equal to
$\V(\Lambda)(k)$.
We denote by
$\V(\Lambda)$ the image of $Y_\Lambda$. 
Let  $\C$ be the smooth and irreducible Fermat curve in $\P_{\bF_p}^2$
given by the 
equation
$$x_0^{p+1}+x_1^{p+1}+x_2^{p+1}=0.$$ 
We show that the varieties $\V(\Lambda)$ are the irreducible
components of $\N_i$  and
prove the following
explicit description of $\N^{red}$ (Thm.~\ref{t12}).

\begin{theorem}\label{3}
Let $(r,s)=(1,2)$.
The  schemes $\N_i^{red}$, with $i\in\Z$ even, are the connected components
of $\N^{red}$ which are all  isomorphic to each other.  
Each irreducible component of $\N^{red}$ is 
isomorphic to $\C$.
Two irreducible components
intersect transversally in at most one superspecial point.
Each irreducible component
contains $p^3+1$ superspecial points and
each superspecial point is the intersection
of $p+1$ irreducible components.

The scheme 
$\N^{red}$ is equi-dimensional of dimension 1 and of complete intersection.
\end{theorem}

Using the uniformization theorem, quoted above, we obtain the
following conclusions
for $\calM^{ss}$ if $C^p$ is small enough (Cor.~\ref{a4}). 

\begin{theorem}\label{4}
Let $(r,s)=(1,2)$.
The supersingular locus $\calM^{ss}$ is equi-dimensional of dimension
1 and of complete intersection. Its
singular points are the superspecial points of $\calM^{ss}$. 
Each superspecial point is the pairwise transversal intersection of $p+1$
irreducible components. Each irreducible component is isomorphic to
$\C$ and contains $p^3+1$ superspecial points. 
Two irreducible components intersect in at most one superspecial point.
\end{theorem}

Let $J$ be the group of similitudes of the isocrystal of
$(\X,\iota,\lambda)$, or equivalently, the group of similitudes of
$(C,\{\cdot,\cdot\})$ (\ref{d24}).  
Denote by $J^0$ the subgroup of all elements $g\in J$ such that the
$p$-adic valuation of the multiplier of $g$ is equal to zero.
Let $C_{J,p}$ and $C_{J,p}'$ be maximal compact  subgroups of $J$ such
that $C_{J,p}$
is hyperspecial and $C_{J,p}'$ is not hyperspecial. We obtain the
following corollary (Prop.~\ref{a3}). 
\begin{corollary}\label{5}
We have
\begin{align*}
\#\{\text{irreducible components of }\calM^{ss}\}&= 
\#(I(\Q)\backslash (J/C_{J,p}\times G(\A_f^p)/C^p)),\\
\#\{\text{superspecial points of }\calM^{ss}\}&= 
\#(I(\Q)\backslash (J/C_{J,p}'\times G(\A_f^p)/C^p)),\\ 
\#\{\text{connected components of }\calM^{ss}\}&= 
\#(I(\Q)\backslash (J^0\backslash J\times G(\A_f^p)/C^p))\\
&=\#(I(\Q)\backslash (\Z\times G(\A_f^p)/C^p)).
\end{align*}
\end{corollary}

This paper is organized as follows. In Section \ref{d} we describe the
set $\N(k)$ for $\GU(r,s)$ using classical Dieudonn\'e theory. From Section
\ref{s} on we assume $r=1$.
Section \ref{s}
contains the 
definition of the sets $\V(\Lambda)(k)$ for a lattice $\Lambda$ in an
index set $\L_i$ for every integer $i$. Furthermore, we prove Theorem \ref{2}.
In the next section, we identify the index set $\L_i$ with the set of
vertices of $\B(H,\Q_p)$ and analyse the incidence relation of the
sets $\V(\Lambda)(k)$
(Thm.~\ref{1}).
Sections \ref{g} and \ref{t} deal with the special case $\GU(1,2)$ and
Theorem \ref{3} is proved in Section \ref{t}. 
Here
our main tool is
the theory of displays of Zink (\cite{Zi2}) which is used to construct
a universal 
display over $\N_0^{red}$.
The last section
contains the transfer of the results on the moduli space $\N$ to the
supersingular locus $\calM^{ss}$ (Thm.~\ref{4}, Cor.~\ref{5}).\\

We now explain why we restrict ourselves to the signature $(1,s)$. In
the case $\GU(r,s)$ with $1<r\leq s$, it is not clear how to obtain a
similar decomposition of $\N(k)$ into subsets $\V(\Lambda)(k)$ as
above. In particular, 
one should not expect a linear closure relation order of strata as
stated in Theorem~\ref{2}.

In the case $r=1$, we expect that the pointwise decomposition of $\N$
given here can be made algebraic. However, it seems not to be
promising to construct a universal display over each variety
$Y_\Lambda$ by
using 
a basis of the isodisplay and the equations defining
$Y_\Lambda$. Indeed, 
for increasing 
$s$, these equations become quite complicated.\\
 
\noindent
{\it Acknowledgements.}
I want to thank everybody who helped me writing this paper.
Special thanks go to M. Rapoport for his advice
and his interest in my work. 
I want to
express my gratitude to U. G\"ortz and E. Viehmann for helpful 
discussions on this
subject. I am indebted to T. Wedhorn for enduring all my
questions. 
 Furthermore, I thank S. S. Kudla, Th. Zink and the referee for useful
remarks on an earlier version of this paper.



\section{Dieudonn\'e lattices in the supersingular isocrystal\\
for $\GU(r,s)$}
\label{d}

\begin{se}\label{e4}
In sections \ref{d} to \ref{t} we depart from the introduction and
denote by $E$ an unramified extension of $\mathbb{Q}_p$ of degree~2 with $p\neq 2$. Let
$\O_E$ 
be its ring of integers.
We fix a positive integer $n$ and nonnegative
integers $r$ and $s$ with $n=r+s$.  Let $\X$ be a supersingular
$p$-divisible group of height $2n$ over $\overline{\mathbb{F}}_p$ with
$\O_E$-action  
$$\mathbb{\iota}:\O_E\rightarrow \End\X$$
 such that $(\X,\iota)$ satisfies the determinant condition of
 signature $(r,s)$, i.e.,
\begin{align}
\charpol_{\overline{\F}_p}(a,\Lie\X)=(T-\varphi_0(a))^r(T-\varphi_1(a))^s\in
\overline{\F}_p[T]\label{ed21}  
\end{align}
for all $a\in\O_E$. Here we denote by $\varphi_0$ and $\varphi_1$  the
different $\Z_p$-morphisms of $\O_E$ to $\overline{\F}_p$. Let
$\lambda$ be a 
$p$-principal quasi-polarization of $\X$ such that the Rosati involution 
induced by 
   $\lambda$ induces the involution $^*$ on $\O_E$.

Consider the moduli space $\N$ over $\Spec \overline{\F}_p$ given by
the following data up to isomorphism for an $\overline{\F}_p$-scheme
$S$.
\begin{itemize}
\item A $p$-divisible group $X$ over $S$ of height $2n$ with
  $p$-principal polarization $\lambda_X$ and $\O_E$-action $\iota_X$ such 
that the Rosati
   involution induced by 
   $\lambda_X$ induces the involution $^*$ on $\O_E$. 
  We assume that $(X,\iota)$ satisfies the determinant condition of
  signature $(r,s)$. 
\item An $\O_E$-linear quasi-isogeny
  $\rho:X\rightarrow\X\times_{\Spec\overline{\F}_p}S$ such that
  $\rho\d\circ\lambda\circ\rho$ is a $\Q_p$-multiple of $\lambda_X$
  in $\Hom_{\O_E}(X,X\d)\otimes_\Z\Q$. 
\end{itemize}
The moduli space $\N$ is represented by a separated 
formal scheme which is
locally formally of finite type over $\overline{\F}_p$ (\cite{RZ}
Thm. 2.16). 
Our goal is to describe the irreducible components of $\N$ and their
intersection behaviour. 
\end{se}

\begin{se}\label{d1}
We will now study for any algebraically closed field extension $k$ of
$\overline{\F}_p$ the set $\N(k)$. Let $W(k)$ be the ring
of Witt vectors over $k$, let $W(k)_\Q$ be its quotient field and let
$\sigma$ be the Frobenius automorphism of $W(k)$. We write $W$ instead
of  $W(\overline{\F}_p)$. Denote by $\bM$ the
Dieudonn\'e module of $\X$ and by $N=\bM\otimes_\Z\Q$ the associated
supersingular isocrystal 
of dimension $2n$ with Frobenius $F$ and Verschiebung $V$. The
$\O_E$-action $\iota$ on $\X$ induces an $E$-action on $N$. 
The $p$-principal polarization $\lambda$ of $\X$ induces a perfect
alternating form  
$$\langle\cdot,\cdot\rangle:N\times N\rightarrow W_\Q$$ 
such that for all $a\in E$ and $x,y$ of $N$
\begin{align}
\langle Fx,y\rangle=\langle x,Vy\rangle^\sigma\label{ed1}
\end{align}
and 
\begin{align}
\langle ax,y\rangle=\langle x,a^*y\rangle.\label{ed2}
\end{align}
Denote by $N_k$ the isocrystal $N\otimes_{W_\Q}W(k)_\Q$. 
For a lattice $M\subset N_k$, we denote by  
\begin{align}
M^\perp=\{y\in N_k\mid \langle y,M\rangle\subset W(k)\}\label{ed7}
\end{align}
the dual lattice of $M$ in $N_k$ with respect to the form
$\langle\cdot,\cdot\rangle$. By  
covariant Dieudonn\'e theory, the tangent space $\Lie\X$ can be identified
with $M/VM$ and we obtain the following proposition.
\end{se}

\begin{prop}
For any algebraically closed field extension $k$ of $\overline{\F}_p$ we obtain the following identification.
\begin{align}
\N(k)=\{&M\subset N_k\text{ a }W(k)\text{-lattice}\mid M\text{ is
}F\text{-},V\text{- and }O_E\text{-invariant, } \notag\\ 
&
\charpol_k(a,M/VM)=(T-\varphi_0(a))^r(T-\varphi_1(a))^s\text{ for all
}a\in\O_E,\label{ed24}\\
&M=p^iM^\perp\text{ for some }i\in\Z\}.\notag
\end{align}
\end{prop}

\begin{se}\label{d13}
We will now analyze the set $\N(k)$ in the form of (\ref{ed24}). 
Consider the decomposition
\begin{align}
E\otimes_{\Q_p} W(k)_\Q&\cong W(k)_\Q\times W(k)_\Q\label{ed3}\\
a\otimes x&\mapsto (\varphi_0(a)x,\varphi_1(a)x)\notag
\end{align}
given by the two embeddings $\varphi_i:E\hookrightarrow W(k)_\Q$. It
induces a $\Z/2\Z$-grading 
\begin{align}
N_k=N_{k,0}\oplus N_{k,1}\label{ed4}
\end{align}
of $N_k$ into free $W(k)_\Q$-modules of rank $n$. Furthermore, each 
$N_{k,i}$ is totally isotropic with respect to $\langle\cdot,\cdot\rangle$ and
$F$ induces a $\sigma$-linear isomorphism $F\colon N_{k,i}\rightarrow N_{k,i+1}$.
We obtain 
$\O_E\otimes_{\Z_p}W(k)\cong W(k)\times W(k)$ analogous to
(\ref{ed3}). Therefore, every $\O_E$-invariant Dieudonn\'e module
$M\subset N_k$ has a decomposition  
$M=M_0\oplus M_1$ such that $F$ and $V$ are operators of degree 1 and
$M_i\subset N_{k,i}$. For an $\O_E$-lattice $M$ in $N_k$ we will
always denote by $M_0\oplus M_1$ such a 
decomposition. Furthermore, for $M_i$ we define the dual lattice of $M_i$ with 
respect to $\langle\cdot,\cdot\rangle$ as
$$M_i^\perp=\{y\in N_{k,i+1}\mid \langle y,M_i\rangle\subset W(k)\}$$
\end{se}

For $W(k)$-lattices $L$ and $L'$ in a finite dimensional $W(k)_\Q$-vector
space, we denote by $[L':L]$ the index of $L$ in $L'$. If $L\subset L'$,
the index is defined as the 
length of the $W(k)$-module $L'/L$. If $[L':L]=m$, we write
$L\overset{m}{\subset}L'$. In general, we define 
$$[L':L]=[L':L\cap
L']-[L:L\cap L'].$$
\begin{lem}\label{d3}
Let $M=M_0\oplus M_1$ be an $\O_E$-invariant lattice of $N_k$. Assume that 
 $M$ is invariant under $F$ and $V$. Then $M$ satisfies the
determinant condition of signature $(r,s)$ if and only if
\begin{align}
pM_0\overset{s}{\subset}FM_1\overset{r}{\subset}M_0\label{ed8}\\
pM_1\overset{r}{\subset}FM_0\overset{s}{\subset}M_1.\label{ed9}
\end{align}
\end{lem}

\begin{proof} Consider the decomposition
$$M/VM=M_0/VM_1\oplus M_1/VM_0.$$
The determinant condition is equivalent to the condition that
 $VM_1$ is of index $r$ in $M_0$ and $VM_0$ is of index
 $s$ in $M_1$. Since $FV=VF=p\id_M$, we obtain that $pM_1$ is of index $r$ in
 $FM_0$ and $pM_0$ is of index 
$s$ in $FM_1$.
\end{proof}

\begin{se}\label{d20}
Let $\bM_k=\bM_{k,0}\oplus\bM_{k,1}$ be the Dieudonn\'e module of
$\X_k$ as in \ref{d1}. 
For a Dieudonn\'e lattice $M\in\N(k)$, denote by $\vol(M)=[\bM_k:M]$  the
volume of $M$.
\end{se}

\begin{lem}\label{d14}
Let $M\in\N(k)$.
If  $M=p^iM^\perp$ for some integer $i$, then
$\vol(M)=ni$. 
\end{lem}

\begin{proof}
Since 
both vector spaces $N_{k,0}$ and $N_{k,1}$ are maximal totally
isotropic with respect to $\langle\cdot,\cdot\rangle$, the condition
$M=p^iM^\perp$ is equivalent to the two conditions  
\begin{align*}
p^iM_0^{\perp}=M_1,\\
p^iM_1^\perp=M_0. 
\end{align*}
By duality the last two conditions are equivalent.
We obtain
\begin{align*}
\vol(M)&=[\bM_{k,0}:M_0]+[\bM_{k,1}:M_1]\\
&=[M_0^\perp:\bM_{k,0}^\perp]+[\bM_{k,1}:M_1]\\
&=[p^{-i}M_1:\bM_{k,1}]+[\bM_{k,1}:M_1]\\
&=ni
\end{align*}
which proves the claim.
\end{proof}

\begin{se}\label{d19}
Let $M\in\N(k)$ and let $(X,\rho)$ be the corresponding 
$p$-divisible group and 
quasi-isogeny. 
Denote by $\height(\rho)$ the height of $\rho$. 
Then $\height(\rho)=\vol(M)$. 
As the hight of a quasi-isogeny of p-divisible groups over $S$ is locally 
constant over $S$, we obtain by Lemma \ref{d14} a morphism
\begin{align} 
\kappa:\N&\rightarrow\Z\label{ed23}\\
(X,\rho)&\mapsto \frac{1}{n}\height(\rho).\notag
\end{align}
For $i\in\Z$ the fibre $\N_i=\kappa^{-1}(i)$ is the open and closed formal
subscheme of $\N$ of quasi-isogenies of height $ni$.
\end{se}

\begin{lem}\label{d17}
If $ni$ is odd, the formal scheme $\N_i$ is empty. 
\end{lem} 

\begin{proof}
Let $M$ be an element of $\N(k)$. 
Since both $\bM_k$ and $M$ satisfy the determinant condition of
signature $(r,s)$, we obtain
by Lemma \ref{d3} that 
\begin{align*}
\vol(M)&=[\bM_{k,0}:M_0]+[\bM_{k,1}:M_1]\\
&=[\bM_{k,0}:M_0]+[F\bM_{k,1}:FM_1]\\
&=2[\bM_{k,0}:M_0]+[F\bM_{k,1}:\bM_{k,0}]+[M_0:FM_1]\\
&=2[\bM_{k,0}:M_0].
\end{align*}
Hence $\vol(M)$ 
is even, i.e., $ni$ is even. 
\end{proof}

\begin{se}\label{d4}
Let $N_k=N_{k,0}\oplus N_{k,1}$ be as in \ref{d13}. To describe the
set $\N(k)$, it is convenient  
to express $\N(k)$ in terms of 
$N_{k,0}$. Let $\tau$ be the $\sigma^2$-linear operator $V^{-1}F$ on
$N_k$. Then $N_{k,0}$ and $N_{k,1}$ are both $\tau$-invariant. 
Denote by $\Q_{p^2}$ the unique
unramified extension of degree 2 of $\Q_p$ in $W_\Q$  and denote by
$\Z_{p^2}$ its ring of 
integers. 
Since the isocrystal $N_k$ is supersingular, $(N_k,\tau)$ 
is isoclinic with slope zero. 
Thus $(N_{k,i},\tau)$ is isoclinic of slope zero for $i=0,1$, i.e.,
there exists a $\tau$-invariant lattice in $N_{k,i}$. For every
$\tau$-invariant lattice $M_i\subset N_{k,i}$, there 
exists a $\tau$-invariant basis of $M_i$
(Thm.~of Dieudonn\'e, \cite{Zi1} 6.26). 
Let $C$ be the $\Q_{p^2}$-vector space of all $\tau$-invariant
elements of $N_{k,0}$ and let $M_0^\tau$ be the $\Z_{p^2}$-module of
$\tau$-invariant elements of $M_0$.  
We obtain 
$$M_0=M^\tau_0\otimes_{\Z_{p^2}}W(k),$$
$M^\tau_0$ is a lattice in $C$ and
$$N_{k,0}=C\otimes_{\Q_{p^2}}W(k)_\Q.$$
Note that the $\Q_{p^2}$-vector space $C$ does not depend on $k$. 
We write $C_k$ for the base change $C\otimes_{\Q_{p^2}}W(k)_\Q$.
\end{se}

\begin{se}\label{d8}
We define a new form on $C_k$ by
$$\{x,y\}:=\langle x,Fy\rangle.$$
This is
a perfect form on $C_k$ linear in the first and $\sigma$-linear in the second 
variable. By (\ref{ed1}) we obtain the following property of $\{\cdot,\cdot\}$
\begin{align}
\{x,y\}=-\{y,\tau^{-1}(x)\}^{\sigma},\label{ed5}
\end{align}
which in turn implies
\begin{align}
\{\tau(x),\tau(y)\}=\{x,y\}^{\sigma^2}.\label{ed6}
\end{align}
For a $W(k)$-lattice $L$ in $C_k$, denote by $L\d$ the dual of $L$
with respect to the form $\{\cdot,\cdot\}$ defined by 
\begin{align}
L\d=\{y\in C_k\mid \{y,L\}\subset W(k)\}.\label{ed10}
\end{align}
We obtain by (\ref{ed5}) that
\begin{align}
(L\d)\d=\tau(L).\label{ed11}
\end{align}
In particular, taking the dual is not an involution on the set of
lattices in $C_k$. 
The identity (\ref{ed6}) implies that 
\begin{align}
\tau(L\d)=\tau(L)\d.\label{ed12}
\end{align}

The form $\{\cdot,\cdot\}$ on $C_k$ induces by restriction to $C$ a perfect
skew-hermitian form  on  $C$ with respect to $\Q_{p^2}/\Q_p$
which we will again denote by $\{\cdot,\cdot\}$.
For the perfect form to be skew-hermitian means that
it is linear in the first 
and $\sigma$-linear in the second variable and we have
$$\{x,y\}=-\{y,x\}^\sigma,$$
where the Frobenius $\sigma$ is an involution on $\Q_{p^2}$.

It is clear that for each $\tau$-invariant lattice $A$ of $C_k$ we obtain
\begin{align}
(A\d)^\tau=(A^\tau)\d,\label{ed15}
\end{align}
where the second dual is taken with respect to the skew-hermitian form
$\{\cdot,\cdot\}$ on $C$.  
\end{se}

\begin{prop}\label{d5}
There is a bijection between $\N(k)$ and
\begin{align}
\D(C)(k)=\{\text{lattices }A\subset C_k\mid
p^{i+1}A\d\overset{r}{\subset}A\overset{s}{\subset}p^iA\d,\text{ for
  some } i\in\Z\}.\label{ed14} 
\end{align}
The bijection is obtained
by associating to $M=M_0\oplus M_1\in\N(k)$ the lattice
$M_0$ in $C_k$.  
\end{prop}

\begin{rem}\label{d6}
Note that by duality and (\ref{ed11}) the chain condition  
\begin{align}
p^{i+1}A\d\overset{r}{\subset}A\overset{s}{\subset}p^iA\d
\end{align}
is equivalent to the chain condition
$$p^{-i}\tau(A)\overset{s}{\subset}A\d\overset{r}{\subset}p^{-i-1}\tau(A),$$
which is equivalent to
\begin{align}
p^{i+1}A\d\overset{r}{\subset}\tau(A)\overset{s}{\subset}p^iA\d.
\end{align}
\end{rem}

\begin{de}
An element $M\in\N(k)$ is called superspecial if $F(M)=V(M)$, i.e., if
$M$ is $\tau$-invariant. A  lattice $A\subset C$ is  superspecial if
and only if $A$ is $\tau$-invariant.
\end{de}
Note that $M$ is superspecial if and 
only if the corresponding lattice $A=M_0$ is superspecial.

\begin{se}\label{d7}
{\it Proof of Proposition \ref{d5}.}\\
Let $M=M_0\oplus M_1$ be an $\O_E$-invariant lattice which is stable
under $F$ and $V$. 

{\it Claim}: $M=p^iM^\perp$ with respect to
$\langle\cdot,\cdot\rangle$ if and only if $FM_1=p^{i+1}M_0\d$.\\  
Indeed,
the lattice $M$ is equal to $p^iM^\perp$ if and only if
$p^iM_0^\perp=M_1$ (proof of Lem.~\ref{d14}).
We have
\begin{align}
F(M_0^{\perp})&=\{y\in N_{k,0}\mid \langle
F^{-1}y,M_0\rangle \subset W(k)\}\notag\\
&= \{y\in N_{k,0}\mid \{p^{-1}y,M_0\}\subset W(k)\}\notag\\
&=pM_0\d\notag
\end{align}
which proves the claim.

Let $M$ be an element of the set $\N(k)$. Since $FM_1$ is equal to
$p^{i+1}M_0\d$ for some $i\in\Z$, we obtain from (\ref{ed8}) for
$M_0=A$ the chain condition 
\begin{align}
p^{i+1}A\d&\overset{r}{\subset}A\overset{s}{\subset}p^iA\d.
\end{align}
Hence $A$ is an element of $\D(C)(k)$.
Conversely, associate to a lattice $A$ of $\D(C)(k)$ the lattice
$A\oplus F^{-1}(p^{i+1}A\d)\subset N_{k,0}\oplus N_{k,1}$. It is an element of
$\N(k)$ by the same arguments. 
\hfill{$\Box$}
\end{se}

\begin{lem}\label{d11}
Let $t\in\Z_{p^2}^\times$ with $t^\sigma=-t$ and let $V$ be a
$\Q_{p^2}$-vector space of dimension $n$. 
Let $I_n$ be the identity matrix of rank $n$ and let $J_n$ be the
matrix
$$J_n=\begin{pmatrix}
p& & &\\
 &1& &\\
 & &\ddots& \\
 & & &1\\
\end{pmatrix}.$$
 There exist 
two perfect skew-hermitian forms on $V$ up to isomorphism. These forms
correspond 
to $tI_n$ and to 
$tJ_n$
respectively. Furthermore, if $M$ is a lattice in $V$ and $i\in\Z$ with
\begin{align}
p^{i+1}M\d\overset{r}{\subset}M\overset{s}{\subset}p^iM\d,\label{ed20}
\end{align}
then $s\equiv ni\mod 2$ in the first case and $s\not\equiv ni\mod 2$
in the second case. 
In particular,
the form $t\begin{pmatrix}
I_r&\\
 &pI_s\\
\end{pmatrix}
$ is isomorphic to the form  $tI_n$ if $s$ is even and is isomorphic
to $tJ_n$ if $s$ is 
odd. 
\end{lem}

\begin{proof}
Let $\{\cdot,\cdot\}$ be a perfect skew-hermitian form on $V$. Let 
$U$ be the unitary group over $\Q_p$ associated  to the pair
$(V,\{\cdot,\cdot\})$. 
As $H^1(\Q_p,U)\cong\Z/2\Z$ (\cite{Ko1} $\S 6$), there exist two non isomorphic skew-hermitian
forms on $V$. 
Choose a basis of $V$. 

{\it Claim:} The skew-hermitian forms on $V$ corresponding to the matrix
$tI_n$ and to the matrix
$tJ_n$ are not isomorphic. 

Indeed,
assume that these forms are
isomorphic. Then there exists a lattice $M$ in $V$ with $M=M\d$ and a
lattice $L$ in $V$ with  
$$pL\d\overset{n-1}{\subset}L\overset{1}{\subset}L\d.$$
Let $k$ be an integer such that $L\d$ is contained in $p^k M$. We obtain
$$p^{-k}M\overset{l}{\subset}L\overset{1}{\subset}L\d\overset{l}{\subset}p^k
M.$$ 
Thus $2kn=2l+1$ which is a contradiction. Therefore, the two forms are
not isomorphic. 

Let $M$ be a lattice in $V$ which satisfies (\ref{ed20}). 
A similar index argument as above shows that $s\equiv ni\mod 2$ if the
skew-hermitian form is isomorphic to $tI_n$ and $s\not\equiv ni\mod 2$
in the other case.  
\end{proof}

\begin{lem}\label{d18}
Let $V$ be a $\Q_{p^2}$-vector space of dimension $n$ and let
 $\{\cdot,\cdot\}$ be a perfect   
 skew-hermitian form on $V$.  Let $M\subset V$ be a $\Z_{p^2}$-lattice with
\begin{align}
pM\d\overset{r}{\subset}M\overset{s}{\subset}M\d.
\end{align}
Then there exists a basis of $M$ such that the form $\{\cdot,\cdot\}$ with respect to
this basis is given by the matrix
$$t\begin{pmatrix}
I_r& \\
 &pI_s\\
\end{pmatrix}.
$$
\end{lem}

\begin{proof}
The lemma is proved by an analogue of the Gram-Schmidt
orthogonalization. 
\end{proof}

\begin{prop}\label{d12}
There exists a basis of $C$ such that the form $\{\cdot,\cdot\}$ with
respect to this basis is given by the 
matrix $(tI_n)$ if $s$ is even and by $(tJ_n)$ if 
$s$ is odd.

In particular, the isocrystal $N$ with $\O_E$-action and perfect form
$\langle\cdot,\cdot\rangle$ as in \ref{d1} is uniquely determined up to
isomorphism.  
\end{prop}

\begin{proof}
The image of the Dieudonn\'e lattice $\bM_k\in\N(k)$  under the
bijection of Proposition \ref{d5} satisfies 
$$p\bM_{k,0}\d\overset{r}{\subset}\bM_{k,0}\overset{s}{\subset}\bM_{k,0}\d.$$
Thus the proposition follows from Lemma~\ref{d18}.
\end{proof}

\begin{se}\label{d24}
Let $J$ be the group of isomorphisms of the isocrystal $N$ with
additional structure, i.e., 
\begin{align*}
\begin{aligned}
J=\{g\in\GL_{E\otimes_{\Q_p}W(\bF_p)_\Q}(N)
\mid   Fg=gF;\ \langle gx,gy\rangle=c(g)\langle x,y\rangle& \\
\text{ with
  }c(g)\in(W(\bF_p)_\Q)^\times\}&. 
\end{aligned}
\end{align*}
Then $J$ is the group $\End^\circ_{\O_E,\lambda}(\X)^\times$ of $\O_E$-linear
quasi-isogenies $\rho$ of the $p$-divisible group $\X$ of \ref{e4} such that
$\lambda\circ\rho$ is a 
$\Q_p^\times$-multiple of 
$\rho\d\circ\lambda$. 

An element $g\in J$ acts on $\N$ by sending an element
$(X,\iota_X,\lambda_X,\rho)$ to $(X,\iota_X,\lambda_X,g\circ\rho)$.
Consider the decomposition $N=N_0\oplus N_1$ of the isocrystal $N$ as
in \ref{d13}. We will show that $J$ can be identified with the group of
unitary similitudes of 
the hermitian space $(C,\{\cdot,\cdot\})$.

Indeed,
let $g\in J$. As $g$ is $\O_E$-linear, it respects
the grading of $N$. Since $g$
commutes with $F$, the action of $g$ on $N$ is uniquely determined by
its action on $N_0$ and we obtain $\{gx,gy\}=c(g)\{x,y\}$ for all
elements $x,y\in N_0$. As $g$ commutes with $\tau$, 
it is defined over $\Q_{p^2}$, i.e., an automorphism of
$N_0^\tau=C$. In particular, $c(g)\in\Q_{p^2}$.  

Let $v_p$ be the $p$-adic valuation on $\Q_{p^2}$. It defines a
morphism $\theta:J\rightarrow\Z$ by sending an element $g\in J$ to
$v_p(c(g))$. 
Denote by $J^0$ the kernel of $\theta$.  

For  $\rho\in\End^\circ_{\O_E,\lambda}(\X)^\times$,
let $g\in J$ be the corresponding automorphism of the isocrystal and
let $\alpha=v_p(c(g))$. We obtain $g\bM=c(g)(g\bM)^\perp$, hence by
Lemma~\ref{d14}, the height of $\rho$ is equal to $n\alpha$. By
\ref{d19} the element $g$ defines an isomorphism of $\N_i$ with
$\N_{i+\alpha}$ for every integer $i$.
\end{se}

\begin{lem}\label{d23}
The image of $\theta$ is equal to $\Z$ if $n$ is even and equal to
$2\Z$ if $n$ is odd. 

In particular, there exists a quasi-isogeny
$\rho\in\End_{\O_E,\lambda}(\X)$ of height $h\in\Z$ if and only if $h$ is
divisible by
$n$ if $n$ is even and if $h$  is divisible by $2n$ if $n$ is
odd. 
\end{lem}

\begin{proof}
For $g=p\id_N$ we have $c(g)=p^2$. 
If $n$ is odd, the same argument as in Lemma~\ref{d17} shows that the
image of $\theta$ is contained in $2\Z$. 
Thus we may assume that $n$ is even. It is sufficient to show that
there exists an element $g\in\GL(C)$ such
that $\{gx,gy\}=p\{x,y\}$ for all $x,y\in C$. 

Let $T_1$ and $T_2$  be the matrices in $GL(C)$ given by
\begin{align*}
T_1&=
\begin{pmatrix}
 &&1\\
&\a&\\
1&&\\
\end{pmatrix}
&
T_2&=
\begin{pmatrix}
&&&p\\
 &&1&\\
&\a&&\\
1&&&\\
\end{pmatrix}.
\end{align*}
By Lemma~\ref{d11} the perfect skew-hermitian form on $C$ induced by $tT_1$
is isomorphic to 
$tI_n$ and the perfect skew-hermitian form induced by $tT_2$ is
isomorphic to $tJ_n$. 
We set $g=\diag(p^{n/2},1^{n/2})$.
Then $g$ satisfies the claim. 
\end{proof}

\begin{se}\label{d22}
For $i\in\Z$ we denote by $\D_i(C)(k)$ the image of $\N_i(k)$ (\ref{d19}) under
the bijection of \ref{d5}, i.e.,
\begin{align}
\D_i(C)(k)=\{A\in\D(C)(k)\mid
p^{i+1}A\d\overset{r}{\subset}A\overset{s}{\subset}p^iA\d\}.\label{es14} 
\end{align}
We have a  decomposition of $\D(C)(k)$ into a disjoint union of the
sets $\D_i(C)(k)$, 
\begin{align}
\D(C)(k)=\biguplus_{i\in\Z}\D_i(C)(k).
\end{align}
The sets $\D_i(C)(k)$ are invariant under the action of $\tau$ on $\D(C)(k)$.
By Lemma~\ref{d17} we know that $\D_i(C)(k)$ is empty if $ni$ is
odd. 
\end{se}

\begin{prop}\label{d21}
Let $i$ be an integer such that $ni$ is even. If $n$ is even, let
$g$ be an element of $J$ such that $v_p(c(g))=-1$ (Lem.~\ref{d23}).
There exists an
isomorphism
\begin{align}
\Psi_i:\N_i\stackrel{\sim}{\longrightarrow}\N_0\label{ed25}
\end{align}
such that the following holds.
If $i$ is even,
$\Psi_i$ induces on $k$-rational points the bijection
\begin{align*}
\Psi_i:\D_i(C)(k)&\stackrel{\sim}{\longrightarrow}\D_0(C)(k)\\
A&\mapsto p^{-\frac{i}{2}}A\notag
\end{align*}
and if $i$ is odd, $\Psi_i$ induces the bijection
\begin{align*}
\Psi_i:\D_i(C)(k)&\stackrel{\sim}{\longrightarrow}\D_0(C)(k)\\
A&\mapsto p^{\frac{-i+1}{2}}g(A).\notag
\end{align*}
\end{prop}

\begin{proof}
If $i$ is even, the multiplication $p^{-i/2}:\X\rightarrow\X$ defines an
isomorphism of $\N_i$ with $\N_0$ which satisfies the claim. If $i$ is
odd,  the quasi-isogeny $p^{(-i+1)/2}g$
induces an isomorphism of $\N_i$ with $\N_0$ which satisfies the claim.
\end{proof}


\section{The set structure of $\N$ for $\GU(1,s)$}
\label{s}

\begin{se}\label{s16}
From now on, we will restrict ourselves to the case of $\GU(1,s)$. 
Our goal is to describe the irreducible components of $\N_i$
 for any integer $i$ with $ni$ even. In this section
we will define irreducible varieties over $\overline{\F}_p$ of
dimension equal to the dimension of $\N_i$ such that the $k$-rational points
of these varieties cover $\N_i(k)$ for every algebraically closed field extension $k$
of
$\overline{\F}_p$. 

We will always denote by $k$ 
an algebraically closed field extension of $\overline{\F}_p$. Let $C_k$ be the
$W(k)_\Q$-vector space of dimension $n=s+1$ with $\sigma^2$-linear
operator 
$\tau$ as in \ref{d4}.  
By \ref{d8} the vector space $C_k$ is  equipped with a
perfect form $\{\cdot,\cdot\}$ linear in the first and
$\sigma$-linear in the second variable which satisfies 
$$\{x,y\}=-\{y,\tau^{-1}(x)\}^\sigma$$
for all $x,y\in C_k$.
By Proposition \ref{d12} the restriction of  $\{\cdot,\cdot\}$ to $C$  is a
perfect skew-hermitian form equivalent to the skew-hermitian form
induced by $tI_n$ if $n$ is odd and equivalent to $tJ_n$ if $n$ is
even. 

Let $i\in\Z$.
The set $\D_i(C)(k)$ as in (\ref{es14}) consists of all lattices $A$ in
$C_k$ such that
\begin{align}
p^{i+1}A\d&\overset{1}{\subset}A\overset{n-1}{\subset}p^iA\d.\label{es24}
\end{align}
Equivalently,
\begin{align}
p^{i+1}A\d&\overset{1}{\subset}\tau(A)\overset{n-1}{\subset}p^iA\d.
\label{es25}  
\end{align}
By Lemma~\ref{d17} we know that $\D_i(C)(k)$ is empty if $ni$ is odd.
If $ni$ is even, the set $\D_i(C)(k)$ is nonempty by Proposition~\ref{d21}.
Therefore, we will always assume that $ni$ is even. 
\end{se}

We need the following crucial lemma.
\begin{lem}\label{s17}
Let $A$ be a lattice in $\D_i(C)(k)$. There exists an integer $d$ with
$0\leq d\leq s/2$ 
and such that 
$$\Lambda=A+\tau(A)+...+\tau^d(A)$$
is a $\tau$-invariant lattice. If $d$ is  minimal
with this condition, we have
\begin{align}
p^{i+1}\Lambda\d\subset
p^{i+1}A\d\overset{1}{\subset}A\overset{d}{\subset}\Lambda
\overset{n-2d-1}{\subset}p^i\Lambda\d\subset 
p^iA\d\label{es4} 
\end{align}
and $p^{i+1}\Lambda\d$ is of index $(2d+1)$ in $\Lambda$.
\end{lem}

The proof of Lemma \ref{s17} will use explicitly that the index of
$p^{i+1}A\d$ in $A$ is equal to~1. It does not work in the case of
general index.

\begin{proof}
For every nonnegative integer $j$, denote by $T_j$ the lattice
\begin{align}
T_j=A+\tau(A)+...+\tau^j(A).\label{es26}
\end{align}
We have 
\begin{align}
T_{j+1}=T_j+\tau(T_j).
\end{align}
Let $d\geq 0$ be minimal with $T_d=\tau(T_d)$, i.e., $\tau(T_j)\neq T_j$ for 
every $0\leq j<d$. Such an integer exists (\cite{RZ} Prop. 2.17). If $d=0$, 
there is nothing to prove, hence we may assume that $d\geq 1$.

\emph{Claim:} If $d\geq 2$, then for $2\leq j\leq d$ 
\begin{alignat}{2}
\tau(T_{j-2})&\overset{1}{\subset}T_{j-1}&\overset{1}{\subset}&T_j
\label{es27}\\ 
\tau(T_{j-2})&\overset{1}{\subset}\tau(T_{j-1})&\overset{1}{\subset}&T_j.
\label{es28} 
\end{alignat}
In particular, $A=T_0$ is of index $j$ in $T_j$.

Indeed, if $j=2$ we obtain $p^{i+1}A\d\overset{1}{\subset}A$ and
$p^{i+1}A\d\overset{1}{\subset}\tau(A)$ by (\ref{es24}) and
(\ref{es25}).  
Hence either $A=\tau(A)$ or $A$ and $\tau(A)$ are both of index 1 in
$T_1=A+\tau(A)$. Since 
$d\geq 1$, the second possibility occurs. As $A$ is of index 1 in
$T_1$, the lattice $\tau(A)$ is of index 1 in $\tau(T_1)$. We obtain
that either $T_1=\tau(T_1)$ or that $T_1$ and $\tau(T_1)$ are both of
index 1 in $T_2=T_1+\tau(T_1)$. Since $d\geq 2$, the second
possibility occurs which proves
the claim for $j=2$.
The general case follows by induction on $j$.

We will now show that $T_d$ is contained in $T_d\d$.
By (\ref{es24}) and (\ref{es25}), we have
\begin{alignat*}{2}
A+\tau(A)&\subset p^iA\d&\subset& p^{-1}\tau(A),\\
\tau(A)+\tau^2(A)&\subset p^i\tau(A)\d&\subset& p^{-1}\tau(A),
\end{alignat*}
hence
\begin{align}
A+\tau(A)+\tau^2(A)\subset p^{-1}\tau(A).\label{es29}
\end{align}
 Using (\ref{es29}) we obtain
\begin{align*}
T_d&=A+...+\tau^d(A)\\
&\subset p^{-1}\tau(A)+...+p^{-1}\tau^{d-1}(A)\subset
p^{-1}T_{d-1} 
\end{align*}
Since $T_d$ is $\tau$-invariant, $T_d$ is contained in
$p^{-1}\tau^l(T_{d-1})$ for every integer $l$, thus 
\begin{align}
T_d&\subset p^{-1}\bigcap_{l\in\Z}\tau^l(T_{d-1}).\notag
\end{align}
Now
\begin{align}
\bigcap_{l\in\Z}\tau^l(T_{d-1})
=\bigcap_{l\in\Z}\tau^l(A).\label{es30}
\end{align}
Indeed, this is clear if $d=1$ since $T_0=A$. If $d\geq 2$, we obtain by
(\ref{es27}) and (\ref{es28}) that
$$T_{d-1}\cap\tau(T_{d-1})=\tau(T_{d-2}),$$
hence
$$\bigcap_{l\in\Z}\tau^l(T_{d-1})=\bigcap_{l\in\Z}\tau^l(T_{d-2})$$
and (\ref{es30}) follows by induction. Since $d\geq 1$, we have
$A\neq\tau(A)$, and hence $A\cap \tau(A)=p^{i+1}A\d$ by
(\ref{es24}) and (\ref{es25}). 
We obtain  
\begin{align}
T_d&\subset p^i \bigcap_{l\in\Z}\tau^l(A)\d.\label{es5}
\end{align}
Dualizing (\ref{es26}) for $j=d$ shows that
$$T_d\d=A\d\cap...\cap \tau^d(A\d),$$
hence by (\ref{es5})
\begin{align*}
T_d&\subset p^i\bigcap_{l\in\Z}\tau^l(A)\d\\
&=p^i\bigcap_{l\in\Z}\tau^l(T_d)\d=p^iT_d\d.
\end{align*}
The last equality is satisfied because
$T_d$, and hence $T_d\d$, are $\tau$-invariant. We obtain
\begin{align}
p^{i+1}T_d\d\subset p^{i+1}A\d\overset{1}{\subset}A\subset T_d\subset
p^iT_d\d\subset p^iA\d.\label{es31} 
\end{align}
 Using (\ref{es31}) and $A\overset{s}{\subset}p^iA\d$ we see
that $d\leq s/2$ which proves the claim.
\end{proof}

\begin{de}\label{s22}
For $i\in\Z$ with $ni$ even, let $\L_i$ be the set of all 
lattices $\Lambda$  in $C$ satisfying 
\begin{align}
p^{i+1}\Lambda\d\subsetneqq \Lambda\subset p^i\Lambda\d.\label{es32}
\end{align} 
Let $\L=\biguplus_{i\in\Z}\L_i$ be the disjoint union of the sets $\L_i$.
We say that $\Lambda\in\L_i$ is of type $l$, if $p^{i+1}\Lambda\d$ is
of index $l$ in $\Lambda$. 
Denote by $\L_i^{(l)}$ the set  
of all lattices of type $l$ in $\L_i$.
For $\Lambda\in\L_i$ let $\Lambda_k=\Lambda\otimes_{\Z_{p^2}}W(k)$. We
define 
\begin{align}
\V(\Lambda)(k)=\{A\subset\Lambda_k\mid
p^{i+1}A\d\overset{1}{\subset}A\}. 
\end{align}
\end{de}

\begin{rem}\label{s18}
\begin{enumerate}
\item[a)] Let $\Lambda\in\L_i$. The type $l$ of $\Lambda$ is always
  an odd integer with 
  $1\leq l\leq n$. 
Indeed, it is clear that $1\leq l\leq n$.
Since $ni$ is even, the integer $n-l$ is even if and only if $n$ is
odd (Lem.~\ref{d11}). 
\item[b)] By \ref{d4} and \ref{d8}, there is a bijection between
  $\L_i$ and the set of all 
$\tau$-invariant lattices in $C_k$ which satisfy (\ref{es32}).
Via this bijection the superspecial lattices in $\D_i(C)(k)$
  correspond to the lattices $\Lambda\in\L_i$ of type 1.
\item[c)] Let $\Lambda\in\L_i$.
By duality a lattice $A$ in $\V(\Lambda)(k)$ satisfies
\begin{align}
p^{i+1}\Lambda_k\subset p^{i+1}\Lambda_k\d\subset p^{i+1}A\d\subset
A\subset\Lambda_k.\label{es6} 
\end{align}
\item[d)] The  isomorphism $\Psi_i$ of Proposition~\ref{d21} induces
 for every $i$ an 
 isomorphism $\Psi_i$ of $\L_i^{(l)}$ with $\L_0^{(l)}$ such that
 $\V(\Psi_i(\Lambda))(k)=\Psi_i(\V(\Lambda)(k))$. 
\end{enumerate}
\end{rem}

\begin{se}
The sets $\V(\Lambda)(k)$ will be identified with the $k$-rational points
of an irreducible smooth variety. 
In Section \ref{t} we will prove that for $\GU(1,2)$ the varieties
corresponding to  lattices  $\Lambda$ of maximal type are isomorphic to
the irreducible components of $\N$. We will  
start with some basic properties of $\V(\Lambda)(k)$.
\end{se}

\begin{prop}\label{s3}
\begin{enumerate}
\item[a)] We have
  $\D_i(C)(k)=\bigcup_{\Lambda\in\L_i}\V(\Lambda)(k)$. In particular,
  $\L_i\neq\emptyset$. 
\item[b)] For $\Lambda\in\L_i$ 
  and $\Lambda'\in\L_j$ with $i\neq j$, we have
  $\V(\Lambda)(k)\cap\V(\Lambda')(k)=\emptyset$.
\item[c)] Let $\Lambda$ and $\Lambda'$ be elements of $\L_i$.
\begin{enumerate}
\item[(i)] If $\Lambda\subset\Lambda'$, then
  $\V(\Lambda)(k)\subset\V(\Lambda')(k)$. 
\item[(ii)] We have 
$$\V(\Lambda)(k)\cap\V(\Lambda')(k)=\begin{cases}
\V(\Lambda\cap\Lambda')(k)&\text{ if }\Lambda\cap\Lambda'\in\L_i\\
\emptyset&\text{ otherwise.}
\end{cases}$$
\end{enumerate}
\end{enumerate}
\end{prop}

\begin{proof}
Statements a), b) and c)(i) follow from the definition and Lemma
\ref{s17}. 

To prove c)(ii), let $A$ be an element of
$\V(\Lambda)(k)\cap\V(\Lambda')(k)$. By  
(\ref{es6}) we have
\begin{align}
p^{i+1}\Lambda\d\subset p^{i+1}A\d\subsetneqq A\subset\Lambda_k\subset p^i\Lambda_k\d\subset p^iA\d
\end{align}
and similarly for $\Lambda_k'$ instead of $\Lambda_k$. We obtain
\begin{align}
p^{i+1}(\Lambda_k\cap\Lambda_k')\d\subset p^{i+1}A\d\subsetneqq
A\subset(\Lambda_k\cap\Lambda_k')\subset\Lambda_k\subset
p^i\Lambda_k\d\subset p^i(\Lambda_k\cap\Lambda_k')\d\subset p^iA\d, 
\end{align}
hence $\Lambda\cap\Lambda'\in\L_i$. A similar calculation shows the equality
\[
\V(\Lambda)(k)\cap\V(\Lambda')(k)=\V(\Lambda\cap\Lambda')(k).\qedhere
\]
\end{proof}

\begin{prop}\label{s12}
Let $\Lambda,\Lambda'\in\L_i$.
\begin{enumerate}
\item[a)] $\V(\Lambda)(k)$ contains a superspecial
  lattice. Furthermore, $\#\V(\Lambda)(k)=1$ if and only if $\Lambda$
  is of type $1$. In this case $\V(\Lambda)(k)=\{\Lambda_k\}$. 
\item[b)] $\V(\Lambda')(k)\subset\V(\Lambda)(k)$ if and only if
  $\Lambda'\subset\Lambda$. In particular,
  $\V(\Lambda')(k)=\V(\Lambda)(k)$ if and only if $\Lambda'=\Lambda$. 
\item[c)] Let $l$ be the type of $\Lambda$. For every odd integer
  $1\leq l'\leq n$, there exists a lattice $\Lambda'\in\L_i$ of type
  $l'$ with $\Lambda'\subset\Lambda$ if $l'\leq l$ and
  $\Lambda\subset\Lambda'$ if $l\leq l'$.  
\end{enumerate}
In particular, the maximal sets $\V(\Lambda)(k)$ are the sets with
$\Lambda$ of type $n$ if n is odd and of type $n-1$ if n is even. 
\end{prop}

\begin{se}\label{s19}
We will prove the proposition in \ref{s21}. 
For this proof we  need a 
description of the sets $\V(\Lambda)(k)$ in terms of linear
algebra. We  
first consider the case $i=0$. 
Let $\Lambda\in\L_0$ be of type $l$. We associate to $\Lambda$ the
$\F_{p^2}$-vector spaces    
  $V=\Lambda/p\Lambda\d$ and $V'=\Lambda\d/\Lambda$.  
 By (\ref{es32}) the vector space $V$ is of dimension $l$ and $V'$
 is of dimension $n-l$. For $z\in\Z_{p^2}$ denote by $\overline{z}$
 its image in $\F_{p^2}$.
The skew-hermitian form $\{\cdot,\cdot\}$ on $C$ induces a perfect
skew-hermitian form $(\cdot,\cdot)$ on $V$ by
$$(\overline{x},\overline{y})=\overline{\{x,y\}}\in\F_{p^2}$$
for $\overline{x},\overline{y}\in V$ and lifts
$x,y\in\Lambda$. Similarly, we obtain a perfect skew-hermitian form on $V'$ by
$$(\overline{x},\overline{y})'=\overline{(p\{x,y\})}\in\F_{p^2}$$
for $\overline{x},\overline{y}\in V'$ and lifts
$x,y\in\Lambda\d$.
 
Let 
$\tau$ be the operator on $V_k=V\otimes_{\F_{p^2}}k$ defined by the
Frobenius of $k$ over $\F_{p^2}$. 
We denote again by $(\cdot,\cdot)$ the induced form on
$V_k$ 
given by
\begin{align*}
V_k\times V_k&\rightarrow k\\
(v\otimes x,w\otimes y)&\mapsto xy^\sigma(v,w).
\end{align*}
This form is linear in the first and $\sigma$-linear in the
second variable and
satisfies 
\begin{align}
(x,y)=-(y,\tau^{-1}(x))^\sigma.\label{es10}
\end{align}
For a subspace $U$ of $V_k$, we denote by $U^\perp$ the orthogonal complement
$$U^\perp=\{x\in V_k\mid(x,U)=0\}.$$
By (\ref{es10}) we obtain  $(U^\perp)^\perp=\tau(U)$ and
$\tau(U)^\perp=(\tau(U))^\perp$ analogous to \ref{d8}.

Let $G$ be the unitary group associated to $(V,(\cdot,\cdot))$. Since
$\hH^1(\F_p,G)=0$, 
there exists up to isomorphism only one skew-hermitian form on $V$ and
similarly for $V'$. 
\end{se}

\begin{prop}\label{s20}
There exists an inclusion preserving
bijection 
\begin{align}
\{\text{lattices }T\subset\Lambda_k\mid
pT\d\overset{j}{\subset}T\subset\Lambda_k\} 
&\rightarrow\{k\text{-subspaces }U\subset V_k\mid \dim
U=\frac{l+j}{2},\ U^\perp\subset U\}\notag\\ 
T&\mapsto \overline{T},\notag
\end{align}
where $\overline{T}$ is equal to $T/p\Lambda_k\d$. 

In particular, the $\tau$-invariant lattices on the left hand side correspond
to the $\tau$-invariant subspaces on the right hand side.
\end{prop}

\begin{proof}
For a lattice $T$ contained in the set on the left hand side, we
obtain from (\ref{es32}) 
\begin{align}
pT\subset p\Lambda_k\overset{n-l}{\subset} p\Lambda_k\d\subset
pT\d\overset{j}{\subset}T\subset\Lambda_k.\label{es33} 
\end{align}
Since $\overline{pT\d}=\overline{T}^\perp$, we obtain from (\ref{es33})
\begin{align*}
\{0\}\subset \overline{T}^\perp\overset{j}{\subset}\overline{T}\subset
V_k 
\end{align*}
which proves the claim.
\end{proof}

\begin{cor}\label{s4}
Let $\Lambda\in\L_0^{(l)}$.
\begin{enumerate}
\item[a)]
The set of lattices $\Lambda_1\in\L_0^{(l_1)}$ with $\Lambda_1\subset\Lambda$
is equal to the set of $\F_{p^2}$-subspaces $U\subset V$ of dimension
$(l+l_1)/2$ with $U^\perp\subset U$.  

In particular, superspecial points in $\V(\Lambda)(k)$ correspond to
subspaces $U\subset V$ of dimension $(l+1)/2$ with
$U^\perp\subset U$. 
\item[b)] The lattices
$\Lambda_1\in\L_0^{(l_1)}$ with $\Lambda\subset\Lambda_1$ correspond
to the $\F_{p^2}$-subspaces $U\subset V'$ of dimension
$n-(l+l_1)/2$ with $U^{\perp'}\subset U$.
\end{enumerate}
\end{cor}

\begin{proof}
Part (a) follows from Proposition \ref{s20}. The superspecial points
are lattices in $\L_0$ of type $1$.  

To prove (b) 
let $\Lambda_1$ be an element of $\L_0^{(l_1)}$ with
$\Lambda\subset\Lambda_1$. We have 
\begin{align}
\Lambda\subset\Lambda_1\overset{n-l_1}{\subset}\Lambda_1\d\subset
\Lambda\d.\label{es19}  
\end{align}
For a lattice $L\subset C$ the dual $L^{\vee'}$ with respect to
$p\{\cdot,\cdot\}$ is equal to $p^{-1}L\d$. Thus (\ref{es19}) is equivalent to 
\begin{align*}
\Lambda=p(\Lambda\d)^{\vee'}\subset
p(\Lambda_1\d)^{\vee'}\overset{n-l_1}{\subset}\Lambda_1\d\subset\Lambda\d. 
\end{align*}
Now (b) follows from
Proposition \ref{s20} with $V'$ instead of $V$ and $p\{\cdot,\cdot\}$
instead of $\{\cdot,\cdot\}$ as the dimension of $V'$ is equal to
$n-l$.
\end{proof}  

\begin{se}
{\it Proof of Proposition \ref{s12}.}\label{s21}\\
We first prove the proposition in the case $i=0$.

Let $\Lambda$ be an element of $\L_0$ of type $l$. Let
$V$ be as in \ref{s19}. By Corollary \ref{s4} the superspecial 
points of $\V(\Lambda)(k)$ correspond to $\F_{p^2}$-subspaces $U\subset V$ of
dimension $(l+1)/2$ with $U^\perp\subset U$.  Such subspaces
always exist, hence $\V(\Lambda)(k)$ always contains superspecial
points. Furthermore, if $\dim V=l\geq 3$, there exists more than one
subspace with these properties. This shows that $\V(\Lambda)(k)$
consists of only one element if and only if $l=1$, hence 
proves a). 

To prove b), let $\Lambda'$ and $\Lambda$ be two elements of $\L_0$ such that
$\V(\Lambda')(k)\subset\V(\Lambda)(k)$. We want to prove that
$\Lambda'\subset\Lambda$. First note that
$\V(\Lambda')(k)=\V(\Lambda')(k)\cap \V(\Lambda)(k)$ is not
empty. Hence by Proposition \ref{s3} c)(ii), we obtain
$\V(\Lambda')(k)=\V(\Lambda\cap\Lambda')(k)$ and
$\Lambda\cap\Lambda'\in\L_0$. Therefore, it is sufficient to prove
that 
for $\Lambda'\subsetneqq\Lambda$ 
the set $\V(\Lambda')(k)$ is strictly contained in $\V(\Lambda)(k)$. 

Let $V$ be as in \ref{s19} and let $V_1$ be the subspace of $V$
corresponding to $\Lambda'$ (Cor.~\ref{s4}). Since $V_1\subsetneqq V$,
there exists a subspace $U\nsubseteq V_1$ of $V$ with
$U^\perp\overset{1}{\subset} U$. Thus there exists an element of
$\V(\Lambda)(k)\setminus\V(\Lambda')(k)$ (Prop.~\ref{s20}). 

Part c) follows from Corollary \ref{s4}.

By Remark~\ref{s18} d) 
the case of arbitrary $i$ follows from the case $i=0$.
\hfill{$\Box$}
\end{se}

\begin{se}\label{s2}
Let $\Lambda\in\L_0$ be of type $l$ and let 
$$d=\frac{l-1}{2}.$$ 
Let $V=\Lambda/p\Lambda\d$ as in \ref{s19}. By Proposition \ref{s20} we have
\begin{align*}
\V(\Lambda)(k)=\{U\subset V_k\mid \text{ dim }U=d+1,\
U^\perp\subset U\}. 
\end{align*}
Denote by $\Grass_{d+1}(V)$ the Grassmannian over
$\F_{p^2}$ of 
$(d+1)$-dimensional subspaces of $V$. 
The set $\V(\Lambda)(k)$ can naturally be endowed with the structure
of a closed 
subscheme of 
$\Grass_{d+1}(V)$. For
every $\F_{p^2}$-algebra $R$,
let $V_R$ be the base change $V\otimes_{\F_{p^2}}R$. By abuse 
of notation, we denote again by $\sigma$ the Frobenius on $R$. Let $U$
be a locally direct summand of $V_R$ of rank $m$. We define the dual
module $U^\perp\subset V_R$ as follows. For an $R$-module 
$M$, let $M^{(p)}=M\otimes_{R,\sigma}R$ be the Frobenius twist of $M$
and let $M^*=\Hom_R(M,R)$. 
Then $(\cdot,\cdot)$ induces
an $R$-linear isomorphism
\begin{align*}
\phi:(V_R)^{(p)}\stackrel{\sim}{\longrightarrow}(V_R)^*.
\end{align*}
Thus $\phi(U^{(p)})$ is a locally direct summand of $(V_R)^*$ of rank
 $m$. Let $\psi_U$ be the composition
 $V_R\stackrel{\sim}{\longrightarrow}(V_R)^{**}\twoheadrightarrow
 \phi(U^{(p)})^{*}$.     
 The orthogonal complement 
\begin{align}
U^\perp:=\ker(\psi_U)\label{es11}
\end{align}
is a locally direct summand of $V_R$ of rank $l-m$. Over an algebraically closed field, this 
definition coincides with the usual definition. 

We denote by $Y_\Lambda$ the closed
subscheme of $\Grass_{d+1}(V)$ defined over
$\F_{p^2}$ given by  
\begin{align}
Y_\Lambda(R)=\{U\subset V_R\text{ a locally direct summand}\mid\rk_R
U=d+1,\ U^\perp\subset U\}\label{es7} 
\end{align}
for every $\F_{p^2}$-algebra $R$.  
Note that for every algebraically closed field extension $k$ of $\overline{\F}_p$ we obtain
$Y_\Lambda(k)=\V(\Lambda)(k)$. By abuse of notation we will again
denote by $Y_\Lambda$ the 
corresponding scheme over $\bF_p$. 
Since there exists only one skew-hermitian form on $V$ up to
isomorphism, $Y_\Lambda$ depends up to isomorphism only on the
dimension of $V$, i.e., the type of $\Lambda$.   
We will show that $Y_\Lambda$ is a smooth irreducible variety over $\F_{p^2}$.
\end{se}

\begin{rem}
Let $\Lambda\in\L_i^{(l)}$.  By Remark~\ref{s18} d) the lattice 
$\Lambda'=\Psi_i(\Lambda)$ is contained in $\L_0^{(l)}$ and $\Psi_i$
induces a bijection  between $\V(\Lambda)(k)$ and
$\V(\Lambda')(k)$. Thus by \ref{s2} the set $\V(\Lambda)(k)$ can be
identified with the $k$-rational points of the variety
$Y_{\Lambda'}$. We set $Y_\Lambda=Y_{\Lambda'}$.
\end{rem}

\begin{se}\label{s7} 
For every $\F_p$-algebra $R$ we define a form 
$\langle\cdot,\cdot\rangle$ on $V\otimes_{\F_p}R$ by 
\begin{align*}
\langle\cdot,\cdot\rangle:(V\otimes_{\F_p}R)\times (V\otimes_{\F_p}R)
&\rightarrow
\F_{p^2}\otimes_{\F_p}R\\
(v\otimes a, w\otimes b)&\mapsto (v,w)\otimes ab, 
\end{align*}
where $(\cdot,\cdot)$ is the skew-hermitian form defined in \ref{s19}.  
The form $\langle\cdot,\cdot\rangle$ is linear in $R$ in both components where
the form $(\cdot,\cdot)_R$, extended to $V\otimes_{\F_{p^2}}R$ 
analogous to \ref{s19}, is 
linear in $R$ in the first and $\sigma$-linear in $R$ in the second component.
Let $G$ be the unitary group over $\F_p$ of \ref{s19}. Then for every 
$\F_p$-algebra $R$ the set $G(R)$ is defined as
$$G(R)=\{g\in\GL(V\otimes_{\F_p}R)\mid \langle gx,gy\rangle = \langle x,y
\rangle \text{ for all }x,y\in V\otimes_{\F_p}R\}.$$
Let T be the matrix 
$$T=\begin{pmatrix}
 &&1\\
&\a&\\
1&&\\
\end{pmatrix}\in\GL_l(\F_{p})$$ 
and let 
$\overline{t}$ be an element of $\F_{p^2}^\times$ with
$\overline{t}^p=-\overline{t}$. As there exists 
only one perfect skew-hermitian form on $V$ up to isomorphism, we can choose
an $\F_{p^2}$-basis $e_1,...,e_l$ of $V$ such that the form $(\cdot,\cdot)$  
with
respect to this basis is 
given by the matrix $\overline{t}T$. 
We identify $V$ via this basis with $(\F_{p^2})^l$.  
\end{se}

\begin{lem}\label{s24}
Let $G_{\overline{\F}_p}$ be the base change of $G$ over $\overline{\F}_p$. Then $G_{\overline{\F}_p}$ is isomorphic to $\GL_{l,\overline{\F}_p}$ with
Frobenius action given by 
$$\Phi(h)=T{^t}\!(h^{(p)})^{-1}T$$
for an element $h\in \GL_{l,\overline{\F}_p}$.
\end{lem}

\begin{proof}
Let $\sigma$ be the Frobenius of $\overline{\F}_p$. Furthermore, denote by 
$V_{\overline{\F}_p,\id}$ the vector space $V\otimes_{\F_{p^2}}\overline{\F}_p$
and denote by $V_{\overline{\F}_p,\sigma}$ the vektor space 
$V\otimes_{\F_{p^2}}\overline{\F}_p$, where in the second case the morphism 
$\F_{p^2}\hookrightarrow\overline{\F}_p$ is given by the Frobenius $\sigma$. 
We obtain an isomorphism
\begin{align*}
V\otimes_{\F_p}\overline{\F}_p&\stackrel{\sim}{\longrightarrow}
V_{\overline{\F}_p,\id}\oplus V_{\overline{\F}_p,\sigma}\\
v\otimes a &\mapsto (v\otimes a, v\otimes a).
\end{align*}
Both $V_{\overline{\F}_p,\id}$ and $V_{\overline{\F}_p,\sigma}$ are totally 
isotropic with respect to the form $\langle\cdot,\cdot\rangle$
and the form defines a perfect $\overline{\F}_p$-linear pairing between 
the two 
$\overline{\F}_p$-vektor spaces $V_{\overline{\F}_p,\id}$ and 
$V_{\overline{\F}_p,\sigma}$.

We now identify both $V_{\overline{\F}_p,\id}$ and 
$V_{\overline{\F}_p,\sigma}$ with $\overline{\F}_p^l$ using the 
$\F_{p^2}$-basis $e_1,...,e_l$ of $V$ of \ref{s7}. Then each $g\in 
G(\overline{\F}_p)$ corresponds to a pair $(g_{\id}, g_\sigma)\in
\GL_l(\overline{\F}_p)\times\GL_l(\overline{\F}_p)$ with 
\begin{align}
{^t}\!g_{\id} T g_\sigma = T.\label{es36}
\end{align} 
In particular, we obtain an isomorphism
\begin{align}
G_{\overline{\F}_p}&\rightarrow \GL_{l,\overline{\F}_p}\label{es35}\\
g&\mapsto g_{\id}.\notag
\end{align}

For an element $g\in G(\overline{\F}_p)$ the Frobenius $\Phi$ of $G$ is 
defined as $(\id_V\otimes\sigma)\circ g \circ(\id_V\otimes\sigma)^{-1}$. It
is easy to see that for the corresponding $(g_{\id}, g_\sigma)$ we obtain
$$\Phi(g_{\id},g_\sigma)=(g_\sigma^{(p)},g_{\id}^{(p)}),$$ 
where we denote for a matrix 
$h\in \GL_{l,\overline{\F}_p}$ by $h^{(p)}$ the matrix obtained by application 
of the Frobenius $\sigma$ to every entry of $h$.
Thus by (\ref{es36}) and (\ref{es35}) the Frobenius $\Phi$ 
may be identified with the morphism 
\begin{align*}
\GL_{l,\overline{\F}_p}&\rightarrow \GL_{l,\overline{\F}_p}\\
h&\mapsto T{^t}\!(h^{(p)})^{-1}T,
\end{align*}
which proves the claim.
\end{proof}

\begin{se}\label{s23}
The diagonal torus $S$ and the standard Borel $B$ of upper triangular
matrices in $G_{\overline{\F}_p}$ are defined over $\F_p$.
Note that $T$ 
corresponds to the longest Weyl group element of $G$ with respect to
$S$. 
Let  $I=(i_1,...,i_c)$ be an ordered partition of $l$. Denote by $P_I$ the
standard parabolic subgroup of $G_{\overline{\F}_p}$ containing  
$B$ such that
$\GL_{i_1,\overline{\F}_p}\times...\times\GL_{i_c,\overline{\F}_p}$ is
the Levi subgroup of $P$ containing $S$.   
We write $\Phi(I)$ for the partition ${(i_c,...,i_1)}$ and
obtain $\Phi(P_{(i_1,...,i_c)})=P_{\Phi(I)}$. Hence $\Phi$
induces a morphism   
$$\Phi:G_{\overline{\F}_p}/P_I\rightarrow
G_{\overline{\F}_p}/P_{\Phi(I)}.$$
To simplify notation, we write $G$ instead of $G_{\overline{\F}_p}$.
For a flag $\calF\in (G/P_I)(R)$, the dual flag $\calF^\perp$ with
respect to $(\cdot,\cdot)$ 
is defined analogous to (\ref{es11}). 
\end{se}

\begin{lem}\label{s6} 
Let $\calF$ be a flag in $G/P_I$.
The Frobenius $\Phi$ and the duality morphism
$\calF\mapsto\calF^\perp$ define the same morphism 
$G/P_I\rightarrow G/P_{\Phi(I)}$, i.e.,
the dual flag $\calF^\perp$ is equal to $\Phi(\calF)$. 
\end{lem}

\begin{proof}
Let $\calF$ be a flag in $G/P_I(R)$. Without loss of generality, we may
assume that the constituents of $\calF$ are free.  
Let $\calF_{I,0}$ be the standard flag corresponding to the parabolic
subgroup $P_I$. 
Then there exists an element $g\in G(R)$ such that
$\calF=g\calF_{I,0}$. We have  
$\Phi(gP_I)=\Phi(g)P_{\Phi(I)}$, hence 
$$\Phi(\calF)=\Phi(g)\calF_{\Phi(I),0}=\Phi(g)(\calF_{I,0})^\perp.$$ 
Furthermore, 
$(\Phi(g)x,gy)=(x,y)$ for all $x,y\in R^l$, hence
$\Phi(g)(\calF_{I,0})^\perp=(g\calF_{I,0})^\perp$ which proves the
claim. 
\end{proof}

\begin{prop}\label{s1}
The closed subscheme $Y_\Lambda$ of
$\Grass_{d+1}(V)$ is  smooth of dimension $d$. 
\end{prop}

\begin{proof}
We use the notation of \ref{s23}. It is sufficient to prove the claim
after base change to $\bF_p$. By Lemma \ref{s6} it is clear
that $Y_\Lambda$ is the intersection of the graph of the Frobenius 
\begin{align*}
G/P_{(d+1,d)}&\hookrightarrow G/P_{(d+1,d)}\times G/P_{(d,d+1)}\\
g&\mapsto (g,\Phi(g))
\end{align*}
with the morphism
\begin{align*}
G/P_{(d,1,d)}&\hookrightarrow G/P_{(d+1,d)}\times G/P_{(d,d+1)}\\
g&\mapsto (g,g).
\end{align*}
It is easy to see that this intersection is transversal of dimension $d$, hence
$Y_\Lambda$ is smooth.
\end{proof}

\begin{se}
Our next goal is to relate $Y_\Lambda$ with generalized
Deligne-Lusztig varieties.
Let $S_l$ be the symmetric group in $l$ elements and
denote by $W_I$ the Weyl group of the Levi subgroup
$\GL_{i_1,\overline{\F}_p}\times...\times\GL_{i_c,\overline{\F}_p}$ of
$P_I$. We identify 
$W_I$ with the corresponding subgroup of $S_l$. 
Let 
\begin{align*} 
\calF&=[0\subsetneqq\calF_1\subsetneqq...\subsetneqq\calF_c\subsetneqq V]\\
\G&=[0\subsetneqq\G_1\subsetneqq...\subsetneqq\G_c\subsetneqq V]
\end{align*}
be two flags  in $G/P_I(R)$. Then $\calF$ and $\G$ are in standard
position if all submodules $\calF_i+\G_j$  are locally direct summands
of $R^l$. Equivalently, the stabilizers of $\calF$ and $\G$ contain a
common maximal torus. 
For $\calF$ and $\G$ in standard position, we recall the definition of the
relative position $\inv(\calF,\G)$. Consider the diagonal action of $G$ on
$G/P_I\times G/P_I$. By the Bruhat decomposition, the orbits of this action
are classified by the quotient $W_{P_I}\backslash
S_l/W_{P_I}$. Then $\inv(\calF,\G)$ is defined as the element of the
constant sheaf  
associated to $W_{P_I}\backslash
S_l/W_{P_I}$ such that locally on $R$ the element $(\calF,\G)\in
G/P_I\times G/P_I$ is 
contained in the corresponding orbit.

For a partition $I$ with $I=\Phi(I)$ and an element $w\in
W_{P_I}\backslash 
S_l/W_{P_I}$,  we recall the generalized
Deligne-Lusztig variety $X_{P_I}(w)$ over $\overline{\F}_p$,  
\begin{align}
X_{P_I}(w)=\{\calF\in G/P_I\mid \inv(\Phi(\calF),\calF)=w\}.\label{es12}
\end{align}
We call this variety a generalized Deligne-Lusztig variety as in
\cite{DL} this variety is only defined in the case of a Borel
subgroup. As noted in \cite{Lu2}, this construction works for
arbitrary parabolic subgroups defined over $\F_p$.

The variety $X_{P_I}(w)$ is the transversal intersection of the graph of the 
Frobenius with the orbit of $(1,w)$ in $G/P_I\times G/P_I$ under the
diagonal action of $G$. Therefore, the variety $X_{P_I}(w)$ is  
smooth and its dimension is equal to the dimension of the subvariety
$P_IwP_I/P_I$ of $G/P_I$. In particular, if $P_I=B$, the dimension is
equal to the length of $w$. 

Denote for $i=0,...,d$ by 
$$w_i=(d+1,...,d+i+1)\in S_l$$ 
the cycle which
maps $d+1$ to $d+2$ etc. and denote by $I_i$ the partition
$(d-i,1^{2i+1},d-i)$. 
To simplify notation, we
write $P_i$ and $W_i$ instead of $P_{I_i}$ and $W_{P_i}$.

We have $w_0=\id$ and $I_0=(d,1,d)$. 
If $i=d$, we
obtain $w_d=(d+1,d+2,...,l-1,l)$ and $I_d=(1^l)$, i.e., $P_d=B$. 
Note that $I_{d-1}=I_d=(1^l)$ but $w_{d-1}\neq w_d$. 
\end{se}

\begin{lem}\label{s13}
For $0\leq i\leq d$
the dimension of $X_{P_i}(w_i)$ is equal to $i$.
\end{lem}

\begin{proof}
We have
$$\dim P_iw_iP_i/P_i=\dim P_iw_iP_i/B-\dim P_i/B.$$ 
Since $I_i$ is equal to $(d-i,1^{2i+1},d-i)$, the Levi subgroup $L$ as
in \ref{s23} of
$P_i$ is isomorphic to $\GL_{d-i}\times
\mathbb{G}_m^{2i+1}\times\GL_{d-i}$. Let $B_{d-i}$ be the standard
Borel of upper triangular matrices 
in $\GL_{d-i}$. Then
\begin{align*}
\dim P_i/B&=\dim L/(B\cap L)\\
&=2\dim\GL_{d-i}/B_{d-i}\\
&=(d-i-1)(d-i).  
\end{align*}

On the other hand, $P_iw_iP_i=BW_iw_iW_iB$. Let $w_i'\in S_l$ be the
longest representative of the residue class of $w_i$ in
$W_i\backslash S_l/W_i$. Then the dimension of
$P_iw_iP_i/B$ is equal to the length of $w_i'$. 

Suppose first that
$i=d$, i.e., $P_d=B$. 
We obtain $W_d=1$. Thus $w_d'=w_d$ and 
$$\dim X_B(w_d)=\text{length}(w_d)=d.$$ 
Now suppose  $i<d$. Then
\begin{align*}
W_i&=\langle (12),...,(d-i-1,d-i),(d+i+2,d+i+3),...,(n-1,n)\rangle\\
&\cong S_{d-i}\times S_{d-i}.
\end{align*}
Since $w_i=(d+1,...,d+i+1)$, we obtain that $w_i'$ is equal to the
product of $w_i$ with the longest element in $W_i$. Therefore, the
length of $w_i'$ is equal to $i+(n-i)(n-i-1)$ which proves the
claim.
\end{proof} 

The next theorem establishes a link between $Y_\Lambda$ and
generalized Deligne-Lusztig varieties.
\begin{thm}\label{s8}
There exists a decomposition of $Y_\Lambda$ over $\bF_p$ into a
disjoint union of 
locally closed subvarieties 
\begin{align}
Y_\Lambda=\biguplus_{i=0}^{d}X_{P_i}(w_i),\label{es21}
\end{align}
such that for every $j$ with $0\leq j\leq d$ the subset
$\biguplus_{i=0}^{j}X_{P_i}(w_i)$
is closed in $Y_\Lambda$.

The variety $X_B(w_d)$ is open, dense and irreducible of
dimension $d$ 
in $Y_\Lambda$. In particular, $Y_\Lambda$ is irreducible of dimension
$d$. 
\end{thm}

\begin{proof} 
 As in \ref{s7} we identify $V$ with $(\F_{p^2})^l$ such that
 the form $(\cdot,\cdot)$ is given by the matrix $tT$.  
For an algebraically closed field extension $k$ of $\overline{\F}_p$ and  a subspace
 $U\subset V_k$, we have  
$(U^\perp)^\perp=\tau(U)$ by \ref{s19}. 
For an arbitrary $\overline{\F}_p$-algebra $R$, we do not have an
 operator $\tau$ on $V_R$, but
to simplify notation, we write
 $\tau(U)$ for $(U^\perp)^\perp$ for all locally direct summands $U$
 of $V_R$.  

Let $U$ be an element of $Y_\Lambda(R)$. We may
assume that $\Spec R$ is connected.
The module $U$ is a locally direct summand of $R^l$ of rank $d+1$ such
that $U^\perp\subset U$ with $U/U^\perp$ locally free of rank 1. Let
$\calF$ be the flag   
$$\calF=[0\subsetneqq U^\perp\subsetneqq U\subsetneqq V].$$
By Lemma \ref{s6} the flag $\Phi(\calF)$ is equal to 
$$\Phi(\calF)=[0\subsetneqq U^\perp\subsetneqq \tau(U)\subsetneqq V].$$ 
If $U+\tau(U)$ is a locally direct summand of $R^l$, i.e., if $\calF$ and
$\Phi(\calF)$  
are in standard position, 
we obtain that $\inv(\Phi(\calF),\calF)$ is either the identity or
$w_1$ in $W_0\backslash S_l/W_0$. Hence we
obtain a disjoint decomposition of $\V(\Lambda)$ into the open
subvariety $X_{P_0}(w_1)$ and the closed subvariety 
$X_{P_0}(\id)$
\begin{align}
Y_\Lambda=X_{P_0}(\id)\uplus X_{P_0}(w_1).\label{es20}
\end{align}

Let $i$ be an integer with $0\leq i\leq d-1$ and let $\calF$ be a flag in
$X_{P_i}(w_i)\uplus X_{P_i}(w_{i+1})$. Then $\calF$ is of the form
\begin{align*}
\calF=[0\subset\calF_{-i-1}\subsetneqq...\subsetneqq\calF_{-1}\subsetneqq
\calF_1\subsetneqq...\subsetneqq\calF_{i+1}\subset
V]. 
\end{align*}
For an integer $j$ with $1\leq j\leq i+1$, the modules $\calF_{-j}$ and
$\calF_j$ are 
locally direct summands of $R^l$ of rank $(d-j+1)$ and $(d+j)$
respectively.   

{\it Claim:} The flag $\calF$ is uniquely determined by $\calF_1$. \\
Indeed, since $\Phi(\calF)$ and $\calF$ are in relative position $w_i$
or $w_{i+1}$, we obtain a commutative diagram
$$
\xymatrix{
\calF_{i+1}^\perp\ar@{=}[d]\ar@{}[r]|{\subsetneqq}&\hdots\ar@{}[r]|{\subsetneqq}&\calF_{1}^\perp\ar@{}[r]|{\subsetneqq}\ar@{=}[d]&\calF_{-1}^\perp\ar@{}[r]|{\subsetneqq}\ar@{}[d]|{\neq}\ar@{^{(}->}[dr]&\calF_{-2}^\perp\ar@{}[r]|{\subsetneqq}\ar@{}[d]|{\neq}&\hdots\ar@{}[r]|{\subsetneqq}&\calF_{-i}^\perp
\ar@{}[r]|{\subsetneqq}\ar@{}[d]|{\neq}\ar@{^{(}->}[dr]&\calF_{-i-1}^\perp\\
\calF_{-i-1}\ar@{}[r]|{\subsetneqq}&\hdots\ar@{}[r]|{\subsetneqq}&\calF_{-1}\ar@{}[r]|{\subsetneqq}&\calF_1\ar@{}[r]|{\subsetneqq}&\calF_2\ar@{}[r]|{\subsetneqq}&\hdots\ar@{}[r]|{\subsetneqq}&\calF_i
\ar@{}[r]|{\subsetneqq}&\calF_{i+1}.\\
}
$$
Furthermore, we have
\begin{align*}
\calF_{-i-1}^\perp=\calF_{i+1}&\qquad\text{ if }\inv(\Phi(\calF),\calF)=w_i,\\ 
\calF_{-i-1}^\perp\neq\calF_{i+1}&\qquad\text{ if
}\inv(\Phi(\calF),\calF)=w_{i+1}.
\end{align*} 
As $\Phi(\calF)$ and $\calF$ are in standard position, we obtain that
for $j$ with $2\leq j\leq i+1$
\begin{align*}
\calF_j&=\calF_{j-1}+\calF_{-j+1}^\perp\\
&=\calF_{j-1}+\tau(\calF_{j-1}).\notag
\end{align*}
By induction we obtain
$$\calF_j=\sum_{m=1}^{j}\tau^m(\calF_1).$$
Therefore, $\calF$ is uniquely determined by $\calF_1$ and the claim is proved.

Let $i$ with $0\leq i\leq d-2$. The claim shows that 
the natural morphisms 
\begin{align*}
X_{P_{i+1}}(w_{i+1})&\rightarrow X_{P_i}(w_{i+1}),\\  
X_{P_{i+1}}(w_{i+2})&\rightarrow X_{P_i}(w_{i+1})
\end{align*}
are monomorphisms.
We now want to show that 
we have a decomposition of $X_{P_i}(w_{i+1})$ into the  closed
subvariety $X_{P_{i+1}}(w_{i+1})$ and the open subvariety
$X_{P_{i+1}}(w_{i+2})$, 
\begin{align}
X_{P_i}(w_{i+1})=X_{P_{i+1}}(w_{i+1})\uplus X_{P_{i+1}}(w_{i+2}).\label{es13}
\end{align}
Indeed, let $\calF$ be a flag in $X_{P_i}(w_{i+1})(R)$.
Since $\Phi(\calF)$ and $\calF$ are in standard position,
$\calF_{i+2}:=\calF_{i+1}+\tau(\calF_{i+1})$ and
$\calF_{-i-2}:=\calF_{i+2}^\perp$ are locally direct summands of $R^l$
of rank $(d+i+2)$ and $(d-i-1)$ respectively.  
We extend $\calF$ to the flag
$$\calF'=[0\subset\calF_{-i-2}\subsetneqq...\subsetneqq\calF_{-1}\subsetneqq
\calF_1\subsetneqq...\subsetneqq\calF_{i+2}\subset V].$$
If $\Phi(\calF')$ and $\calF'$ are in standard position, we obtain
that $\inv(\Phi(\calF'),\calF')$ is either $w_{i+1}$ or $w_{i+2}$ in
$W_{i+1}\backslash S_l/W_{i+1}$. Thus
$\calF'$ is either an element of $X_{P_{i+1}}(w_{i+1})$ or of
$X_{P_{i+1}}(w_{i+2})$. This proves (\ref{es13}). 

The disjoint sum decomposition in (\ref{es21}) follows by
induction from 
(\ref{es20}) and (\ref{es13}).  
The subset
$\biguplus_{i=0}^{j}X_{P_i}(w_i)$
is closed in $Y_\Lambda$ by construction. 
By Lemma~\ref{s13} the variety $X_{P_i}(w_i)$ is of
dimension $i$. Since 
$Y_\Lambda$ is smooth of dimension $d$, the open subvariety
$X_B(w_d)$ is dense. It remains to show that $X_B(w_d)$ is irreducible.

Consider the set
$D=\{(1,2),...,(l-1,l)\}$  of simple reflections of $G$ with respect
to $B$. The Frobenius action on $D$ with respect to
$G$ is given by 
$\sigma((i,i+1))=(l-i,l-i+1)$. Since $l$ is odd, there are $d$ orbits of
the action of $\sigma$ on $D$. Each orbit is of the form
$\{(i,i+1),(l-i,l-i+1)\}$ with $d+1\leq i\leq l-1$. A Coxeter element of
$S_l$ is a product of the form
$\nu_1\cdots\nu_d\in S_l$, where $\nu_1,...,\nu_d\in D$ are
representatives of the different orbits
(\cite{Lu1}). The element
$w_d=(d+1,d+2)\cdots (l-1,l)\in S_l$ is a Coxeter element. It is shown
in (\cite{Lu1} 
4.8) that $X_B(w)$ is irreducible for each Coxeter element $w$. In
particular, $X_B(w_d)$ is irreducible.  
This proves the theorem.
\end{proof}

\begin{prop}\label{s11}
Let $\Lambda'\in\L_0$ with $\Lambda'\subset\Lambda$ and denote by $V'$
the $\F_{p^2}$-vector space $\Lambda'/p(\Lambda')\d$. Let $l'$ be the 
type of $\Lambda'$ and let $d'=(l'-1)/2$. Then 
the variety $Y_{\Lambda'}$ is a closed subvariety of $Y_\Lambda$ of
dimension equal to $d'$. 
\end{prop}

\begin{proof}
By definition (\ref{es7}) we have
\begin{align}
Y_{\Lambda'}(R)=\{U\subset V'_R\text{ a locally direct
  summand}\mid\rk_R U=d'+1,\ U^\perp\subset U\},\notag
\end{align}
hence by \ref{s20}
\begin{align}
Y_{\Lambda'}(R)
=\{U\in Y_\Lambda(R)\mid U\subset (\Lambda'/p\Lambda\d)_R\}.\notag
\end{align}
Thus $Y_{\Lambda'}$ is a closed subvariety of $Y_\Lambda$.
By Proposition \ref{s1} the dimension of $Y_{\Lambda'}$ is equal to
$d'$.
\end{proof}

\begin{se}
For $U\in Y_\Lambda(k)$ let $0\leq i_U\leq d$ be the minimal integer such
that $U+...+\tau^i(U)$ is
$\tau$-invariant. For a lattice $A\in\V(\Lambda)(k)$, let $i_A$ be the
minimal integer such that $A+...+\tau^i(A)$ is $\tau$-invariant. 
Suppose that $A$ corresponds to $U$ via the bijection of Proposition \ref{s20},
then $i_A=i_U$.

Denote by $\V(\Lambda)(k)^\circ\subset\V(\Lambda)(k)$  the set of
lattices $A\in\V(\Lambda)(k)$  
 with $i_A=d$, i.e., $A+...+\tau^d(A)=\Lambda_k$.
\end{se}

\begin{cor}
For every integer $i$ with $0\leq i\leq d$ we have
\begin{align*}
X_{P_i}(w_i)(k)&=\{U\in Y_\Lambda(k)\mid i_U=i\}\\
&=\{A\in\V(\Lambda)(k)\mid i_A=i\}.
\end{align*}
In particular,  $X_B(w_d)(k)=\V(\Lambda)(k)^\circ$. 

Hence 
\begin{align}
\biguplus_{j=0}^{i}X_{P_j}(w_j)(k)=\bigcup_{\substack{\Lambda'\subset\Lambda\\
\Lambda'\in\L_0^{(l')}\\
l'\leq 2i+1}}
\V(\Lambda')(k)
=\bigcup_{\substack{\Lambda'\subset\Lambda\\
\Lambda'\in\L_0^{(2i+1)}}}
\V(\Lambda')(k).
\end{align}
There exists a  scheme theoretical decomposition
\begin{align}
Y_\Lambda=\bigcup_{\substack{\Lambda'\subsetneqq\Lambda\\
\Lambda'\in\L_{0}}}
Y_{\Lambda'}\uplus X_B(w_d),\label{es17}
\end{align}
where the first summand is closed and the second summand is open.
\end{cor}

\begin{proof}
The corollary follows from the proof of Theorem \ref{s8} and from
Proposition \ref{s11}.
\end{proof}


\section{The combinatorial intersection behaviour of the sets $\V(\Lambda)$}
\label{i}

\begin{se}\label{i1}
Let $C$ be a $\Q_{p^2}$-vector space of dimension $n$ with perfect
skew-hermitian  form $\{\cdot,\cdot\}$
as in \ref{s16}. Denote by $t$ an element of $\Z_{p^2}^\times$ with
$t^\sigma=-t$. As in \ref{s16} we assume that the
skew-hermitian form $\{\cdot,\cdot\}$ is equivalent to the form  
induced by the matrix $T$ with  $T=tI_n$ if $n$ is odd and equivalent
to the form induced  
by
$T=tJ_n$ if $n$ is even. Denote by $H$ the special unitary group over
$\Q_p$ with respect to $(C,\{\cdot,\cdot\})$, i.e., for a
$\Q_p$-algebra $R$, we have 
\begin{align*}
H(R)=\{g\in\SL_{\Q_{p^2}\otimes_{\Q_p}R} (C\otimes_{\Q_p}R)\mid
\{ gx,gy\}=\{ x,y\}\text{ for all }x,y\in C\otimes_{\Q_p}R\}. 
\end{align*}
Let $k$ be an algebraically closed field extension of $\overline{\F}_p$. 
In this section we will prove that the incidence relation of the
sets $\V(\Lambda)(k)$ of $\D_i(C)(k)$
can be read off from the 
combinatorial simplicial structure of the Bruhat-Tits building $\B(H,\Q_p)$
associated to $H$. As in Section~\ref{s}, we assume that $ni$ is even.

Propositions \ref{s3} and \ref{s12} show that the intersection
behaviour of the sets $\V(\Lambda)(k)$ only depends on
$\Lambda\in\L_i$. Therefore, we write $\V(\Lambda)$ and $\D_i(C)$ instead of
$\V(\Lambda)(k)$ and $\D_i(C)(k)$. 
\end{se}

\begin{de}\label{i3}
Let $\B_i$ be the abstract simplicial complex given by the following
data. Let $m$ be a nonnegative integer. 
An $m$-simplex is a subset $S\subset\L_i$ of $m+1$ elements which
satisfies the following condition. There exists an ordering
$\Lambda_0,...,\Lambda_m$ of the elements of $S$ such that  
\begin{align}
p^{i+1}\Lambda_m\d\subsetneqq\Lambda_0\subsetneqq\Lambda_1\subsetneqq...\subsetneqq
\Lambda_m.\label{ei3} 
\end{align}
A vertex is defined as a 0-simplex.
\end{de}

\begin{rem}\label{i8}
Let $S$ be an $m$-simplex of $\B_i$ and let $\Lambda_0,...,\Lambda_m$ be
an ordering of the elements of $S$ which satisfies (\ref{ei3}).  
Since all $\Lambda_j$ are in $\L_i$, we obtain from (\ref{ei3}) the
more precise chain of inclusions
\begin{align}
p^{i+1}\Lambda_0\d\subsetneqq\Lambda_0\subsetneqq...\subsetneqq\Lambda_m\subset
p^i\Lambda_m\d \subsetneqq...\subsetneqq p^i\Lambda_0\d.\label{ei5} 
\end{align}
Obviously, we have $0\leq m\leq (n-1)/2$.
\end{rem}

\begin{se}\label{i10}
Let $\{\Lambda\}$ be a vertex of $\B_i$. 
We say that $\{\Lambda\}$ is of type $l$ if $\Lambda$ is of type $l$.
A vertex $\{\Lambda'\}\in\B_i$ is a neighbour of $\{\Lambda\}$ if
$\Lambda\neq\Lambda'$ and if there exists
a simplex $S\in\B_i$   
such that $\Lambda$ and $\Lambda'$ are contained in $S$. This is equivalent to
the condition that $\{\Lambda,\Lambda'\}$ is a 1-simplex of $\B_i$, i.e.,
$\Lambda\subsetneqq\Lambda'$ or $\Lambda'\subsetneqq\Lambda$.  
Let $l$ and $l'$ be the type of $\Lambda$ and  $\Lambda'$ respectively. In the
first case we obtain $l< l'$ and in the second case $l'>l$.
\end{se}

\begin{prop}\label{i4}
\begin{enumerate}
\item[a)]
Let $\Lambda$ be a lattice in $\L_i$ of type $l$. 
Then the set
$$\{\Lambda'\in\L_i\mid \V(\Lambda')\subsetneqq\V(\Lambda)\}$$
corresponds to the set of neighbours
of $\{\Lambda\}$ of type $l'< l$. 
\item[b)]
Let $\Lambda$ and $\Lambda'$ be in $\L_i$ of type $l$ and $l'$
respectively such that $\Lambda\neq\Lambda'$.  Then
$$\V(\Lambda)\cap\V(\Lambda')\neq\emptyset$$ 
if and only if $\{\Lambda\}$
is  a neighbour of $\{\Lambda'\}$ or if there exists
a vertex $\{\widetilde{\Lambda}\}$ in $\B_i$ of type $c<\min\{l,l'\}$ such
that $\{\widetilde{\Lambda}\}$  is a
common neighbour of $\{\Lambda\}$ and $\{\Lambda'\}$.
\end{enumerate}
\end{prop}

\begin{proof} 
The proposition follows from Proposition \ref{s3} c)(ii).
\end{proof}

\begin{thm}\label{i5}
Let $\B(H,\Q_p)_\text{simp}$ be the abstract simplicial complex of
the Bruhat-Tits 
building of $H$. Then there exists a simplicial bijection 
between $\B_i$ and
$\B(H,\Q_p)_\text{simp}$.  
\end{thm}

\begin{proof}
We choose a basis of $C$ such that the form $\{\cdot,\cdot\}$ is given by
$T$. For a matrix $g$, denote by $g^{(\sigma)}$ the
matrix obtained by applying the Frobenius $\sigma$ of $\Q_{p^2}/\Q_p$
to every entry of $g$.
Then $H$ is isomorphic to $\SL(C)$ over $\Q_{p^2}$ with
Frobenius 
$$\Phi(g)=T(^t\!g^{(\sigma)})^{-1}T$$
 for
$g\in\SL(C)$ (proof analogous to proof of \ref{s24}).  
The simplicial complex $\B(H,\Q_p)_\text{simp}$ is equal to the
fixed points of $\Phi$ on the  
simplicial complex $\B(\SL(C),\Q_{p^2})_\text{simp}$
 (\cite{Ti} 2.6.1).

An $m$-simplex of $\B(\SL(C),\Q_{p^2})_\text{simp}$ is a set
$\{[\Lambda_0],...,[\Lambda_m]\}$ of homothety classes of lattices 
with the following property. There
exist representatives
$\Lambda_j\in[\Lambda_j]$ which form after renumbering the lattices
$\Lambda_0,...,\Lambda_m$ an infinite
lattice chain
\begin{align}
...\subsetneqq
p\Lambda_m\subsetneqq\Lambda_0\subsetneqq...\subsetneqq\Lambda_m\subsetneqq
p^{-1}\Lambda_0\subsetneqq...\label{ei2} 
\end{align}
We first consider the case $i=0$.
We define a simplicial morphism
\begin{align}
\varphi:\B_0&\rightarrow\B(\SL(C),\Q_{p^2})_\text{simp}\label{ei6}\\
\{\Lambda_0,...,\Lambda_m\}&\mapsto\{[\Lambda_0],...,[\Lambda_m],[\Lambda_m\d]
...[\Lambda_0\d]\}.\notag
\end{align} 
This morphism is well defined as (\ref{ei2}) follows from
(\ref{ei5}). To show that $\varphi$ induces an isomorphism onto
$\B(H,\Q_p)_\text{simp}$, it is sufficient to prove the following claim.

{\it Claim:} Each simplex of $\B(H,\Q_p)_\text{simp}$ can be written
uniquely as 
$$\{[\Lambda_0],...,[\Lambda_a],[\Lambda_a\d],...,[\Lambda_0\d]\}$$ 
such that there exist representatives
$\Lambda_0,...,\Lambda_a$ satisfying
\begin{align}
p\Lambda_0\d\subsetneqq\Lambda_0\subsetneqq...\subsetneqq\Lambda_a\subset
\Lambda_a\d\subsetneqq...\subsetneqq\Lambda_0\d.\label{ei4} 
\end{align}

Indeed, a simplex 
$\{[\Lambda_0],...,[\Lambda_m]\}$ of $\B(\SL(C),\Q_{p^2})_\text{simp}$
is a fixed point under the action of 
$\Phi$ if the lattice chain
(\ref{ei2}) 
coincides with its dual chain
\begin{align*}
...\subsetneqq
p\Lambda_0\d\subsetneqq\Lambda_m\d\subsetneqq...\subsetneqq\Lambda_0\d\subsetneqq
p^{-1}\Lambda_m\d\subsetneqq... 
\end{align*}
This means that there exist integers $j$ and $a$ with $0\leq a\leq m$ such that
$$
\xymatrix{
p^j\Lambda_{a+1}\d\ar@{}[r]|{\subsetneqq}\ar@{=}[d]&p^j\Lambda_a\d\ar@{}[r]|{\subsetneqq}\ar@{=}[d]&\hdots\ar@{}[r]|{\subsetneqq}&p^j\Lambda_0\d\ar@{}[r]|{\subsetneqq}\ar@{=}[d]&p^{j-1}\Lambda_m\d\ar@{}[r]|{\subsetneqq}\ar@{=}[d]&\hdots\ar@{}[r]|{\subsetneqq}&p^{j-1}\Lambda_{a+1}\d\ar@{=}[d]\\
p\Lambda_m\ar@{}[r]|{\subsetneqq}&\Lambda_0\ar@{}[r]|{\subsetneqq}&\hdots\ar@{}[r]|{\subsetneqq}&\Lambda_a\ar@{}[r]|{\subsetneqq}&\Lambda_{a+1}\ar@{}[r]|{\subsetneqq}&\hdots\ar@{}[r]|{\subsetneqq}&\Lambda_m.\\
}
$$
An easy index calculation shows that there exists an ordering of the
homothety classes and representatives such that we obtain a lattice
chain 
\begin{align*}
p\Lambda_0\d\subsetneqq\Lambda_0\subsetneqq\Lambda_1\subsetneqq...\subsetneqq
\Lambda_{[\frac{m+1}{2}]} 
\subset\Lambda_{[\frac{m+1}{2}]}\d\subsetneqq...\subsetneqq\Lambda_1\d\subsetneqq
\Lambda_0\d.   
\end{align*}
This proves the claim in the case of $i=0$.

Now consider the general case. By Remark~\ref{s18} d) the isomorphism
$\Psi_i$ of Proposition~\ref{d21} induces an inclusion and type
preserving isomorphism of $\L_i$ with $\L_0$. Thus $\Psi_i$ induces a
simplicial bijection
\begin{align}
\theta_i:\B_i\stackrel{\sim}{\longrightarrow}\B_0\label{ei8}
\end{align}
which proves the theorem.
\end{proof}

\begin{prop}\label{i6}
Let $\Lambda,\Lambda'\in\L_i$. There exist a positive integer $u$ and elements
$\Lambda_1=\Lambda,\Lambda_2,...,\Lambda_{u-1},\Lambda_u=\Lambda'$ in $\L_i$
such that  
$$\V(\Lambda_j)\cap\V(\Lambda_{j+1})\neq\emptyset$$
for every $j$ with $1\leq j\leq u-1$. 
\end{prop}

\begin{proof}
The building $\B(H,\Q_p)$ is connected, hence the simplicial complex
$\B(H,\Q_p)_{simp}$ is connected.
Then the proposition follows from Theorem
\ref{i5} and Proposition~\ref{i4}.
\end{proof}

\begin{prop}\label{i7}
In the case of $\GU(1,2)$, i.e., if $n=3$, the simplicial complex $\B_i$
is a tree. It has two different kind of vertices. Vertices   of type 1
correspond to lattices in $\L_i^{(1)}$, i.e., superspecial lattices in
$\D_i(C)$, and have $p+1$ neighbours  of type 3. Vertices of type 3
correspond to lattices in $\L_i^{(3)}$ 
and have $p^3+1$ neighbours of type 1.
\end{prop}

\begin{proof}
As $n$ is equal to 3, Remark~\ref{i8} shows that there exist only
0-simplices and 1-simplices. The building $\B(H,\Q_p)$ is 
contractible, hence its simplicial complex is a tree. Thus by
Theorem~\ref{i5}, the simplicial complex $\B_i$ is a tree.
The type of a lattice $\Lambda\in\L_i$ is equal to 1 or 3 (Rem.~\ref{s18}).
By construction of $\B_i$, each vertex of type 1 has only neighbours of
type 3 and each vertex of 
type 3 has only neighbours of type 1. 

Now consider the case $i=0$.
Let $\Lambda$ be in $\L_0^{(3)}$. Denote by $V$ the $\F_{p^2}$-vector
space $V=\Lambda/p\Lambda$ with perfect skew-hermitian form $(\cdot,\cdot)$
as in \ref{s19}.  
By Corollary~\ref{s4} a) the neighbours of $\{\Lambda\}$ correspond to
$\F_{p^2}$-subspaces $U$ of $V$ of dimension 2  with
$U^\perp\subset U$. As all skew-hermitian forms on $V$ are equivalent, there
exists a basis of $V$ such that $(\cdot,\cdot)$ is given by the matrix
$\overline{t}I_3$. By duality the neighbours of $\{\Lambda\}$
correspond to the totally isotropic subspaces of $V$ of dimension 1,
i.e., to the $\F_{p^2}$-rational points of the 
Fermat curve $\C$ in $\P^2_{\F_{p^2}}$ given by the equation
$$x_0^{p+1}+x_1^{p+1}+x_2^{p+1}.$$
An easy calculation shows that the number of $\F_{p^2}$-rational
points of $\C$ is equal to $p^3+1$.

Analogously, the neighbours of a 0-simplex $\{M\}$ correspond by
Corollary \ref{s4} b) to the totally isotropic 1-dimensional subspaces
of the 2-dimensional $\F_{p^2}$-vector space $V'=M\d/M$. An easy 
calculation shows the claim.

Now consider the general case.
By \eqref{ei8} there exists a simplicial bijection between $\B_i$ and
$B_0$, hence  the claim follows from the case $i=0$.
\end{proof}

\begin{rem}\label{i9}
Let $n=3$ and let $\Lambda,\Lambda'\in\L_i$. As $\B_i$
is a tree, 
there exists a unique lattice chain for $\Lambda$ and $\Lambda'$ as in
Proposition~\ref{i6} of 
minimal length. We call its length $u$ the distance
$u(\Lambda,\Lambda')$ of  $\Lambda$ and $\Lambda'$.  The distance
of two lattices $\Lambda$ and $\Lambda'$ of the same type is even. 
\end{rem}



\section{The local structure of $\N^{red}$ for $\GU(1,2)$}
\label{g}

\begin{se}\label{g1}
In the next two sections, we will describe the scheme theoretic structure of
$\N^{red}$ in the case of the unitary group $\GU(1,2)$. We again assume
that $p\neq 2$. We will first
describe the scheme theoretic structure of the open and closed
subscheme $\N_0^{red}$ of quasi-isogenies of height 0.
\end{se}

\begin{se}\label{g17}
We will always denote by
$k$ an algebraically closed field extension of 
$\overline{\F}_p$.
By abuse of notation, we will mostly identify the elements of $\N(k)$
with their corresponding Dieudonn\'e modules as in \ref{d1}. We say
that a $p$-divisible group in $\N_0(k)$ is superspecial if the
corresponding Dieudonn\'e module is superspecial.
Let $N=N_0\oplus N_1$ be the
isocrystal over $W(\overline{\F}_p)_\Q$ with perfect alternating form
$\langle\cdot,\cdot\rangle$ as in Section \ref{d}. 
In Sections 2 and 3, we have denoted by $C$ the
$\Q_{p^2}$-vector space  $N_0^\tau$ of $\tau$-invariant elements of
$N_0$ (\ref{d4}). 
By abuse of notation, we will identify the sets $\N_0(k)$ and
$\D_0(C)(k)$ using the bijection of Proposition~\ref{d5}.
By Remark~\ref{s18} we will identify the
superspecial points of $\N_0(k)$ with the lattices in
$\L_0^{(1)}$ . 

We fix some notation. 
As in Lemma \ref{d11}, we fix an element $t\in\Z_{p^2}^\times$ with
$t^\sigma=-t$ and denote 
by $\overline{t}$ its image in $\F_{p^2}^\times$. 
For an element $x$ in a ring $R$, we write $[x]$ 
for the Teichm\"uller representative $(x,0,...)\in W(R)$. The map $[\
]:R\rightarrow W(R)$ is  multiplicative and injective.
As $p\neq 2$, we have $1+[-1]=0$ which we will frequently use in the sequel.

For every scheme $S$ over
$\F_{p^2}$, we denote again by $S$ the corresponding scheme over $\bF_p$. 
\end{se}

\begin{se}\label{g26}
For $\Lambda\in\L_0$ we have defined a
closed
subset $\V(\Lambda)(k)$ of $\N_0(k)$. 
By \ref{s2} there exists a smooth,
irreducible and proper variety $Y_\Lambda$ over $\F_{p^2}$ such that
$Y_\Lambda(k)=\V(\Lambda)(k)$ for every algebraically closed field
extension $k$ of 
$\overline{\F}_p$.  
The varieties $Y_\Lambda$ depend up to isomorphism  only on the type
$l$ of $\Lambda$ (\ref{s2}).    
By Remark~\ref{s18} the type $l$ is equal to 1 or 3. If $l=1$,  
the variety $Y_\Lambda$ consists of only one point.
If
$l=3$, the variety 
$Y_\Lambda$ is smooth and irreducible of dimension~1 (Prop.~\ref{s1},
Thm.~\ref{s8}). 

By Proposition~\ref{s3} we know that 
$$\N_0(k)=\bigcup_{\Lambda\in\L_0^{(3)}}\V(\Lambda)(k).$$  
If the intersection of two different sets $\V(\Lambda)(k)$ and
$\V(\Lambda')(k)$ is nonempty, they intersect at one superspecial
point, the lattice
$\Lambda\cap\Lambda'$
(Prop.~\ref{s3}, Prop.~\ref{s12}). 
Each set $\V(\Lambda)(k)$ contains $p^3+1$ superspecial points and
each superspecial point is contained in $p+1$ sets $\V(\Lambda)(k)$
(Prop.~\ref{i7}). 
\end{se}

Let $\M\in\N_0(k)$ be a superspecial Dieudonn\'e
lattice. 
 The following lemma follows directly from
Sections \ref{d} and \ref{s}. 
\begin{lem}\label{g3}
The lattice $\M$ is already defined over $\Z_{p^2}$. 
There exists a $\tau$-invariant basis $e_1,e_2,e_3,f_1,f_2,f_3$ of
$\M$ 
such that 
\begin{align*}
\M_0&=\langle e_1,e_2,e_3\rangle_{W(k)}\\
\M_1&=\langle f_1,f_2,f_3\rangle_{W(k)}.
\end{align*}
The matrix of $F$ with respect to the above basis of $\M$ is given by
\begin{align*}
\begin{pmatrix}
&&&1&&\\
&&&&p&\\
&&&&&1\\
p&&&&&\\
&1&&&&\\
&&p&&&\\
\end{pmatrix}.
\end{align*}
The form $\{\cdot,\cdot\}$ on $\M_0$ is given by the matrix
\begin{align}
t\begin{pmatrix}
p&&\\
&1&\\
&&p\\
\end{pmatrix}.\label{eg23}
\end{align}
In particular, the
only nonzero values of 
the alternating form $\langle\cdot,\cdot\rangle$ on the basis of $\M$
are given by 
$\langle e_i,f_j\rangle=-\langle f_i,e_j\rangle=t\delta_{ij}$. 
\end{lem}

\begin{proof}
As $\M$ is superspecial, it is $\tau$-invariant, hence defined over $\Z_{p^2}$.
By Lemma \ref{d18} there exists a basis $e_1,e_2,e_3$
of $\M_0^\tau$ such that the form $\{\cdot,\cdot\}$ is given by the matrix
(\ref{eg23}). 
We have $p\M_0\d=\langle e_1,pe_2,e_3\rangle$.
The proof of Proposition \ref{d5} shows that $\M_1$ is equal to
$F^{-1}(p\M_0\d)$. For $i=1,3$ we define $f_i=F^{-1}e_i$ and set
$f_2=F^{-1}(pe_2)$. This basis satisfies the conditions of the
lemma. 
\end{proof} 

\begin{se}  
Let $\J$ be the set 
$$\J=\{[\lambda:\mu]\in\P^1(\F_{p^2})\mid\lambda^{p+1}+\mu^{p+1}=0\}.$$
Note that $\J$ has $p+1$ elements.
\end{se}

\begin{lem}\label{g4}
There exists a bijection between the sets $\V(\Lambda)(k)$ with
$\Lambda\in\L_0^{(3)}$ which contain  $\M$ and the elements
$[\lambda:\mu]\in 
\J$. For $[\lambda:\mu]\in \J$  
the corresponding lattice $\Lambda\in\L_0^{(3)}$ is given by
\begin{align}
\Lambda=\langle
e_1,e_2,e_3,p^{-1}([\lambda]e_1+[\mu]e_3)\rangle_{W(\overline{\F}_p)}.
\label{eg2} 
\end{align}
In particular,  $\M$ is contained in $p+1$ sets $\V(\Lambda)(k)$.
\end{lem}

\begin{proof}
Let $\Lambda\in\L_0^{(3)}$ be a lattice with $\M_0^\tau\subset\Lambda$. 
We have to prove that $\Lambda$ is of the form (\ref{eg2}) for unique
$[\lambda:\mu]\in \J$. 
We have
\begin{align*}
p\Lambda\overset{1}{\subset}p(\M_0^\tau)\d\overset{1}{\subset}\M_0^\tau
\overset{1}{\subset}\Lambda=\Lambda\d\overset{1}{\subset}(\M_0^\tau)\d.
\end{align*}
Since $\M_0^\tau$ is of index 1 in $\Lambda$, there exist elements
$\lambda,\mu,\nu\in\F_{p^2}$ such that 
\begin{align*}
\Lambda=\langle e_1,e_2,e_3,p^{-1}([\lambda]e_1+[\nu]e_2+[\mu]e_3)
\rangle_{W(\F_{p^2})}. 
\end{align*}
As $\Lambda$ is totally isotropic, $\nu=0$ and $\lambda^{p+1}+\mu^{p+1}=0$.  
The element $[\lambda:\mu]\in\P^1(\F_{p^2})$ is uniquely determined by
$\Lambda$.

On the other hand,
it is clear that every lattice defined by (\ref{eg2}) is contained in
$\L_0^{(3)}$ and contains $\M_0^\tau$.
\end{proof}

\begin{se}\label{g5}
The Frobenius $\sigma$ acts on the set $\J$. 
We choose a set $\tilde{\J}=\{(\lambda_i,\mu_i)\}_{0\leq i\leq p}$ of
representatives of the different elements of
$\J$ such that $(\lambda_i^p,\mu_i^p)\in\tilde{\J}$ for every $i$.
We denote by $\sigma(i)$ the unique integer $j$ with $0\leq j\leq p$
and $(\lambda_i^p,\mu_i^p)=(\lambda_j,\mu_j)$.

We fix an integer $i$ with $0\leq i\leq p$.
Let $\Lambda_i$ be the 
lattice corresponding to $(\lambda_i,\mu_i)$ and let
$$e_\lmi=p^{-1}([\lambda_i]e_1+[\mu_i]e_3)\in\Lambda_i.$$ 
Then
$\{e_1,e_2,e_\lmi\}$ is a $W(\F_{p^2})$-basis of
$\Lambda_i$.   
Let $V$ be the $\overline{\F}_p$-vector space
$\Lambda_i/p\Lambda_i$ with 
induced basis $\{\oe_1,\oe_2,\oe_\lmi\}$. Denote by $(\cdot,\cdot)$ the
perfect form on $V$ induced by $\{\cdot,\cdot\}$. 
With respect to the above basis the form $(\cdot,\cdot)$ is given by the matrix
\begin{align*}
t\begin{pmatrix}
&&\lambda_i^p\\
&1&\\
\lambda_i&&\\
\end{pmatrix}\in\GL_3(\F_{p^2}).
\end{align*}

In \ref{s2} the variety $Y_{\Lambda_i}$ is defined as the
 closed subvariety of 
 $\Grass_2(V)$ over $\F_{p^2}$ given by
\begin{align*}
Y_{\Lambda_i}(R)=\{U\subset V_R\text{ a locally direct
  summand}\mid\rk_R 
U=2,\ U^\perp\subset U\}. 
\end{align*}
The superspecial lattice $\M$ corresponds via the bijection of
Proposition \ref{s20} to 
$\bU=\langle 
\oe_1,\oe_2\rangle$. Let $\calU_{\M,i}$ be the open neighbourhood
of $\bU$ in $Y_{\Lambda_i}$ given for each
$\F_{p^2}$-algebra $R$ 
by 
\begin{align*}
\calU_{\M,i}(R)=\{U_{a,b}=\langle
\oe_1+a\oe_\lmi,\oe_2+b\oe_\lmi\rangle\mid a,b\in
R\}\cap\V(\Lambda_i)(R). 
\end{align*}
\end{se}

\begin{lem}\label{g6}
We have
\begin{align*}
\calU_{\M,i}(R)=\{U_{a,b}=\langle
\oe_1+a\oe_\lmi,\oe_2+b\oe_\lmi\rangle\mid
a^p\lambda_i^p+a\lambda_i-b^{p+1}\lambda_i^{p+1}=0;\  a,b\in R\}. 
\end{align*}
In particular, $\calU_{\M,i}$ is isomorphic to $T_i=\Spec R_i$, where 
\begin{align}
R_i=\F_{p^2}[a,b]/(a^p\lambda_i^p+a\lambda_i-b^{p+1}\lambda_i^{p+1}).
\label{eg1}  
\end{align}
\end{lem}

\begin{proof}
For $U_{a,b}$ to be contained in $\V(\Lambda_i)(R)$, it is
necessary and 
sufficient that $U_{a,b}^\perp\subset U_{a,b}$. An easy computation
shows that $U_{a,b}\cap U_{a,b}^\perp\neq 0$ if and only if
$a^p\lambda_i^p+a\lambda_i-b^{p+1}\lambda_i^{p+1}=0$. In this case we obtain 
\[
U_{a,b}^\perp=\langle \oe_1-b^p\lambda_i^p\oe_2-a^p\lambda_i^{p-1}
\oe_\lmi\rangle. \qedhere
\]
\end{proof}

\begin{rem}\label{g25}
As the Frobenius twist $T_i^{(p)}$ of $T_i$ is isomorphic to $T_{\sigma(i)}$,
we will always identify $T_i^{(p)}$ with $T_{\sigma(i)}$. We denote by
$\Fr_{T_i}:T_i\rightarrow T_{\sigma(i)}$ the relative Frobenius.
\end{rem}

\begin{lem}\label{g2}
The variety $Y_{\Lambda_i}$ is isomorphic to the smooth
and projective 
curve $C_{\lambda_i}$ in $\P^2_{\F_{p^2}}$ given by the equation 
$$a^p\lambda_i^pd+a\lambda_i d^p-b^{p+1}\lambda_i^{p+1}=0$$
for $[a:b:d]\in\P^2_{\F_{p^2}}$.
\end{lem}

\begin{proof}
By Lemma \ref{g6} the projective nonsingular
curves $Y_{\Lambda_i}$ and $C_{\lambda_i}$ are
isomorphic over the 
open subvarieties
$\calU_{a,b}$ and $\{d\neq 0\}$ respectively. Therefore, they are
isomorphic.  
\end{proof}  

\begin{rem}\label{g24}
\begin{enumerate}
\item[a)]
The complement of $\calU_{\M,i}$ in  $Y_{\Lambda_i}$ contains only
one point. This point is
$\F_{p^2}$-rational. Indeed, this follows from the same result for
$C_{\lambda_i}$. 

\item[b)]
The curve $Y_{\Lambda_i}$ is isomorphic to the Fermat
curve $\C$ in $\P^2_{\F_{p^2}}$ given by the equation
$$x_0^{p+1}+x_1^{p+1}+x_2^{p+1}=0.$$
Indeed, as there exists up to isomorphism only one perfect skew-hermitian form
on $(\F_{p^2})^3$, the Fermat curve $\C$ is isomorphic to the curve
$C_{\lambda_i}$ of Lemma \ref{g2}.
\end{enumerate}
\end{rem}

\begin{se}\label{g7} 
For $a,b\in k$ denote by $c\in k$ the element $c=-\lambda_i\mu_i^{-1}a$.
Let 
$$f_\lmi=p^{-1}([\lambda_i^p]f_1+[\mu_i^p]f_3)\in N_1.$$
We define the following elements of the isocrystal $N_k$ which depend
on $a,b$.
\begin{align}\label{eg30}
\begin{aligned}
\e_1&=e_1-[b^{p^{-1}}\lambda_i^p]e_2+([a]-[b^{\frac{p+1}{p}}\lambda_i^p])
e_\lmi\\  
\e_2&=e_2+[b]e_\lmi\\
\e_3&=e_3-[b^{p^{-1}}\mu_i^p]e_2+([c]-[b^{\frac{p+1}{p}}\mu_i^p])e_\lmi\\
\f_1&=f_1-[b\lambda_i]p^{-1}f_2-[a\lambda_i^{1-p}]f_\lmi\\
\f_2&=f_2+[b^{p^{-1}}]pf_\lmi\\
\f_3&=f_3-[b\mu_i]p^{-1}f_2-[c\mu_i^{1-p}]f_\lmi.
\end{aligned}
\end{align}
Denote by $\delta_{\M,i}$ the open immersion 
\begin{align}
\delta_{\M,i}:T_i&\hookrightarrow Y_{\Lambda_i}\label{eg32}\\
(a,b)&\mapsto U_{a,b}\notag
\end{align}
of Lemma~\ref{g6}. 
By Proposition~\ref{s20} we have a bijection
$\theta_i:Y_{\Lambda_i}(k)\leftrightarrow\V(\Lambda_i)(k)$. 
Consider the following diagram
\begin{align}\label{eg31}
\begin{aligned}
\xymatrix{
&Y_{\Lambda_i}(k)\ar@{<->}[r]_-{\theta_{\Lambda_i}}&\V(\Lambda_i)(k)\\
T_i(k)\ar@/^35pt/[urr]^{\Psi_{\M,i}(k)}\ar@{^{(}->}[ur]_-{\delta_{\M,i}}\ar@{->}[r]_-\sim&\calU_{\M,i}(k)\ar@{^{(}->}[u]\ar@{<->}[r]^-{1:1}&\calS_{\M,i}(k)\ar@{^{(}->}[u]
}
\end{aligned}
\end{align}
By Remark~\ref{g24} a), the complement of the set
$\calS_{\M,i}(k)\subset\V(\Lambda_i)(k)$ defined in the above diagram
consists of only one superspecial point.  
\end{se}

\begin{prop}\label{g8}
The map $\Psi_{\M,i}(k)$ is given by
\begin{align*}
\Psi_{\M,i}(k):T_i(k)&\stackrel{1:1}{\longrightarrow}\calS_{\M,i}(k)\subset\V(\Lambda)(k)\\
(a,b)&\mapsto M_{a,b}=M_0\oplus M_1
\end{align*}
with
\begin{align}
M_0&=\langle \e_1,\e_2,\e_3\rangle_{W(k)}\label{eg17}\\
M_1&=\langle \f_1,\f_2,\f_3\rangle_{W(k)},\label{eg18}
\end{align}
and $V(M_{a,b})=\langle \e_1,p\e_2,\e_3,p\f_1,\f_2,p\f_3\rangle_{W(k)}$.
\end{prop}

\begin{proof}
For $(a,b)\in T_i(k)$ we have $\delta_{\M,i}(a,b)=U_{a,b}$.
Let
$\pi:\Lambda_i\rightarrow\Lambda_i/
p\Lambda_i$ 
be the natural 
projection and let $\pi_k$ be the base change of $\pi$ with $W(k)$. 
Then by Proposition~\ref{s20} we obtain $M_0=\pi_k^{-1}(U_{a,b})$ and
$M_1=F^{-1}(pM_0\d)$. We have   
\begin{align}
M_0&=\langle e_1+[a]e_\lmi, e_2+[b]e_\lmi, pe_\lmi\rangle\notag\\
&=\langle e_1+[a]e_\lmi, e_2+[b]e_\lmi,
e_3+[c]e_\lmi\rangle.\label{eg15} 
\end{align}
Similarly, we obtain for $pM_0\d=\pi_k^{-1}(U_{a,b}^\perp)$
\begin{align}
pM_0\d&=\langle
e_1-[b^p\lambda_i^p]e_2-[a^p\lambda_i^{p-1}]e_\lmi, pe_2,
pe_\lmi\rangle\notag\\ 
&=\langle e_1-[b^p\lambda_i^p]e_2-[a^p\lambda_i^{p-1}]e_\lmi,
pe_2, e_3-[b^p\mu_i^p]e_2-[c^p\mu_i^{p-1}]e_\lmi
\rangle\label{eg16} 
\end{align}
Note that for the equalities (\ref{eg15}) and (\ref{eg16}) we use
again $p\neq 2$. From (\ref{eg16}) we obtain 
\begin{align}
M_1=\langle f_1-[b\lambda_i]p^{-1}f_2-[a\lambda_i^{1-p}]f_\lmi,
f_2, f_3-[b\mu_i]p^{-1}f_2-[c\mu_i^{1-p}]f_\lmi\rangle. \label{eg26}
\end{align}
As $pf_\lmi=[\lambda_i^p]\f_1+[\mu_i^p]\f_3$, 
equality (\ref{eg18}) follows from (\ref{eg26}). The equality
(\ref{eg17}) follows from the equations 
\begin{align*}
\e_1&=(e_1+[a]e_\lmi)-[b^{p^{-1}}\lambda_i^p]\e_2\\
\e_3&=(e_3+[c]e_\lmi)-[b^{p^{-1}}\mu_i^p]\e_2.
\end{align*}
An easy calculation shows 
that $V(M_{a,b})$ has the desired
basis.
\end{proof} 

\begin{rem}
The map $\Psi_{\M,i}(k)$ is not a morphism. Indeed, the formulas
\eqref{eg30} show that the module $M_{a,b}$ is not defined over a
non perfect ring. 
\end{rem}

\begin{se}\label{g9}
Our goal is to
find an affine scheme $T_\M$ such that locally at $\M$ the variety
$\N_0^{red}$ is 
isomorphic to $T_\M$. 
Consider the polynomial ring in $a_i,b_i$ for $0\leq i\leq p$,
$$A=\F_{p^2}[a_i,b_i]_{0\leq i\leq p}.$$
For $0\leq i\leq p$ let $h_i\in A$ be the polynomial 
\begin{align*}
h_i=a_i^p\lambda_i^p+a_i\lambda_i-b_i^{p+1}\lambda_i^{p+1} 
\end{align*}
and let
$\mathfrak{a}\subset A$ be the ideal 
\begin{align*}
\mathfrak{a}=(h_i, a_ia_j, a_ib_j, b_ib_j)_{0\leq i\neq j\leq p}.
\end{align*}
Let $A'=A/\mathfrak{a}$ and let $Z$ be the affine scheme $\Spec A'$. 
Denote by $Z_i$ the closed subscheme $\Spec A'/(a_j,b_j)_{j\neq i}$ of $Z$.
\end{se}

\begin{lem}\label{g15}
The closed subschemes $Z_i$ are the irreducible components of
$Z$. They intersect 
transversally at the origin 0 and each $Z_i$ is isomorphic to 
$\calU_i$. Furthermore, $Z$ is reduced and its
tangent space at the origin has dimension $p+1$. 
\end{lem}

\begin{proof}
For $0\leq i\leq p $
the closed subscheme $Z_i$ is isomorphic to
$\calU_i$.  We obtain
$$Z\setminus \{0\}=\biguplus_{i=0}^p Z_i\setminus\{0\},$$ 
hence the $Z_i$ are the irreducible components of $Z$.  In particular, 
$Z$ is smooth away from the origin. 

We now prove that $A'$ is reduced. 
Since
$\calU_i$ is irreducible, it is clear that $h_i$ is
irreducible.
Let $f\in A$ such that $f^r\in\mathfrak{a}$ for an integer $r\geq
1$. 
We want to show that $f\in\mathfrak{a}$. 
By the definition of $\mathfrak{a}$,
we may 
assume that $f=\sum_{i=0}^p f_i$ with polynomials $f_i\in A$ which
depend only on $a_i,b_i$. Then  
\begin{align*}
f^r\equiv \sum_{i=0}^p f_i^r \mod \mathfrak{a}.
\end{align*}
As $f^r\in\mathfrak{a}$, we obtain that $f_i^r\in\mathfrak{a}$ for
every $i$ by the definition of $\mathfrak{a}$. Therefore, $f_i^r$ is
divisible by $h_i$. Since $h_i$ is 
irreducible, it divides $f_i$, which proves that $A'$ is reduced.

The tangent space at the origin is given by the equations $da_i=0$ for
$i=0,...,p$, where 
 $da_i$ and $db_i$    are the differentials of $a_i$ and $b_i$ respectively.
 Therefore, the 
tangent
space is of dimension $p+1$. In particular, all irreducible components
intersect transversally at the origin.
\end{proof}
 
\begin{se}\label{g20}
The moduli space $\calM$ of abelian varieties defined in the introduction 
is smooth of dimension 2. By the uniformization theorem of Rapoport
and Zink, the moduli space $\N$ is 
locally isomorphic to the closed subvariety $\calM^{ss}$ of
the special fibre of $\calM$ if $C^p$ is small enough
(comp.~Sec.~\ref{a}). Thus the 
tangent space of $\N^{red}$ at 
each closed 
point is at most of dimension 2.  Therefore, the local structure of
$\N_0^{red}$ at a 
supersingular point cannot be given by $Z$.
We will define a modification of the ring $A'$ such that the tangent
space at the origin has dimension 2 and prove that
the local structure of $\N_0^{red}$ is given by this modification.  

Consider
$R=\F_{p^2}[a_i,x,y]_{0\leq i\leq p}$. 
For an integer $k$ with $0\leq k\leq p$ let $g_k\in R$ be the polynomial
\begin{align*}
g_k=\sum_{i=0}^p(a_i^p\lambda_i^{p-k}\mu_i^k+a_i\lambda_i^{1-k}\mu_i^k)-
x^{p+1-k}y^k. 
\end{align*}
Define
\begin{align}
R_\M=R/(x^{p+1}+y^{p+1},a_ia_j, a_i(\lambda_iy-\mu_ix), g_k)_{
0\leq i\neq j\leq p,\ 
0\leq k\leq p}\label{eg3}
\end{align}
and denote by $T_\M=\Spec R_\M$ the corresponding affine scheme. Let
$$R_{\M,i}=R_\M/(a_j, \lambda_iy-\mu_ix)_{j\neq i}$$ 
and let $T_{\M,i}$ be the
corresponding closed subscheme of $T_\M$.

Let $R_i$ be as in (\ref{eg1}). We write $a_i,b_i$
instead of $a,b$ for the indeterminates of $R_i$.
We have the following equality in $R_i$
\begin{align*}
(\lambda_i^{-1}\mu_i)^kh_i=a_i^p\lambda_i^{p-k}\mu_i^k+a_i\lambda_i^{1-k}
\mu_i^k-  
(b_i\lambda_i)^{p+1-k}(b_i\mu_i)^k.
\end{align*}
Therefore, the morphism
\begin{align}
\eta_i: R_i&\rightarrow R_{\M,i}\label{eg28}\\
a_i&\mapsto a_i\notag\\
b_i&\mapsto \lambda_i^{-1}x\notag
\end{align}
is an isomorphism.
\end{se}

\begin{prop}\label{g10}
The closed subschemes $T_{\M,i}$ are the $p+1$ irreducible components of
$T_\M$ which intersect pairwise transversally at the
origin 0. Furthermore, 
$T_{\M,i}$ is isomorphic to 
$\calU_i$ for $0\leq i\leq p$, hence is smooth of
dimension 1.

In particular, $T_\M$ is of dimension $1$ and smooth away from the origin.
The tangent space of $T_\M$ at the origin is
$2$-dimensional. 
\end{prop}

\begin{proof}
As $R_{\M,i}$ is isomorphic to $R_i$ by \eqref{eg28},
the subscheme $T_{\M,i}$ is isomorphic to  $\calU_i$. Thus it is smooth of
dimension 1. The
schemes $T_\M$ and $T_{\M,i}$ coincide on the open locus $\{a_i\neq 0\}$.   
We have 
$$T_\M\setminus\{0\}=\biguplus_{i=0}^pT_{\M,i}\setminus\{0\}.$$
This shows that $T_\M$ is of dimension 1 and
smooth away from the origin.  

We now compute the tangent space of $T_\M$ at the origin. By
(\ref{eg3}) the tangent space at the origin in terms of the
differentials $da_0,...,da_p,dx,dy$ is given by the equations 
\begin{align*}
\sum_{i=0}^p \lambda_i^{1-k}\mu_i^kda_i=0
\end{align*}
for $0\leq k\leq p$. This system of equations has rank $p+1$
(Vandermonde determinant), hence we obtain $da_i=0$ for all $i$. This
proves that the tangent space has 
dimension 2. Furthermore, we see that 
any two irreducible components intersect transversally at the
origin. 
\end{proof} 

\begin{prop}\label{g11}
The scheme $T_\M$ is reduced.
\end{prop}

\begin{se}
To prove Proposition \ref{g11}, 
we first consider the homomorphism 
\begin{align*}
\psi:R=\F_{p^2}[a_i,x,y]_{0\leq i\leq p}&\rightarrow A=\F_{p^2}[a_i,b_i]_{0\leq i\leq p}\\
a_i&\mapsto a_i\\
x&\mapsto \sum_{i=0}^p b_i\lambda_i\\
y&\mapsto \sum_{i=0}^p b_i\mu_i.
\end{align*}
It is easy to see that
the homomorphism $\psi$ induces a homomorphism
$$\psi': R_\M\rightarrow A'.$$
We will  show that
$\psi'$ is injective. This will prove
that $R_\M$ is reduced, as $A'$ is reduced by Lemma \ref{g15}.

Let $\tilde{A}$ be the ring
$$\tilde{A}=A/(a_ia_j,  b_ib_j, a_ib_j)_{0\leq i\neq j\leq p}$$
and denote by $\tilde{R}$ the ring
$$\tilde{R}=R/(x^{p+1}+y^{p+1},a_ia_j, a_i(\lambda_iy-\mu_ix))_{0\leq
  i\neq j\leq p}.$$ 
We obtain a commutative diagram      
\begin{align*}
\xymatrix{
R\ar@/_1pc/[dd]_-{\alpha_3}\ar@{->}[r]^-{\psi}\ar@{->>}^-{\alpha_1}[d]&A\ar@{->>}_-{\beta_1}[d]\ar@/^1pc/[dd]^-{\beta_3}\\ 
\tilde{R}\ar@{->}[r]^-{\tilde{\psi}}\ar@{->>}^-{\alpha_2}[d]&\tilde{A}\ar@{->>}[d]_-{\beta_2}\\
R_\M\ar@{->}[r]^-{\psi'}&A'.\ 
}
\end{align*}
Note that all rings in this diagram are graded rings with respect to
the grading induced by the indeterminates $a_i$, $b_i$, $x$ and $y$.
All
homomorphisms respect the gradings.
\end{se}

\begin{lem}\label{g16}
The morphism $\tilde{\psi}:\tilde{R}\rightarrow \tilde{A}$ is
injective. 
\end{lem}

\begin{proof}
Let $f\in R$ such that
$\beta_1\circ\psi(f)=0$. We want to show that $\alpha_1(f)=0$. We may
assume that $f$ is of the form 
\begin{align}
f=\sum_{i=0}^p a_if_i(a_i,x,y)+\tilde{f}(x,y).\label{eg5}
\end{align}  
Here $f_i,\tilde{f}$ are polynomials in $R$ such that $f_i$ depends only on
$a_i,x,y$ and $\tilde{f}$ depends only on $x,y$. Since
$\beta_1\circ\psi(f)=0$, we obtain for all $i$ 
\begin{align}
\tilde{f}(b_i\lambda_i,b_i\mu_i)&=0\in A\\
a_if_i(a_i,b_i\lambda_i,b_i\mu_i)&=0\in A.\label{eg4}
\end{align}
Let $\tilde{f}_m(x,y)$ be the homogeneous component of degree $m$ of
$\tilde{f}$. Then 
\begin{align*}
\f_m(\lambda_ib_i,\mu_ib_i)=b_i^m\tilde{f}_m(\lambda_i,\mu_i)=0
\end{align*}
for all $i$.
Therefore, $\f_m(x,y)=0$ for all $(x,y)\in \Proj
\F_{p^2}[x,y]/(x^{p+1}+y^{p+1})$. Since this
scheme is reduced, 
we obtain that $x^{p+1}+y^{p+1}$ divides $\tilde{f}_m$. Thus
$\alpha_1(\tilde{f})=0$.  

Now consider the equation (\ref{eg4}). We may write  
\begin{align*}
f_i(a_i,x,y)=\sum_{r,m}a_i^rf_{r,m}(x,y),
\end{align*}
where $f_{r,m}$ is a polynomial in $x,y$ homogeneous of degree $m$. By
(\ref{eg4}) we obtain $f_{r,m}(\lambda_i,\mu_i)=0$. An analogous
argument as above shows that the polynomial $\mu_ix-\lambda_iy$
divides $f_{r,m}$, hence $a_i(\mu_ix-\lambda_iy)$ divides
$f_i$. Therefore, $\alpha_1(f_i)=0$ and $\tilde{\psi}$ is
injective.
\end{proof}

\begin{se}
{\it Proof of Proposition \ref{g11}.}\\
As $A'$ is reduced by Lemma \ref{g15}, it suffices to prove
that $\psi'$ is injective. 

Let $f\in R$ such that the
image of  $\beta_3\circ\psi(f)=0$. We may assume that $f$ is of the
form (\ref{eg5}). We obtain 
\begin{align} 
\tilde{\psi}(\alpha_1(f))=\sum_{i=0}^p
(a_if_i(a_i,b_i\lambda_i,b_i\mu_i)+\tilde{f}(b_i\lambda_i,b_i\mu_i)).
\label{eg6}  
\end{align}
As $\beta_3\circ\psi(f)=0$, there exist polynomials $w_i(a_i,b_i)$ and
$z_i(b_i)$ in $A$ with $z_i(b_i)$ depending only on $b_i$ such that
\begin{align}
\tilde{\psi}(\alpha_1(f))=\sum_{i=0}^p
(a_iw_i(a_i,b_i)+z_i(b_i))h_i.\label{eg24} 
\end{align}
Since $h_i$ is a polynomial of degree $p+1$ in $b_i$, it is easy to see that
every monomial in $\tilde{f}$ is of degree greater or equal than $p+1$. 
We write 
\begin{align}
\tilde{f}=\sum_{k=0}^p x^{p+1-k}y^k\tilde{f}_k(x,y).\label{eg7}
\end{align}
Define
\begin{align*}
f'=\sum_{i=0}^p a_iw_i(a_i,\lambda_i^{-1}x)g_0-\sum_{k=0}^p\tilde{f}_k(x,y)g_k.
\end{align*}

{\it Claim:} We have $\tilde{\psi}(\alpha_1(f'))=\tilde{\psi}(\alpha_1(f))$.\\
An easy calculation shows that 
\begin{align}
\tilde{\psi}(\alpha_1(f))-\tilde{\psi}(\alpha_1(f'))=\sum_{i=0}^ph_i[z_i(b_i)+
\sum_{k=0}^p
\lambda_i^{-k}\mu_i^k\tilde{f}_k(b_i\lambda_i,b_i\mu_i)].\label{eg27} 
\end{align}
We obtain from (\ref{eg6}) to (\ref{eg7}) with $a_i=0$ for all $i$ that 
\begin{align*}
-z_i(b_i)=\sum_{k=0}^p\lambda_i^{-k}\mu_i^k\tilde{f}_k(b_i\lambda_i,b_i\mu_i)
\end{align*}
which proves the claim by (\ref{eg27}).

Since $\tilde{\psi}$ is injective by Lemma \ref{g16}, we obtain that
$\alpha_1(f')=\alpha_1(f)$. By construction $\alpha_3(f')=0$, hence
$\psi'$ is injective. 
\hfill{$\Box$}
\end{se}

\begin{lem}\label{g27}
The completion $\hat{\O}_{T_\M,0}$ of the local ring of $T_\M$ at the
origin $0$ is isomorphic
to 
$$ \F_{p^2}[[x,y]]/\prod_{i=0}^p(\lambda_iy-\mu_ix).$$
In particular, $T_\M$ is of complete intersection.
\end{lem}

\begin{proof}
We use the notations of \ref{g20}.
An elementary calculation shows that
$\prod_{i=0}^p(\lambda_iy-\mu_ix)$ is contained in the ideal
$(g_k)_{0\leq k\leq p}$ of $R$. Thus the morphism 
\begin{align*}
\varphi^\circ:\F_{p^2}[[x,y]]/\prod_{i=0}^p(\lambda_iy-\mu_ix)&\rightarrow \hat R_{\M,0}\\
x&\mapsto x\\
y&\mapsto y
\end{align*}
is well defined. The corresponding schemes are both of dimension 1 and 
have $(p+1)$
irreducible components. Let $\varphi$ be the morphism of schemes corresponding to $\varphi^\circ$. It is easy to see that the tangent map of
$\varphi$ is an 
isomorphism. Therefore, $\varphi^\circ$ is surjective. Furthermore,
$\varphi$ induced an isomorphism on irreducible components, hence
$\varphi^\circ$ is injective. 
\end{proof}

\begin{se}\label{g13}

From now on, we use the theory of displays as in \cite{Zi2}. 
First of all, we recall the definition of a display. 
Let
$R$ be a ring of characteristic $p$. We denote  by $I_R$ the
image of the Verschiebung $\tau'$ on $W(R)$.
A $3n$-display over $R$ is a tuple $\calP=(P,Q,F,V^{-1})$ such that $P$ is
a finitely generated projective $W(R)$-module, $Q$ is a submodule of
$P$, and $F:P\rightarrow P$ and $V^{-1}:Q\rightarrow P$ are
$\sigma$-linear maps such that the following conditions are satisfied 
(\cite{Zi2} Def.~1).
\begin{enumerate}
\item[a)] $I_R P\subset Q\subset P$ and $P/Q$ is a direct summand of
  the $W(R)$-module $P/I_RP$.
\item[b)] $V^{-1}:Q\rightarrow P$ is a $\sigma$-linear epimorphism,
  i.e., the map $W(R)\otimes_{F,W(R)}Q\rightarrow P$ with 
$w\otimes m\mapsto wm$ is surjective.
\item[c)] For $x\in P$ and $w\in W(R)$ we have $V^{-1}(\tau'(w)x)=wF(x)$. 
\end{enumerate}
A $3n$-display $\calP$ is called a display if $\calP$ satisfies the
nilpotent condition of \cite{Zi2}~Def.~13. 

Let $M$ be a Dieudonn\'e module over a perfect field $k$. Then
$\calP_M=(M,V(M), F, V^{-1})$ is a $3n$-display. It is a display if
and only if $V$ is topologically nilpotent
on $M$ for the $p$-adic topology, i.e., if the associated $p$-divisible
group has no \'etale part. 

There exists a functor from the category of
displays over $R$ to the category of formal $p$-divisible groups over
$R$ (\cite{Zi2}, Thm.~81). 
If $R$ is of finite type over a field of characteristic $p$, this
functor is an equivalence of categories (\cite{Zi2}
Thm.~103). For such a ring $R$, we will identify the elements of
$\N_0(R)$ with the corresponding displays with additional structure.

For a $3n$-display $\calP=(P,Q,F,V^{-1})$, denote by $N$ the base 
change $N= P\otimes\Q$. 
Then $N$ is a projective $W_\Q(R)$-module and $F$ induces a $\sigma$-linear
operator on $N$ which we will again denote by $F$. The pair $(N,F)$ is
called the isodisplay of $\calP$ (\cite{Zi2} Def.~61). 

Now assume that
$R$ is torsion free
as an abelian group and that $P$ is of rank $n$. 
Then the morphism $W(R)\rightarrow W_\Q(R)$ is injective, 
hence $P$ is a $W(R)$-submodule of $N$. 
Let $P'$ be a projective 
$W(R)$-submodule of $N$ of rank $n$ 
and $Q'$ be a $W(R)$-submodule of $P'$ such that $P/Q$ is a direct summand of
the $W(R)$-module $P/I_RP$. 
Denote by $V^{-1}$ the operator 
$p^{-1}F$ on $N$. If $P'$ is $F$-invariant and $V^{-1}(Q')\subset P'$, 
the tuple $\calP'=(P',Q',F,V^{-1})$ is a $3n$-display. If $\calP$ is a 
display, the $3n$-display $\calP'$ is a display. 
Let $\Qisg_R$ be the category of displays over $R$ up to isogeny,
i.e., the objects are displays and the morphisms are 
$\Hom(\calP,\calP')\otimes_\Z\Q$.
Then the functor $\Qisg_R\rightarrow(Isodisplays)_R$ is fully faithful 
(\cite{Zi2} Prop.~66).
\end{se}

\begin{se}
We again use the notation of \ref{g17}. Let $\M\in\N_0(k)$ be a superspecial 
Dieudonn\'e lattice. As $\M$
is $\tau$-invariant, we obtain $V^2=\sigma^{-2}\id_\M$, hence the 
corresponding 
3n-display is a display.
Let $R_\M$ be the ring defined in (\ref{eg3}). By base change we obtain
a display $\calP_\M=(P_\M,Q_\M,F,V^{-1})$ over $R_\M$ with 
$$P_\M=\M\otimes_{W(\F_{p^2})}W(R_\M),$$ 
\begin{align}
Q_\M=\langle e_1,f_2,e_3\rangle_{W(R_\M)}\oplus I_{R_\M}\langle
f_1,e_2,f_3 \rangle_{W(R_\M)} \subset P_\M\label{eg14}
\end{align} 
and induced $F$ and $V^{-1}$ (\cite{Zi2} Def.~20).
The display $\calP_\M$ is equipped with an action of $\O_E$, i.e., a
$\Z/2\Z$-grading 
\begin{align*}
P_\M&=P_{\M,0}\oplus P_{\M,1}\\
&=\langle e_1,e_2,e_3\rangle_{W(R_\M)}\oplus \langle
f_1,f_2,f_3\rangle_{W(R_\M)}
\end{align*}
and an induced grading  $Q_\M=Q_{\M,0}\oplus Q_{\M,1}$. The perfect
alternating form $\langle\cdot,\cdot\rangle$ on $\M$ induces a perfect 
alternating
form on the display $\calP_\M$, i.e., a perfect alternating form
$$\langle\cdot,\cdot\rangle:P_\M\times P_\M\rightarrow W(R_\M)$$
such that 
$\tau'(\langle V^{-1}(x_1),V^{-1}(x_2)\rangle)=\langle x_1,x_2\rangle$
for all $x_1,x_2$ in $Q_\M$ (\cite{Zi2} Def.~18).

Let $N_\M=P_\M\otimes\Q$ and denote by $(N_\M,F)$ the isodisplay of $\calP_\M$
(\ref{g13}). As $R_\M$ is reduced (Prop.~\ref{g11}),
the morphism $P_\M\hookrightarrow N_\M$ is injective and we obtain a
$\Z/2\Z$-grading and a perfect alternating form on $N_\M$. As in \ref{g13}, we
denote by $V^{-1}$ the
operator $p^{-1}F$ on $N_\M$. 
\end{se}

\begin{se}\label{g28}
We now construct a display $\calP$ over $R_\M$ together
with a quasi-isogeny $\rho:\calP\rightarrow\calP_\M$ such that $\calP$
is equipped with all the data of $\N_0$. Then by base change this display will
define a morphism $T_\M\rightarrow \N_0$ (\ref{g21}). 

Let $c_i=-\lambda_i\mu_i^{-1}a_i$. Consider  the following elements of
$N_\M$  analogously to
\ref{g7} 
\begin{align}
\e_1&=e_1-[x]e_2+([\sum_{i=0}^p a_i^p\lambda_i^p]-[x^{p+1}])p^{-1}e_1+([\sum_{i=0}^p a_i^p\mu_i^p]-[xy^p])p^{-1}e_3\label{eg8}\\
\e_2&=e_2+[x^p]p^{-1}e_1+[y^p]p^{-1}e_3\label{eg9}\\
\e_3&=e_3-[y]e_2+([\sum_{i=0}^p c_i^p\lambda_i^p]-[x^py])p^{-1}e_1+([\sum_{i=0}^p c_i^p\mu_i^p]-[y^{p+1}])p^{-1}e_3\label{eg10}\\
\f_1&=f_1-[x^p]p^{-1}f_2-[\sum_{i=0}^p a_i^p\lambda_i^p]p^{-1}f_1-[\sum_{i=0}^p a_i^p\lambda_i^{p-1}\mu_i]p^{-1}f_3\label{eg11}\\
\f_2&=f_2+[x]f_1+[y]f_3\label{eg12}\\
\f_3&=f_3-[y^p]p^{-1}f_2-[\sum_{i=0}^p c_i^p\lambda_i\mu_i^{p-1}]p^{-1}f_1-[\sum_{i=0}^p c_i^p\mu_i^p]p^{-1}f_3.\label{eg13}
\end{align}
If we identify $R_{\M,i}$ with $R_i$ via $\eta_i$ as in \eqref{eg28},
we obtain over $R_{\M,i}$ the elements of \ref{g7} up to a Frobenius twist. 

Let $P=P_0\oplus P_1\subset N_\M$ be the module given by
\begin{align}\label{eg20}
\begin{aligned}
P_0&=\langle \e_1,\e_2,\e_3\rangle_{W(R_\M)}\\
P_1&=\langle \f_1,\f_2,\f_3\rangle_{W(R_\M)}. 
\end{aligned}
\end{align}
Denote by $Q$ the submodule $Q=\langle \e_1,\f_2,\e_3\rangle_{W(R_\M)}\oplus
I_{R_\M}\langle \f_1,\e_2,\f_3\rangle_{W(R_\M)}$ of $P$.
We define a morphism
\begin{align*}
\rho':\calP &\rightarrow \calP_\M\\
\e_i&\mapsto p\e_i\\
\f_i&\mapsto p\f_i.
\end{align*}
\end{se}

\begin{prop}\label{g12}
The tuple
$\calP:=(P,Q,F,V^{-1})$ is a subdisplay of the isodisplay $N_\M$. It is
invariant under the action of $\O_E$ and satisfies the determinant
condition of signature $(1,2)$. 
The module $P$ is free of rank $6$ and the form $\langle\cdot,\cdot\rangle$
is given by the matrix
\begin{align}
t\begin{pmatrix}
&I_3\\
-I_3&\\
\end{pmatrix}\label{eg21}
\end{align}
with respect to the above basis of $P$. 

The morphism $\rho'$ is a quasi-isogeny of height $6$ of displays with
$\O_E$-action such that for all $x_1,x_2\in \calP$ we have
$\langle\rho'(x_1),\rho'(x_2)\rangle=p^2\langle x_1,x_2\rangle$. 
\end{prop}

\begin{proof}
The display $\calP$ is $\O_E$-invariant because the $\Z/2\Z$-grading
on $N_\M$ induces the $\Z/2\Z$-grading  $P=P_0\oplus P_1$ and a
$\Z/2\Z$-grading $Q=Q_0\oplus Q_1$.    
To show that $\calP$ is a subdisplay of $N_\M$, we have to prove that
$P$ is invariant under $F$ and that $V^{-1}Q\subset P$.
Consider
\begin{align*}
L&=\langle \e_1,\f_2,\e_3\rangle_{W(R_\M)}\\
T&=\langle \f_1,\e_2,\f_3\rangle_{W(R_\M)}.
\end{align*}
We will show that $P=L\oplus T$ is a normal decomposition of $P$. It
is sufficient to prove that $FT\subset P$ and $V^{-1}L\subset P$.

As an example we will check that  $V^{-1}(\e_1)\in
P$. The other inclusions can be proved by a similar calculation. 
By (\ref{eg8}) we have 
\begin{align*}
V^{-1}(\e_1)=f_1-[x^p]p^{-1}f_2+([\sum_{i=0}^p
a_i^{p^2}\lambda_i]-[x^{p(p+1)}])p^{-1}f_1+([\sum_{i=0}^p
a_i^{p^2}\mu_i]-[x^py^{p^2}])p^{-1}f_3. 
\end{align*}
By (\ref{eg11}) we obtain
\begin{align*}
V^{-1}(\e_1)=\f_1&+([\sum_{i=0}^p a_i^p\lambda_i^p]+[\sum_{i=0}^p
a_i^{p^2}\lambda_i]-[x^{p(p+1)}])p^{-1}f_1\\ 
&+([\sum_{i=0}^p a_i^p\lambda_i^{p-1}\mu_i]+[\sum_{i=0}^p
a_i^{p^2}\mu_i]-[x^py^{p^2}])p^{-1}f_3. 
\end{align*}
Now an easy calculation using the relations of $R_\M$ (\ref{eg3}) shows that
\begin{align*}
V^{-1}(\e_1)=\f_1&+([\sum_{i=0}^p a_i^p\lambda_i^p]+[\sum_{i=0}^p
a_i^{p^2}\lambda_i]-[x^{p(p+1)}])p^{-1}\f_1\\ 
&+([\sum_{i=0}^p a_i^p\lambda_i^{p-1}\mu_i]+[\sum_{i=0}^p
a_i^{p^2}\mu_i]-[x^py^{p^2}])p^{-1}\f_3. 
\end{align*}
By definition of $R_\M$ (\ref{eg3}) we know that 
\begin{align*}
[\sum_{i=0}^p a_i\lambda_i]+[\sum_{i=0}^p
a_i^p\lambda_i^p]-[x^{p+1}]&\in I_{R_\M}\\ 
[\sum_{i=0}^p a_i\lambda_i^{1-p}\mu_i^p]+[\sum_{i=0}^p
a_i^p\mu_i^p]-[xy^p]&\in I_{R_\M},  
\end{align*}
hence their images under $\sigma$ are elements of $pW(R_\M)$. Thus
$V^{-1}(\e_1)\in P$. 
 
The following statements follow by a straightforward calculation.
The matrix of $\langle\cdot,\cdot\rangle$ on $P$ with respect to the basis
in (\ref{eg20}) is given by the matrix
$$
t\begin{pmatrix}
&I_3\\
-I_3&\\
\end{pmatrix},
$$
hence the form is perfect on $P$. 
As $\det (\e_1,\e_2,\e_3,\f_1,\f_2,\f_3)=1$, the module $P$ is free of rank 6.
We have $\rho'(Q)\subset Q_\M$, hence $\rho'$ is a quasi-isogeny. The
height of $\rho'$ is equal to 6. 

To prove the determinant condition of signature $(1,2)$, note that
$Q/I_{R_\M}P$ is 
isomorphic to the 
dual of the Lie algebra of the corresponding $p$-divisible group. By
construction the dimension of $Q_0/I_{R_\M}P_0$ is equal to 2 and the
dimension of $Q_1/I_{R_\M}P_1$ is equal to 1,
hence the determinant condition is satisfied.
\end{proof}

\begin{se}\label{g21}
Let $\calP_\bM$ be the display over $R_\bM$  of the $p$-divisible group $\X$ 
of \ref{d1}. Denote by $\rho_\M:\calP_\M\rightarrow\calP_\bM$ the
quasi-isogeny of height 0 induced by the two lattices $\bM$ and $\M$
in $N$.  
Let $\rho=p^{-1}\rho_\M\circ\rho':\calP\rightarrow\calP_\bM$. Then by
Proposition \ref{g12}, the 
morphism $\rho$ is a 
quasi-isogeny of height 0 of polarized displays. 
The pair $(\calP,\rho)$ defines by base change a morphism
\begin{align}
\Phi_\M:T_\M\rightarrow \N_0^{red}.\label{eg34}
\end{align}
Here we denote again by $T_\M$ the corresponding scheme over $\bF_p$. 
By construction of the display $\calP$ of Proposition~\ref{g12}, the
morphism $\Phi_\M$ depends on the choice of the 
basis $e_1,...,f_3$ of $\M$. 
We denote by $\Phi_{\M,i}$ the restriction of $\Phi_\M$ to the
irreducible component 
$T_{\M,i}$ of $T_\M$. 

Let $\xi$ be a $k$-rational point of $T_{\M,i}$ and let
$(\underline{a},x,y)\in 
k^{p+3}$ be the coordinates of $\xi$ with
$\underline{a}=(a_0,...,a_p)$. We have $a_j=0$ for $j\neq i$ and 
$y=\lambda_i^{-1}\mu_ix$.
We identify $T_{\M,i}$ with $T_i$ via the
isomorphism $\eta_i$ of \eqref{eg28}, i.e., the element $\xi$ corresponds to
$(a_i,b_i)$ with $b_i=x\lambda_i^{-1}$.
Consider the set 
$\calS_{\M,\sigma(i)}(k)\subset\V(\Lambda_{\sigma(i)})(k)$
as in \ref{g7} and consider the relative Frobenius
$\Fr_{T_{\M,i}}:T_{\M,i}\rightarrow 
T_{\M,\sigma(i)}$ (Rem.~\ref{g25}).
\end{se}

\begin{lem}\label{g23}
On $k$-rational points
$\Phi_{\M,i}$ induces a bijection between $T_{\M,i}$ and
$\calS_{\M,\sigma(i)}(k)$,
\begin{align}
\Phi_{\M,i}:T_{\M,i}(k)&\rightarrow \calS_{\M,\sigma(i)}(k)\label{eg29}\\
(a_i,b_i)&\mapsto M_{a_i^p,b_i^p},\notag
\end{align}
i.e., we have on $k$-rational points
\begin{align}
\Phi_{\M,i}=\Psi_{\M,\sigma(i)}(k)\circ \Fr_{T_{\M,i}}.\label{eg33}
\end{align}

In particular, the morphism $\Phi_{\M,i}$ is universally injective.
\end{lem}

\begin{proof}
Denote by $(\calP_\xi,\rho_\xi)$
the image of $\xi$ under $\Phi_\M$.
If
$\xi$ is equal to zero, i.e., $\underline{a}=\underline{0}$ and
$x=y=0$, then $(\calP_\xi,\rho_\xi)$ is 
equal to the superspecial
point $\M$. 

If we compare the formulas for $P_\xi$ in \ref{g28} with the
formulas of
$M_{a_i^p,b_i^p}\in\calS_{\M,\sigma(i)}(k)$ in
\ref{g7}, we find that 
$P_\xi=M_{a_i^p,b_i^p}$ as sublattices of the isocrystal
$N_k$. Therefore, $\Phi_{\M,i}$ is given by \eqref{eg29} and is equal to
the composition $\Psi_{\M,\sigma(i)}\circ \Fr_{T_{\M,i}}$ (Prop.~\ref{g8}). 
As $\Fr_{T_{\M,i}}$ and $\Psi_{\M,\sigma(i)}$ are universally
bijective on $k$-rational points, the morphism \eqref{eg29} is
universally bijective.
\end{proof}

\begin{prop}\label{g22}
The morphism $\Phi_\M$ is universally injective and the tangent morphism
at each closed point is injective. 
\end{prop}

\begin{proof}
By Lemma~\ref{g23} the morphisms $\Phi_{\M,i}$ are universally
injective for each $i$. 
The Frobenius induces a bijection on the set $\tilde{\J}$ of \ref{g5}, hence
$\Lambda_{\sigma(i)}\neq\Lambda_{\sigma(j)}$ for
$i\neq j$. By Lemma~\ref{g4} and \ref{g26}, we obtain
$\V(\Lambda_{\sigma(i)})(k)\cap\V(\Lambda_{\sigma(j)})(k)
=\{\M\}$. 
Therefore, the images of two irreducible components of $T_\M$
intersect only at the superspecial point $\M$. This proves that $\Phi_\M$ is
universally injective.

Now we want to show that $\Psi$ is injective on tangent spaces. 
Let $\xi$ be a $k$-rational point of
$T_\M$ and let $(\calP_\xi,\rho_\xi)$ be its image under $\Phi_\M$.
We know that
$(\calP_\xi,\rho_\xi)$ has 
no nontrivial automorphisms. 
Let $\alpha,\beta$ be two elements of the tangent
space of $T_\M$ at $\xi$. We denote by $(\calP_\alpha,\rho_\alpha)$
and $(\calP_\beta,\rho_\beta)$ the images of $\alpha$ and $\beta$
under $\Phi_\M$
respectively. We assume that 
$\Phi_\M(\alpha)=\Phi_\M(\beta)$. Then there exists an $\O_E$-linear
isomorphism 
$\varphi:\calP_\alpha\stackrel{\sim}{\rightarrow}\calP_\beta$ of
polarized displays such 
that $\rho_\alpha=\varphi\circ\rho_\beta$. We have to show that
$\alpha=\beta$.  

By (\ref{eg20}) the modules $P_\alpha$ and $P_\beta$ are graded
$W(k[\epsilon])$-modules. 
As $\varphi$ is $\O_E$-invariant, it preserves the grading.  
We denote by
$\e_{\alpha,i},\f_{\alpha,i}$ the basis of $\calP_\alpha$ as in
\eqref{eg20} and similarly for $\calP_\beta$.
With respect to
these bases, $\varphi$ is given by a matrix
\begin{align}
\begin{pmatrix}
A&\\
&B\\
\end{pmatrix}\in\GL_6(W(k[\epsilon])).\label{eg25}
\end{align}
The polarizations on $P_\alpha$ and
$P_\beta$ with respect to the above bases are given by the matrix (\ref{eg21}).
As $\varphi$ respects the polarizations, we obtain
$^t\!AB^\sigma=\, ^t\! BA^\sigma=I_3$. Thus $A=A^{\sigma^2}$ and
$B=B^{\sigma^2}$. In particular, the matrices $A$ and $B$ do not depend on
$\epsilon$, i.e., $A$ and $B$ are 
elements of $\GL_3(W(k))$. As
$\calP_\alpha\otimes_{W(k[\epsilon])}W(k)\cong \calP_\xi$ and
$\calP_\beta\otimes_{W(k[\epsilon])}W(k)\cong \calP_\xi$, the base change of
the isomorphism $\varphi$ induces an isomorphism $\overline{\varphi}$
of $\calP_\xi$. Since $A$ and $B$ do not depend on $\epsilon$, the
morphism $\overline{\varphi}$ is given by the matrix
(\ref{eg25}). Since $(\calP_\xi,\rho_\xi)$ has no 
nontrivial automorphisms,
this shows that $A=B=I_3$.

Now let $(\underline{a},x,y)\in k^{p+3}$  be the coordinates of $\xi$. The
computation of the tangent space of $T_\M$ (Prop.~\ref{g10}) shows
that $da_i=0$ for $i=0,...,p$. Thus the coordinates 
of $\alpha$ and $\beta$ are given by $(\underline{a},x_\alpha,y_\alpha)$ 
 and $(\underline{a},x_\beta,y_\beta)$ in
$(k[\epsilon])^{p+3}$ with 
\begin{alignat*}{2}
x_\alpha&\equiv x_\beta&&\equiv x\mod (\epsilon),\\
y_\alpha&\equiv y_\beta&&\equiv y\mod (\epsilon).
\end{alignat*} 
Since $\varphi$ is equal to the identity, we obtain that
$\f_{\alpha,2}=\f_{\beta,2}$, hence $x_\alpha=x_\beta$
and  $y_\alpha=y_\beta$. Therefore, $\alpha=\beta$
which proves the 
claim. 
\end{proof}



\section{The  global structure of $\N^{red}$ for $\GU(1,2)$}
\label{t}

\begin{se}\label{t1}
We use the notation of Section \ref{g}.
In this section we will construct a scheme $\T$ by gluing together
open subsets of the
varieties $T_\M$ for every 
superspecial point $\M\in\L_0^{(1)}$. We will 
prove that $\T$ is isomorphic to the supersingular  locus
$\N_0^{red}$. 

As the intersection behaviour of the sets $\V(\Lambda)(k)$ is given
by the simplicial complex $\B_0$ of Section \ref{i}, we will inductively
glue together the
varieties $T_\M$  using the complex $\B_0$.  
For $\Lambda\in\L_0^{(l)}$ we will for simplicity write $\Lambda$ instead of
$\{\Lambda\}$  for the corresponding vertex of $\Lambda$ in $\B_0$.

As in \ref{d4}, let $C$ be the $\Q_{p^2}$-vector space $N_0^\tau$ with
perfect skew-hermitian form $\{\cdot,\cdot\}$. Denote by $G$ the unitary
group of $(C,\{\cdot,\cdot\})$ with respect to the extension $\Q_{p^2}/\Q_p$. 
\end{se}

\begin{lem}\label{t3}
Let $\M\neq\M'$ be in $\L_0^{(1)}$ such that the corresponding vertices
 in $\B_0$ have a common neighbour $\Lambda$
(\ref{i10}).
Then there exists an element $g\in G(\Q_p)$ such that
$g\M=\M'$ and such that
$\Lambda$ is invariant under $g$. 
\end{lem}

\begin{proof}
Let $\G$ be the unitary group of $(\Lambda,\{\cdot,\cdot\})$
with respect to $\Z_{p^2}/\Z_p$. Then $G$ is isomorphic to the generic
fibre of $\G$.  Let
$V=\Lambda/p\Lambda$ and let $\pi:\Lambda\rightarrow V$ be the natural
projection. Then 
$\G_{\F_p}$ is equal to the unitary group of $(V,(\cdot,\cdot))$. Since $\G$ is
smooth over $\Spec\Z_p$, 
the canonical map 
$\varphi:\G(\Z_p)\rightarrow\G(\F_p)$ is surjective. 
As $\Lambda$ is a common neighbour of $\M$ and $\M'$,
the lattices $\M$ and $\M'$ are contained in
$\V(\Lambda)(\F_{p^2})$ (Prop.~\ref{i4}). They correspond to subspaces
$U,U'$ of $V$ of dimension 2  satisfying $U^\perp\subset U$ and
$(U')^\perp\subset U'$ respectively. As $(\cdot,\cdot)$ is a nondegenerate
skew-hermitian form on $V$, there exists an element
$\overline{g}\in\G(\F_p)$ with $\overline{g}(U)=U'$. Let
$g$ be a lift of $\overline{g}$ in $\G(\Z_p)$. As $\M=\pi^{-1}(U)$
and $\M'=\pi^{-1}(U')$, the automorphism $g$ of $\Lambda$ satisfies
the claim.
\end{proof}

\begin{se}\label{t2}
Let $\M,\M',\Lambda$ and $g$ be as in Lemma~\ref{t3}.
We fix a basis $e_1,...,f_3$ of $\M$ as in
Lemma~\ref{g3}. Then $\Lambda$ is equal to $\Lambda_i$ for some $i$ 
with $0\leq i\leq p$.
We define 
\begin{align}\label{et1}
\begin{array}{ll}
e'_j=ge_j&\ \ \ \text{for }j=1,2,3,\\
f'_j=F^{-1}(ge_j)&\ \ \   \text{for }j=1,3,\\  
f'_2=F^{-1}(pg(e_2)).& 
\end{array}
\end{align}
Let $e'_\lmi=p^{-1}([\lambda_i]e'_1+[\mu_i]e'_3)$ and let 
$\Lambda'_i=\langle
e'_1,e'_2,e'_3,e'_\lmi\rangle$.
\end{se}

\begin{lem}\label{t4}
The elements
$e'_1,...,f'_3$ form a basis of $\M'$ which 
satisfies the conditions of Lemma \ref{g3}. 
Furthermore,
$\Lambda$ is
equal to $\Lambda'_i$.
\end{lem}

\begin{proof}
This follows from Lemma \ref{t3}.
\end{proof}

\begin{se}\label{t16}
For  $\M\in\L_0^{(1)}$ consider the affine variety $T_\M$ and
its closed subvarieties $T_{\M,i}$ as in \ref{g20}. Geometrically, we
will glue
to each
$\F_{p^2}$-rational point  of $T_\M$  a variety $T_{\M'}$
such that $T_{\M,i}$ and $T_{\M',i}$ coincide on the open subsets of
non-$\F_{p^2}$-rational points. We will prove that the scheme $\T$
obtained by iterating this process is isomorphic to $\N_0^{red}$.

Let
$$\T_\M=T_\M\setminus\{\F_{p^2}\text{-rational points }\neq 0\}.$$ 
It is
an open subvariety of $T_\M$. Denote by $\T_{\M,i}$ the intersection of
$T_{\M,i}$ with $\T_\M$. Then $\T_{\M,i}$ is a closed subvariety of
$\T_\M$. Furthermore, let  
$$\T_{\M,i}^\circ=\T_{\M,i}\setminus\{0\}=T_{\M,i}\setminus
\{\F_{p^2}\text{-rational points}\}.$$   
We will always identify $T_{\M,i}$ with $T_i$ via \eqref{eg28}. Then
the variety 
$\T_{\M,i}^\circ$ 
is equal to $T_i^\circ=\Spec R_i^\circ$ with
\begin{align}
R_i^\circ=\F_{p^2}[a_i,b_i,(a_i^{p^2-1}-1)^{-1},(b_i^{p^2-1}-1)^{-1}]/
(a_i^p\lambda_i^p+a_i\lambda_i-b_i^{p+1}\lambda_i^{p+1}),\label{et2}
\end{align}
where we denote by $a_i,b_i$ instead of $a,b$ the indeterminates of $R_i$.
\end{se}

\begin{se}\label{t15}
We will now inductively glue together the varieties $\T_{\M}$ with
$\M\in\L_0^{(1)}$ along the open subsets $\T_{\M,i}^\circ$. We choose
a starting point $\hM\in\L_0^{(1)}$.  
For each
$\M\in\L_0^{(1)}$, we denote by $u_\M$ the distance of $\M$ to
$\hM$ (Rem.~\ref{i9}).
Let $\I$ be the set
of pairs
$(\M,\M')$ such that $\M,\M'\in\L_0^{(1)}$ and such that $\M$ and $\M'$ have a
common neighbour $\Lambda$. We assume that
$u_\M<u_{\M'}$.  

For the rest of this section we fix the following elements. 
For each $(\M,\M')\in\I$, 
we choose an element $g\in
G(\Q_p)$ such that $g(\M)=\M'$ and
$g(\Lambda)=\Lambda$, where $\Lambda$ is the common neighbour  of $\M$
and $\M'$.  
We choose a basis of $\hM$ as in Lemma~\ref{g3}. 
For each $\M\in\L_0^{(1)}$
we choose
inductively a basis $e_1,...,f_3$ of $\M$ 
which satisfies 
the conditions of Lemma~\ref{g3} such that for each $(\M,\M')\in\I$ 
the chosen basis of $\M'$ is given by \eqref{et1}. 
This is possible as $\B_0$ is a tree by Proposition \ref{i7}.

For $\M\in\L_0^{(1)}$
consider the morphism $\Phi_\M:T_\M\rightarrow \N_0^{red}$ of
\eqref{eg34} with respect to the fixed basis $e_1,...,f_3$ of $\M$.  
We denote by the same symbol $\Phi_\M$ its
restriction
to $\T_\M$. Let $\Phi_{\M,i}$ be the restriction of $\Phi_\M$
to $\T_{\M,i}$. 

Let $(\M,\M')\in\I$ and  let the common
vertex $\Lambda$ of $\M$ and $\M'$ be of the form
$\Lambda=\Lambda_i$ with respect to the chosen bases of $\M$ and $\M'$. 
\end{se}

\begin{lem}\label{t6}
There exists an isomorphism 
$\theta_g:T_{\M,i}^\circ\stackrel{\sim}{\longrightarrow}T_{\M',i}^\circ$
such that the restriction of the map 
$\Psi_{\M,i}(k)$ of \eqref{eg31}
to $T_{\M,i}^\circ(k)$ is equal to $\Psi_{\M',i}(k)\circ\theta_g$.
\end{lem}

\begin{proof}
Denote by
$Y_{\Lambda_i}^\circ$ the 
open subvariety of $Y_{\Lambda_i}$,
$$Y_{\Lambda_i}^\circ=Y_{\Lambda_i}\setminus\{\F_{p^2}\text{-rational
  points}\}.$$   
Then the morphism \eqref{eg32} induces by Remark~\ref{g24} a) an
isomorphism
$$\delta_{\M,i}:T_i^\circ\stackrel{\sim}{\longrightarrow}Y_{\Lambda_i}^\circ$$
and similarly for $\delta_{\M',i}$. Note that $\delta_{\M,i}$ and
$\delta_{\M',i}$ depend 
on the chosen bases of $\M$ and $\M'$ respectively.
We use the notation of the proof of Lemma~\ref{t3}.
The element $(\overline{g})^{-1}\in\G(\F_p)$ induces an
isomorphism on $Y_{\Lambda_i}$ which preserves the set of $\F_{p^2}$-rational
points. There exists an automorphism $\theta_g$ of $T_i^\circ$ such that
the following 
diagram of isomorphisms commutes
$$\xymatrix{
Y_{\Lambda_i}^\circ\ar@{->}[r]^-{(\overline{g})^{-1}}&Y_{\Lambda_i}^\circ\ar@{->}[r]^-{\overline{g}}&Y_{\Lambda_i}^\circ\\
T_i^\circ\ar@{->}[u]^-{\delta_{\M,i}}\ar@{->}[r]^-{\theta_g}&T_i^\circ\ar@{->}[u]^-{\delta_{\M,i}}\ar@{->}[ur]_-{\delta_{\M',i}}&.
}
$$
The above triangle commutes 
by definition of $g$ and the choice of the basis of $\M'$.
Thus $\delta_{\M,i}=\delta_{\M',i}\circ\theta_g$. 
As $T_{\M,i}^\circ$ and $T_{\M',i}^\circ$ are equal to $T_i^\circ$ by
\ref{t16}, we obtain an
isomorphism $\theta_g:T_{\M,i}^\circ\stackrel{\sim}{\longrightarrow}
T_{\M',i}^\circ$. 
The claim follows
from diagram \eqref{eg31}.
\end{proof} 
 
\newpage
\begin{prop}\label{t7}
Let
$\theta_g^{(p)}$ 
be the Frobenius pullback of the isomorphism $\theta_g$ of Lemma~\ref{t6}.
Then the diagram
\begin{align*}
\xymatrix{
T_{\M,\sigma(i)}^\circ\ar[dr]^-{\Phi_{\M,\sigma(i)}}\ar[dd]^-{\theta_g^{(p)}}_-\sim&\\
&\N_0^{red}\\
T_{\M',\sigma(i)}^\circ\ar[ur]_-{\Phi_{\M',\sigma(i)}}&
}
\end{align*}
commutes.
\end{prop}

\begin{proof}
It is sufficient to prove the claim on $k$-rational points.
Lemma~\ref{g23} shows that on $k$-rational points
$\Phi_{\M,\sigma(i)}$ is equal to $\Psi_{\M,i}(k)\circ
\Fr_{T_{\M,\sigma(i)}^\circ}$ 
Similarly, we obtain 
$\Phi_{\M',\sigma(i)}=\Psi_{\M',i}(k)\circ
\Fr_{T_{\M',\sigma(i)}^\circ}$. As
$\theta_g^{(p)}$ is defined over $\F_{p^2}$, we obtain 
\begin{align*}
\Fr_{T_{\M',\sigma(i)}^\circ}\circ
\theta_g^{(p)}=\theta_g\circ\Fr_{T_{\M,\sigma(i)}^\circ}.
\end{align*} 
Thus 
\begin{align*}
\Phi_{\M',i}\circ \theta_g^{(p)}=\Psi_{\M',i}(k)\circ
\theta_g\circ 
\Fr_{T_{\M,\sigma(i)}^\circ}
\end{align*}
and the claim follows from Lemma~\ref{t6}.
\end{proof}

\begin{prop}\label{t9}
There exists a reduced scheme $\T$ locally of finite type over
$\F_{p^2}$ of dimension 1 and a morphism 
\begin{align}
\Phi:\T\rightarrow \N_0^{red}
\end{align}
which satisfies the following
conditions.
\begin{enumerate}
\item[a)]
For each $\M\in\L_0^{(1)}$, the scheme $\T_\M$ can be identified with
an open subscheme of $\T$ such that the restriction of  $\Phi$ to
$\T_\M$ is equal to $\Phi_\M$. 
\item[b)] The open subschemes $\T_\M$ with $\M\in\L_0^{(1)}$ form an open
  covering of $\T$.
\item[c)] Let $\M,\M'\in\L_0^{(1)}$. The open subschemes $\T_\M$ and
  $\T_{\M'}$ of 
  $\T$ intersect if and only if $\M$ and $\M'$ have a common neighbour
  in $\B_0$. In this case, there exists an integer $i$ such that the
  intersection $\T_\M\cap\T_{\M'}$ is equal to the open subschemes
  $\T_{\M,i}^\circ$ and $\T_{\M',i}^\circ$ of $\T_\M$ and $\T_{\M'}$
  respectively.  
\end{enumerate}
\end{prop}

\begin{proof}
We glue the schemes $\T_\M$ together along the open subsets
$\T_{\M,i}^\circ$ by induction over the distance $u_\M$  of 
$\M$ from $\hM$. Let $\alpha$ be an even nonnegative integer.
Assume that the schemes $\T_\M$ with $\M$ of distance $u_\M\leq \alpha$
have been glued together to a scheme 
$\T_\alpha$ with morphism $\Phi_\alpha$ which satisfies  condition a)
and the 
condition corresponding to c) of 
the proposition. Let $\M'\in\L_0^{(1)}$ be of distance $\alpha+2$. 
As $\B_0$ is a tree, there exists exactly one element
$\M\in\L_0^{(1)}$ of distance $\alpha$ such that $\M$ and
$\M'$ have a common neighbour $\Lambda$ in $\B_0$. Let
$\Lambda=\Lambda_i$ with respect to the chosen bases of $\M$ and $\M'$.
We glue together $\T_\alpha$ and $\T_{\M'}$ along the open subschemes
$\T_{\M,\sigma(i)}^\circ$ and $\T_{\M',\sigma(i)}^\circ$ respectively via the
isomorphism $\theta_g^{(p)}$ of Proposition~\ref{t7}. 
The same 
proposition shows that  the morphisms  
$\Phi_\alpha$ and $\Phi_\M$ glue together to a morphism. 

As $\B_0$ is a tree and 
the glueing is defined inductively, no cocycle
condition has to be checked.
\end{proof}

\begin{se}
For $\Lambda\in\L_i^{(3)}$ let $\M\in\L_0^{(1)}$ be a neighbour of
$\Lambda$ in $\B_0$. Let $\Lambda=\Lambda_i$ with respect to the
chosen basis of 
$\M$. Denote by $\T_\Lambda$ 
the closure of the open subvariety $\T_{\M,i}$ in $\T$. 
\end{se}

\begin{prop}\label{t13}
The variety $\T_\Lambda$ only depends on $\Lambda$ and is isomorphic
to $Y_\Lambda$.
The scheme $\T$ is connected and the varieties $\T_\Lambda$ with
$\Lambda\in\L_0^{(3)}$ are its irreducible components. 
\end{prop}

\begin{proof}
As $\B_0$ is a tree, the scheme $\T$ is connected. The claim follows
from the construction of $\T$.
\end{proof}

\begin{se}
Denote by $\Phi_\Lambda$ the restriction of $\Phi$ to the irreducible component
$\T_\Lambda$. As $\T_\Lambda$ is isomorphic to $Y_\Lambda$ by
Proposition~\ref{t13}, it is projective.  
By construction and Lemma~\ref{g23}, the morphism $\Phi_\Lambda$
is universally injective and the
image of  
$\Phi_\Lambda$ on  $k$-rational points is equal to $\V(\Lambda)(k)$. 
The moduli space $\N_0^{red}$
is separated (\ref{e4}), hence $\Phi_\Lambda$ is finite.
We denote by
$\V(\Lambda)$ the scheme theoretic image of 
$\Phi_\Lambda$. To show that $\Phi_\Lambda$ induces an isomorphism
onto 
$\V(\Lambda)$, we need the following lemma.
\end{se}

\begin{lem}\label{t11}
Let $f:X\rightarrow X'$
be a morphism of
reduced schemes of finite type over an algebraically closed field
$k$. We assume that $f$ is finite, 
universally bijective and that the  tangent
morphism is injective at every 
closed point. Then $f$ is an isomorphism. 
\end{lem}

\begin{proof}
Obviously, the morphism $f$ is a homeomorphism.
We will
show that 
$$f^\#_x:\O_{X',f(x)}\rightarrow \O_{X,x}$$
is surjective for every
closed point $x$ of $X$. Then $f$ will be
a closed immersion and, 
as $X'$ is reduced, $f$ will be an
isomorphism. 

We may assume that $X=\Spec R$ and $X'=\Spec R'$ are affine. 
Denote
by $\mathfrak{m}$ a maximal ideal of $R$ and by
$\mathfrak{m}'=f(\mathfrak{m})$  its image in $R'$.  
The morphism
\begin{align}
f_{R'_{\mathfrak{m}'}}:\Spec(R\otimes_{R'}R'_{\mathfrak{m}'})\rightarrow\Spec
R'_{\mathfrak{m}'}  
\end{align}
is finite, hence $R\otimes_{R'}R'_{\mathfrak{m}'}$ is a semi-local
ring. As $f_{R'_{\mathfrak{m}'}}$ is universally bijective,
$R\otimes_{R'}R'_{\mathfrak{m}'}$ is a local ring and we obtain
$R_\mathfrak{m}=R\otimes_{R'}R'_{\mathfrak{m}'}$. Thus
$R_\mathfrak{m}$ is a finite 
$R'_{\mathfrak{m}'}$-module.

Furthermore, the tangent morphism at $\mathfrak{m}$ is injective,
hence the morphism
\begin{align*}
\mathfrak{m}'/(\mathfrak{m}')^2\twoheadrightarrow \mathfrak{m}/\mathfrak{m}^2
\end{align*}
is surjective. We obtain  a surjective morphism
$\hat{R}'_{\mathfrak{m}'}\twoheadrightarrow \hat{R}_\mathfrak{m}$.  
Since $R_\mathfrak{m}$ is a finite
$R'_{\mathfrak{m}'}$-module,
$\hat{R}_\mathfrak{m}=R_\mathfrak{m}\otimes_{R'_{\mathfrak{m}'}}
\hat{R}'_{\mathfrak{m}'}$.  
Therefore, the morphism
$f^\#_\mathfrak{m}:R'_{\mathfrak{m}'}\rightarrow R_\mathfrak{m}$ is
surjective. 
\end{proof} 

\begin{thm}\label{t10}
The morphism $\Phi$ is an isomorphism  which induces  an isomorphism of
$\T_\Lambda$ onto $\V(\Lambda)$ for every
$\Lambda\in\L_0^{(3)}$. 
\end{thm}

\begin{proof}
By Proposition~\ref{g22} and Proposition~\ref{t9}, the tangent
morphism of $\Phi_\Lambda$ is injective at 
every closed point. Therefore, the morphism $\Phi_\Lambda$ induces an
isomorphism of $\T_\Lambda$ onto $\V(\Lambda)$ 
by Lemma~\ref{t11}. 

Now we will prove that $\Phi$ is an isomorphism. By construction of $\Phi$, the
$\F_{p^2}$-rational points of $\T$ correspond to the superspecial
points of $\N_0^{red}$ and
the intersection behaviour of the varieties $\V(\Lambda)$ is equal
to the intersection behaviour of the irreducible components
$\T_\Lambda$ of $\T$. Thus $\Phi$ is
universally bijective and
the varieties $\V(\Lambda)$ are the irreducible components of
$\N_0^{red}$.
We obtain that $\Phi$ is an isomorphism locally at every point of $\T$ which is
not $\F_{p^2}$-rational.

Let $x$ be an $\F_{p^2}$-rational point of $\T$ and let
$\M=\Phi(x)$ be the  corresponding  superspecial point. The variety
$\T_\M$ is an open neighbourhood of $x$ in $\T$.
Denote by $\oT_\M$ the closure of $\T_\M$ in $\T$ and 
denote by
$Z_\M$ the image of $\oT_\M$ under $\Phi$, i.e., the union of the
varieties $\V(\Lambda)$ which contain 
$\M$. 
The induced
morphism
$$\Phi\vert_{\oT_\M}:\oT_\M\rightarrow Z_\M$$
is bijective and injective on tangent spaces by Proposition~\ref{g22}
and Proposition~\ref{t9}. The morphism $\Phi\vert_{\oT_\M}$ is
universally closed as its 
restriction to every irreducible 
component is finite, hence $\Phi\vert_{\oT_\M}$ is finite. 
By Lemma~\ref{t11} it is an isomorphism which proves the claim.
\end{proof}

\begin{se}\label{t14}
Let $i$ be an even integer and let $\Lambda\in\L_i^{(3)}$. The isomorphism 
$\Psi_i:\N_i\stackrel{\sim}{\longrightarrow}\N_0$ of \eqref{ed25} maps
$\V(\Lambda)(k)$  to 
a set $\V(\Lambda')(k)$ for a lattice $\Lambda'\in\L_0^{(3)}$
(Rem.~\ref{s18} d)). We denote by 
$\V(\Lambda)$ the preimage of the variety $\V(\Lambda')$ by $\Psi_i$. 
\end{se}

\begin{thm}\label{t12}
The  schemes $\N_i^{red}$ with $i\in\Z$ even are the connected components
of $\N^{red}$ which are all isomorphic to each other.
The varieties $\V(\Lambda)$
with $\Lambda\in\L_i^{(3)}$ 
are the irreducible components of $\N_i^{red}$. 

The singular points of $\N^{red}$
are  the superspecial 
points. 
Each $\V(\Lambda)$ contains $p^3+1$ superspecial points and each
superspecial point is the intersection 
of $p+1$ irreducible components $\V(\Lambda)$. 
Two irreducible components
intersect transversally in at most one superspecial point and the intersection 
graph of $\N_i^{red}$ is a tree.

Each variety $\V(\Lambda)$ is
isomorphic to the Fermat curve $\C$ (Rem.~\ref{g24}). The superspecial
points of 
$\V(\Lambda)$ correspond to 
the $\F_{p^2}$-rational points of $\C$. 

The scheme
$\N^{red}$ is equi-dimensional of dimension 1 and of complete intersection.
\end{thm}

\begin{proof}
We first consider the case $i=0$.
The 
incidence relation of the varieties $\V(\Lambda)$
follows from Proposition~\ref{i7}. 
The geometric statements follow from 
Theorem~\ref{t10}, Proposition~\ref{t13} and Lemma~\ref{g27}.

Now consider the general case.
By \ref{t14} the varieties 
$\V(\Lambda)$ with
$\Lambda\in\L_i^{(l)}$ correspond to the varieties $\V(\Lambda)$ with
$\Lambda\in\L_0^{(l)}$ under the isomorphism $\Psi_i$. The claim
follows from the case 
$i=0$.
\end{proof}


\section{The structure of the supersingular locus of $\calM$ for $\GU(1,2)$}
\label{a}

\begin{se}
Let $\calM$ be the moduli space of abelian varieties defined in the
introduction. 
In this section, we carry over our results on the moduli space $\N$ to
$\calM$ in
the case of $\GU(1,2)$. 

We use the notation of the introduction.
In particular, we now denote by $E$ an
imaginary quadratic extension of $\Q$ such that $p$ is inert in $E$
and denote by $E_p$ the completion of $E$ with respect to the $p$-adic
topology. 
Consider the supersingular locus $\calM^{ss}$ of the special fibre
$\calM_{\F_{p^2}}$ of 
$\calM$. It is a closed subscheme of $\calM_{\F_{p^2}}$ which contains
an $\bF_p$-rational point
(\cite{BW} Lem.~5.2). We view $\calM^{ss}$ as a scheme over
$\bF_p$. We say that a 
point of $\calM^{ss}$ is superspecial if the underlying abelian
variety is superspecial.
Let
$(A\otimes_\Z\Z_{(p)},\iota\otimes_\Z\Z_{(p)},\overline{\lambda},
\overline{\eta})$ 
be an $\bF_p$-valued point of $\calM^{ss}$ with corresponding
$p$-divisible group $(\X,\iota,\lambda)$. Let $\N$ be the moduli space
of quasi-isogenies of \ref{e4} with respect to $(\X,\iota,\lambda)$. 
\end{se}

\begin{se}
Let $J$ be the group of similitudes of the isocrystal $N$ of $\X$ with
additional structure as
in \ref{d24}.
We denote by $I(\Q)$ the group of quasi-isogenies in
$\End_{\O_E}(A)\otimes\Q$ which respect the homogeneous polarization
$\overline{\lambda}$. It is a subgroup of $J$. 
Using the level structure of $A$,
one can define an injective morphism of $I(\Q)$ into
$G(\A_f^p)$ (\cite{RZ} 6.15).
By \cite{RZ} Theorem~6.30, the set $I(\Q)\backslash G(\A_f^p)/C^p$ is finite.
Denote by $g_1,...,g_m\in G(\A_f^p)$ representatives  of the
different elements of $I(\Q)\backslash G(\A_f^p)/C^p$. For every
integer $j$ with $1\leq j\leq m$, let $\Gamma_j$
be the group 
$$\Gamma_j=I(\Q)\cap g_j^{-1}C^pg_j.$$
We view $\Gamma_j$ 
as a subgroup of $J$.
\end{se}

\begin{se}\label{a1}
We recall the uniformization theorem of Rapoport and Zink in case of
$\GU(1,2)$. We will formulate this theorem only 
for the underlying 
schemes, not for the formal schemes.

There exists an isomorphism of schemes over $\Spec\bF_p$
\begin{align}
I(\Q)\backslash(\N^{red}\times G(\A_f^p)/C^p)\stackrel{\sim}{\longrightarrow}
\calM^{ss}.\label{ea1}
\end{align}
The left hand side is isomorphic to the disjoint union of the quotients
$\Gamma_j\backslash\N^{red}$ for $1\leq j\leq m$.
Each group $\Gamma_j\subset J$ is
discrete and cocompact modulo center. If $C^p$ is small enough,
$\Gamma_j$ is torsion free. 

Indeed,
the proof of \cite{RZ} Theorem 6.30 shows that \eqref{ea1}
is an isomorphism if $G$ satisfies the Hasse principle, i.e., if the
kernel of the Hasse map with respect to $G$,
\begin{align}
\hH^1(\Q,G)\rightarrow \prod_{v\text{ place of }\Q}\hH^1(\Q_v,G),\label{ea6}
\end{align}
is trivial. By \cite{Ko} Section~7, the kernel of \eqref{ea6} is equal
to the kernel of the Hasse map with respect to the group
$\Res_{E/\Q}(\bG_m)$. But $\hH^1(\Q,\Res_{E/\Q}(\bG_m))$ is trivial,
thus $\Res_{E/\Q}(\bG_m)$, and
hence $G$, satisfies the Hasse principle.
 
The decomposition into the schemes $\Gamma_j\backslash\N^{red}$
follows from the proof of 
\cite{RZ} Theorem~6.23. The properties of the subgroups $\Gamma_j$
are proved as well.
\end{se}

\begin{se}
Let $J^0$ be the subgroup of $J$ as in \ref{d24} and 
denote by $\Gamma_j^0$ the intersection of
$\Gamma_j$ with $J^0$.
Consider the morphism
\begin{align}
\Psi:\coprod_{j=1}^m\N^{red}\rightarrow\calM^{ss}\label{ea5}
\end{align}
induced by \eqref{ea1}.
\end{se}

\begin{thm}\label{a2}
Let $C^p$ be small enough.
The morphism $\Psi$ is surjective and an isomorphism locally at each
point.
The restriction of $\Psi$ to each irreducible
component $\V(\Lambda)$ of $\N^{red}$ is a closed immersion and the
images of two
irreducible components of $\N^{red}$ in $\calM^{ss}$ intersect in at most one
superspecial point. 
\end{thm}

\begin{proof}
By \ref{a1} we have to show that
\begin{align}
\Psi_j:\N^{red}\rightarrow\Gamma_j\backslash\N^{red}\label{ea3}
\end{align}
is an isomorphism locally at each point for every integer $j$ with
$1\leq j\leq m$. 
We use the notation of \ref{d24}.
Let $g\in J$ and let $\alpha=v_p(c(g))$. 
By  \ref{d24}  the action of $g$
defines an isomorphism of
$\N_i$ with $\N_{i+\alpha}$ for every integer $i$. Note that $\alpha$
is even (Lem.~\ref{d23}) and that $\N_i$ is empty if $i$ is odd
(Lem.~\ref{d17}). 
The schemes $\N_i^{red}$
with $i$ even are the connected components of $\N^{red}$
(Thm.~\ref{t12}) which are all isomorphic to each other (Prop.~\ref{d21}).
Thus we obtain
\begin{align}
\Gamma_j\backslash\N^{red}\stackrel{\sim}{\longrightarrow}
\coprod_{(\Gamma_jJ^0)\backslash J}\Gamma_j^0\backslash\N_0^{red}.\label{ea4}
\end{align}
In particular, the index of $(\Gamma_jJ^0)$ in $J$ is finite as
$\calM^{ss}$ is of finite type over $\Spec\bF_p$.

We now want to understand the action of $\Gamma_j^0$ on
$\N_0^{red}$. 
As in \ref{d24},
we view $J$ as group of similitudes of $(C,\{\cdot,\cdot\})$.
For $l=1,3$ the group $\Gamma_j^0$ acts on the set $\L_0^{(l)}$
(Def.~\ref{s22}).  We will show that the action of $\Gamma_j^0$ on
$\L_0^{(l)}$ has no  
fixed points.

Indeed,
for $\Lambda\in\L_0^{(l)}$ denote by $\Stab(\Lambda)$ the stabilizer of
$\Lambda$ in $J^0$. This is a compact open subgroup of $J^0$, hence
$\Gamma_j^0\cap\Stab(\Lambda)$ is finite. If $C^p$ is small enough,
$\Gamma_j^0$ has no torsion. Thus the intersection of $\Gamma_j^0$ with
$\Stab(\Lambda)$ is trivial. 

We have proved that each element of $\Gamma_j^0$ maps every
irreducible component $\V(\Lambda)$ of $\N_0^{red}$ with
$\Lambda\in\L_0^{(3)}$ onto a different irreducible
component. Furthermore, it fixes no superspecial points. As two
irreducible components of $\N_0^{red}$ intersect in at most one
superspecial point, the 
action of $\Gamma_j^0$ on $\N_0^{red}$ has no fixed points. Thus the
morphism $\N_0^{red}\rightarrow\Gamma_j^0\backslash\N_0^{red}$ is an
isomorphism locally at every point. 

The action of $\Gamma_j^0$ on $\L_0^{(1)}$ and $\L_0^{(3)}$ induces an
action of 
$\Gamma_j^0$ on the simplicial complex $\B_0$ of Definition~\ref{i3}. As
$\B_0$ describes the incidence relation of the irreducible components
and the superspecial points of $\N_0^{red}$, we can choose $C^p$ small
enough such that for every $g\in\Gamma_j^0$ and every $\Lambda\in\L_0$
the distance $u(\Lambda,g\Lambda)$ (Rem.~\ref{i9}) is greater or
equal than 6. In this case, the restriction of the morphism
$\N_0^{red}\rightarrow\Gamma_j^0\backslash\N_0^{red}$ to every
irreducible component is a closed immersion and the images of two
irreducible components of $\N_0^{red}$ in
$\Gamma_j^0\backslash\N_0^{red}$ intersect in at most one point.
This proves the theorem.
\end{proof}

\begin{cor}\label{a4}
Let $C^p$ be small enough.
The supersingular locus $\calM^{ss}$ is locally isomorphic to
$\N^{red}$. It is equi-dimensional of dimension 1 and of complete
intersection. 

The 
singular points of $\calM^{ss}$ are the superspecial points. 
Each superspecial point is the pairwise transversal intersection of $p+1$
irreducible components. 
The irreducible
components are isomorphic to 
the Fermat curve $\C$ (Rem.~\ref{g24}) and contain
$p^3+1$ superspecial points.
Two irreducible components
intersect in at most one superspecial point. 
\end{cor}

\begin{proof}
The claim  follows from Theorem~\ref{a2} and Theorem~\ref{t12}.
\end{proof}

\begin{prop}\label{a3}
Let $C_{J,p}$ and $C_{J,p}'$ be maximal compact  subgroups of $J$ such
that $C_{J,p}$ 
is hyperspecial and $C_{J,p}'$ is not hyperspecial, i.e., $C_{J,p}$ is the
stabilizer in $J$ of a lattice
$\Lambda\in\L_0^{(3)}$ and $C_{J,p}'$ is the stabilizer in $J$ of a
lattice  $M\in\L_0^{(1)}$.  
If $C^{J,p}$ is small enough, we have
\begin{align*}
\#\{\text{irreducible components of }\calM^{ss}\}&= 
\#(I(\Q)\backslash (J/C_{J,p}\times G(\A_f^p)/C^p)),\\
\#\{\text{superspecial points of }\calM^{ss}\}&= 
\#(I(\Q)\backslash (J/C_{J,p}'\times G(\A_f^p)/C^p)),\\ 
\#\{\text{connected components of }\calM^{ss}\}&= 
\#(I(\Q)\backslash (J^0\backslash J\times G(\A_f^p)/C^p))\\
&=\#(I(\Q)\backslash (\Z\times G(\A_f^p)/C^p)).
\end{align*}
All numbers are finite.
\end{prop}

\begin{proof}
All numbers are finite as $\calM^{ss}$ is of finite type over $\Spec\bF_p$.

The number of connected components follows from \ref{a1} and
the decomposition
\eqref{ea4} since $\N_0^{red}$ is connected. 

To compute the number of irreducible components and superspecial
points, we have to count these objects in
$\Gamma_j^0\backslash\N_0^{red}$. The number of irreducible components
is equal to the number of orbits of the action of $\Gamma_j^0$ on
$\L_0^{(3)}$. The group $J^0$ acts transitively on $\L_0^{(3)}$
by Lemma~\ref{d18}, hence the number of irreducible components of
$\Gamma_j^0\backslash\N_0^{red}$ is 
equal to $\#
(\Gamma_j^0\backslash J^0/C_{J,p})$. An easy calculation shows the
above formula. 
  
A similar argument proves the claim in case of the superspecial points.
\end{proof}

\begin{rem}
The numbers in Proposition~\ref{a3} can be expressed in terms of class numbers.
\end{rem}



\vskip1cm
\noindent
\textsc{Inken Vollaard\\
Mathematisches Institut\\
Universit\"at Bonn\\
Beringstr. 1\\
53115 Bonn\\
Germany\\
Email Address:} vollaard@math.uni-bonn.de

\end{document}